\documentstyle[12pt,epic,amssymb,amscd,amsthm,eepic]{amsart}
\def\draft{n}

\theoremstyle{plain}

\newtheorem{theorem}{Theorem}
\newtheorem{proposition}{Proposition}[section]
\newtheorem{lemma}[proposition]{Lemma}
\newtheorem{corollary}[proposition]{Corollary}

\theoremstyle{definition}
\newtheorem{definition}[proposition]{Definition}

\newtheorem{question}{Question}

\theoremstyle{remark}

\newtheorem{exercise}[proposition]{Exercise}

\newtheorem{remark}[proposition]{Remark}

\def\printname#1{
	\if\draft y
		\smash{\makebox[0pt]{\hspace{-0.5in}
			\raisebox{8pt}{\tt\tiny #1}}}
	\fi
}

\newlength{\standardunitlength}
\catcode`\@=11
\long\def\@makecaption#1#2{%
    \vskip 10pt
    \setbox\@tempboxa\hbox{%\ifvoid\tinybox\else\box\tinybox\fi
      \small\sf{\bfcaptionfont #1. }\ignorespaces #2}%
    \ifdim \wd\@tempboxa >\captionwidth {%
        \rightskip=\@captionmargin\leftskip=\@captionmargin
        \unhbox\@tempboxa\par}%
      \else
        \hbox to\hsize{\hfil\box\@tempboxa\hfil}%
    \fi}
\font\bfcaptionfont=cmssbx10 scaled \magstephalf
\newdimen\@captionmargin\@captionmargin=2\parindent
\newdimen\captionwidth\captionwidth=\hsize
\catcode`\@=12

\def\lbl#1{\label{#1}\printname{#1}}

\def\BN{\Bbb N}
\def\BZ{{\Bbb Z}}
\def\BQ{{\Bbb Q}}

\def\BC{{\Bbb C}}

\def\K{{\cal{K}}}
\def\L{{\cal{L}}}

\def\T{{\cal{T}}}
\def\M{{\cal{M}}}
\def\F{{\cal{F}}}

\def\N{\cal{N}}

\def\aa{\alpha}

\def\l{\lambda}

\def\as{algebraically split}
     
\def\ihs{integral homology 3-sphere}          

\def\fti{finite type invariant}                                              
\def\uf{unit-framed}                          

\def\ASA{$AS$-admissible}
\def\BA{$B$-admissible}
\def\BLA{$BL$-admissible}

\def\Fb#1{{{{\cal F}}^{\rm b}_{#1}\cal M}}  %%%% filtrations on \cal M.
\def\Fbl#1{{{{\cal F}}^{\rm bl}_{#1}\cal M}}
\def\Fas#1{{{{\cal F}}^{\rm as}_{#1}\cal M}}
\def\FT#1{{{\cal F^{ T}}_{#1}\cal M}}
\def\FHT#1{{{\cal F^{\rm H\cal T}}_{#1}\cal M}}
\def\FK#1{{{\cal F^{ K}}_{#1}\cal M}}
\def\FHK#1{{{\cal F^{\rm H\cal K}}_{#1}\cal M}}
\def\FL#1{{{\cal F^{\cal L}}_{#1}\cal M}}
\def\FHL#1{{{\cal F^{\rm H\cal L}}_{#1}\cal M}}

\def\Gbl#1{{{\cal G^{\rm bl}}_{#1}\cal M}}
\def\Gas#1{{{\cal G^{\rm as}}_{#1}\cal M}}

\def\Tg{{\cal T_{\rm g,1}}}
\def\Kg{{\cal K_{\rm g,1}}}

\def\p{\prime}
\def\sub{\subseteq}

\def\Lg{{\cal L_{\rm g,1}}}
\def\Lgl{{\cal L_{\rm g,1}^{\rm L}}}

\def\blg{\bar{\cal L}_{\rm g,1}}

\def\S{\Sigma}

\def\Z{{\Bbb Z}}
\def\Q{{\Bbb Q}}

\def\lk{{\text{lk}}}
\def\gg{\gamma}

\begin{document}

%%%%%%%%%%%%%%%%%%%%%%\inplude{page1}

\title[Finite type invariants, the mapping class group and blinks]
{Finite type 3-manifold invariants,  the mapping class group and
blinks}

\author{Stavros Garoufalidis}
\address{Department of Mathematics \\
         Harvard university \\
         1 Oxford Street \\
         Cambridge, MA 02318  U.S.A.}
\email{stavros@math.harvard.edu}
\thanks{The  authors were partially supported by NSF grants 
       DMS-95-05105 and DMS-93-03489 respectively. 
This and related preprints can also be obtained 
by
       accessing the WEB in the address 
 {\tt 
http:\linebreak[0]//\linebreak[0]www.\linebreak[0]math.\linebreak[
0]brown.\linebreak[0]edu/\linebreak[0]$\sim$stavrosg/}}
\author{Jerome Levine}
\address{Department of Mathematics\\
        Brandeis University\\
        Waltham, MA 02254-9110}
\email{levine@max.math.brandeis.edu}

\date{This edition: August 18, 1996 \hspace{0.5cm} First edition: 
      March 10, 1996 \\
      Fax numbers: (617) 495 5436 \hspace{0.5cm} (617) 736 3085 \\
  Email: {\tt stavros@math.harvard.edu} \hspace{0.5cm} {\tt 
levine@max.math.brandeis.edu}}

\maketitle

\begin{abstract}
The goal of the present paper is to find higher genus surgery formulas
for the set of \fti s of \ihs s, and to develop a theory of \fti s
which will be applied in a subsequent publication \cite{GL3}
in the study of subgroups of the mapping class group.
The main result is to show that six filtrations 
on the vector space generated by oriented \ihs s
(three coming from surgery on special classes of links and three
coming from subgroups of the mapping class group) 
are equal. En route we introduce the notion of blink (a special
case of a link) and of a new subgroup of the mapping class group.
\end{abstract}

\tableofcontents

%%%%%%%%%%%%%%%%%%%%%%%\inplude{intro}

\section{Introduction}

\subsection{Motivation}
\lbl{sub.moti}

The motivation/goal for the present paper is threefold:
\begin{itemize}
\item
To find higher genus surgery
formulas for the set of \fti s of \ihs s. 
\item
To develop a 
theory of \fti s that has applications in the study of subgroups of
the mapping class group. 
\item
To propose a philosophical explanation of duality in   the recent idea of 
$p$-branes  in string theory.
\end{itemize}

En route to achieving the above goals we came across the notion of a blink
(a special kind of link, see definition \ref{def.blink}) and across a new 
subgroup of the mapping class
group, see section \ref{sub.fil}.

Finite type invariants of \ihs s, though introduced less than a year 
ago
by T. Ohtsuki \cite{Oh},  play a crucial role in understanding 
the quantum invariants of 3-manifolds and may play an important 
role
 in the {\em interaction
between arithmetic, combinatorics, low dimensional topology
and mathematical physics}. They may also play
 an important 
role
in understanding the properties of the Chern-Simons path 
integrals ,
in a way that diverges from the measure-theoretic analytic 
interpretation of
the path integrals. And finally, they may shed some light on the 
obscure and not at all understood arithmetic properties of the 
quantum 
invariants of links and 3-manifolds.
Finite type invariants are defined in terms of decreasing filtrations
of the vector space generated by oriented \ihs s, for a review see 
section
\ref{sub.ftifil}.
 In the present paper, there are two sources for
 such filtrations: one from cutting and pasting along embedded
surfaces in \ihs s, (see section \ref{sub.fil})
 and the other from doing surgery on (framed) links
in \ihs s (see section \ref{sub.fill}).
An equality of such filtrations as shown in corollary \ref{cor.equiv}
answers the first of the goals of the paper, and in fact characterizes
\fti s of \ihs s in terms of their higher genus surgery properties.

The above equality of filtrations can be  used  
\cite{GL3} to study subgroups of  the mapping class group generalizing 
the work of S. Morita \cite{Mo}, \cite{Mo1}.  Due to the length of the
present paper, we need to postpone the above study in a subsequent
publication, see \cite{GL3}.

Finally, a word about the third motivation for the present paper:
$p$-branes were introduced very recently in the  physics
literature, see \cite{W}. Their role in explaining duality phenomena
in string(?) theory and field theory has been exhibited in a number of 
ways.
Since \fti s are related to Chern-Simons field theory (a gauge theory
in three dimensions with  the Chern-Simons function as Lagrangian 
and 
(colored) knots in $3$-manifolds as the observables) we may learn 
something
about duality of gauge theories in three  dimensions by studying
equivalences of \fti s coming from surgery on one-dimensional (links) or 
two-dimensional (surfaces) objects 
in $3$-manifolds. This is a valid thought since \fti s can be thought
of as the partition function of observables of quantum field theories.

\subsection{Finite type invariants and filtrations on $\M$}
\lbl{sub.ftifil}

All 3-manifolds are  oriented and smooth and all
diffeomorphisms are orientation preserving unless otherwise 
mentioned. 
Let $\M$ denote the vector space over $\BQ$ on the set of 
orientation
preserving diffeomorphism classes of \ihs s. Any decreasing 
filtration
$\F_{\ast}$ on $\M$ defines a notion of {\em finite type invariants}
of \ihs s as follows: a map $v: \M \to \BQ$ is called of $\F$-type 
$m$, if
$v(\F_{m+1} \M)=0$.
Examples of such filtrations were originally introduced by T. Ohtsuki
\cite{Oh} and by S. Garoufalidis \cite{Ga}. For a review of
them see section \ref{sub.fill}  below.

In the present paper we
 describe  two main sources for such filtrations: 
 completions of the (group rings) of subgroups of the mapping 
class
group and surgery on special  classes of (framed) links.
The first source is essentially two dimensional, examples of which
will be $\FT {\ast}$ and $\FK {\ast}$ and $\FL {\ast}$
(described in the next section
\ref{sub.fil}).
The second source is one dimensional, examples of which will be
$\Fas {\ast}, \Fb {\ast}, \Fbl {\ast} $, described in section
\ref{sub.fill}.
The equality of all such filtrations may exhibit a duality between
$1$-branes and $2$-branes in 3-dimensional gauge field theory 
which
has not yet been discovered by the physicists.

\subsection{Filtrations on $\M$ from embedded surfaces}
\lbl{sub.fil}
We begin by describing the filtrations on $\M$ that come from 
subgroups
of the mapping class group.
We recall first some well known facts about mapping class groups 
from
the work of D. Johnson \cite{Jo}, \cite{Jo2} and S. Morita \cite{Mo}.
Let $\Gamma_{g,1}$ denote the group of isotopy classes of 
orientation
preserving diffeomorphisms of closed oriented genus $g$ surfaces
$\Sigma_g$ which are the identity on  a disk $D_g 
\subseteq\Sigma_g$.
This group (which we refer to as the {\em mapping class group})
acts on the fundamental group $\pi$ of the open surface $\Sigma_g -
D_g$. Note that $\pi$ is a free group in $2g$ generators. For the sake 
of simplicity in notation, we will suppress the 
dependence of $\pi$ on $g$; we hope that this will not cause any 
confusion.
The above discussion defines a map:
\begin{equation}
\Gamma_{g,1} \to Aut( \pi)
\end{equation}

For a group $G$, let $G_k$ denote the lower central series, 
defined inductively by
$G_1=G$, and $G_{m+1}=[G, G_m]$. Here for $H,K$
subsets of $G$, $[H,K]$ denotes
the subgroup of $G$ generated by $[\eta, \kappa]=\eta^{-1} \kappa^{-
1}
\eta \kappa$, for $\eta \in H$, $ \kappa \in K$. 
%Let $G^k$ denote the
%(nilpotent) quotient $G/G_k$.  
The above action
of the mapping class group induces (for every non-negative integer 
$k$)
an action on $\pi/\pi_{k+1}$: 
\begin{equation}
\Gamma_{g,1} \to Aut( \pi/\pi_{k+1})
\end{equation}

Let $(\Gamma_{g,1})_{[k]}$ denote the kernel of the above map.
It is obvious that $\{ (\Gamma_{g,1})_{[k]} \}_{k \geq 0}$
is a decreasing sequence of normal subgroups of $\Gamma_{g,1}$.
 Note that $(\Gamma_{g,1})_{[k]}$ is denoted by 
$\Gamma_{g,1}(k+1)$
in \cite{Mo}. The reason for shifting the index by one in our present
notation is to make the statements of question 1 (in section 
\ref{sub.que})
easier. 
Much attention has been paid to the first three  members of 
the above sequence. The first, $(\Gamma_{g,1})_{[0]}$
coincides with the mapping class group $\Gamma_{g,1}$ itself.
 The second, $(\Gamma_{g,1})_{[1]}$
is the {\em Torelli group} (i.e., the kernel of the map:
$\Gamma_{g,1} \to Sp(2g, \BZ)$, in other words all diffeomorphisms 
of
the surface that act trivially on the homology), and will from now
on be denoted by $\Tg$. The third, $(\Gamma_{g,1})_{[2]}$ was
 studied extensively in
\cite{Jo}, \cite{Jo2} and \cite{Mo}
and, following their notation, we will denote it by $\Kg$.

In an alternative view,  it was shown by Johnson
(\cite{Jo} and \cite{Jo2}) that $\Tg$ (resp. $\Kg$)
is the subgroup of the mapping
class group generated by Dehn twists on cobounding (resp. bounding)
 simple closed curves. We will find this alternative view very 
useful in the present paper.
%
%and was shown to coincide with the group
% generated by Dehn twists along bounding simple closed curves
%in $\Sigma_g$. Following the notation of \cite{Jo} and \cite{Mo},
%from now on  we will denote this group by $\Kg$.
Note that all the above groups and maps behave well with respect
to an inclusion of a lower genus open surface into a higher genus 
one,
and with respect to the action of the mapping class groups given by
conjugation.

     Consider the lower central series subgroups  
$(\Tg )_n$ or 
$(\Kg )_n$ and their 
''rational closures'' which we denote by $(\Tg )_{(n)}$, $(\Kg )_{(n)}$. 
Here, for a
group $G$ define $G_{(n)}$ to be the normal subgroup consisting of 
all elements
$g$ such that  $g^k\in G_n$ for some $k>0$. 
Recall that, in a nilpotent group, the set of all elements of finite 
order forms a normal subgroup. 
 Thus, for every non-negative integer $n$, we can consider  three 
interesting sequences of 
normal
subgroups of the mapping class group: $(\Gamma_{g,1})_{[n]}$,
 $(\Tg )_{n}$
and  $( \Tg)_{(n)}$.
It was pointed out by Johnson \cite{Jo2}
that $(\Tg )_n \subseteq(\Gamma_{g,1})_{[n]}$. 
In fact, the following, somewhat stronger, inclusion is 
true:
\begin{equation}
\lbl{eq.que1}
(\Tg)_{(n)} \subseteq (\Gamma_{g,1})_{[n]}
\end{equation}
because $(\Gamma_{g,1})_{[n]}$ is the kernel of a homomorphism
from $(\Gamma_{g,1})_{[n-1]}$ into a torsion-free abelian group 
(see
\cite{Jo2}). 
Note that Johnson \cite{Jo2} has shown that:
$(\Tg )_{(2)}=\Kg =(\Gamma_{g,1})_{[2]}$. He asked whether 
$(\Gamma_{g,1})_{[n]}=(\Tg )_{(n)}$ for every $n$, but this was 
answered in the negative by Morita~\cite{Mo} for $n=3$ and by Hain
\cite{Ha2} for $n \geq 3$. 

We now introduce one more subgroup of $\Gamma_{g,1}$  which 
has apparently not been considered in the literature up to now.
 Recall first that
$H_1(\S_g)$ is a symplectic $2g$-dimensional
vector space, the symplectic form being the
intersection form on the homology. Let 
$L\subseteq H_1 (\S_g )$ be any {\em Lagrangian}, i.e., a direct 
summand of rank $g$ on which the intersection pairing vanishes. For 
example, if $M$ is an \ihs , $i:\S_g\subseteq M$ is any embedding 
in $M$ and $M_{\pm}$ is the 
closure of one of the two components of $M -i(\S_g )$, then 
$\ker\{ H_1 (\S_g )\to H_1 (M_{\pm})\}$ is a Lagrangian. Furthermore any 
two Lagrangians are conjugate via an isometry of $H_1 (\S_g )$ and 
so any Lagrangian arises from an embedding this way. 
If $i:\S_g\subseteq M$ is an embedding as above, then we set
$L^i_{\pm}=\ker\{ H_1 (\S_g )\to H_1 (M_{\pm})\}$, where
 $M_+$ is the closure of the {\em positive } 
component of $M-i(\S_g )$, i.e., the component into which the positive
normal vector to $i(\S_g )$ points, and the other component
 is $M_-$. For a fixed 
Lagrangian $L$ define $\Lgl\subseteq\Gamma_{g,1}$ to be the 
subgroup generated by Dehn twists on simple closed curves 
representing elements of $L$; we will call these {\em $L$-twists}. 
Note that if $h\in\Gamma_{g,1}$ then $\Lg^{h_\ast (L)}
=h^{-1}\Lgl h$. Thus 
$\Lgl$ depends upon the choice of $L$ but any two choices give 
conjugate subgroups.  Moreover $\Kg\subseteq\Lgl$ for any choice of 
$L$ and the intersection of all the conjugates of $\Lgl$ is contained 
in $\Tg$ since every element of $H_1 (\S_g )$ belongs to some 
Lagrangian. We will often just use
the notation $\Lg$ for $\Lgl$ when no confusion will arise. 
For an embedding $i:\S_g\subseteq M$ we will use the notation
$\Lg^i =\Lg^{L^i_+}$

 See the Appendix  for more remarks on $\Lg$.

In the present paper we will concentrate on the subgroups 
$\Kg ,\Tg$ and $\Lg$ of the mapping class group.

For a group $G$, let $\BQ G$ denote the rational group algebra of $G$ 
and let $IG$ denote the augmentation ideal in $\BQ G$ (generated
by all elements of the form $g-1$, for all $g \in G$).
Let us now define two decreasing filtrations on $\M$ as follows:
Let $M$ be an \ihs\ and $i:\Sigma \hookrightarrow M$ an embedded,
oriented, connected, separating genus $g$ surface in $M$. Such a 
surface will be
called {\em admissible} in $M$.
Given any element $f$ of the mapping class group of $\Sigma$, let
$M_f$ denote the 3-manifold obtained by cutting $M$ along 
$\Sigma$,
twisting by $f$ and gluing back.  If $f \in \Lg^i$  and $M$ is
an \ihs\ then it is easy to the resulting manifold will also 
be an \ihs . The assignment $f\to M_f$ defines maps 
$\BQ\Lg^{\rm i}\to\M$
and $\BQ\Tg\to\M$. We will be interested in their restrictions 
to the $m^{th}$ power of the augmentation ideals 
(for every non-negative integer $m$):
\begin{equation}
(I \Lg^i )^m \to \M ,\quad 
(I \Tg )^m \to \M \text{\hspace{0.2cm} and  \hspace{0.2cm}}
(I \Kg )^m \to \M
\end{equation}

We now propose the following three filtrations on $\M$:
\begin{definition}
\lbl{def.surf}
Let, respectively, $\FL m$, $\FT m$ and $\FK m$ denote the
 span of the images of the
above maps for all admissible surfaces $\Sigma$ in all \ihs s $M$.
\end{definition}

For the sake of motivation, we
 make a few remarks which will be proved later in sections 
\ref{sub.iadic}
and \ref{sub.heeg}:

\begin{remark}
\lbl{rem.new}
 We will show, in theorem \ref{thm.T2bl}, that the filtration
 $\FK m$ is equivalent to one  
considered by the first author in \cite{Ga}.
The other ones are apparently new.
\end{remark}
We now show how to describe these filtrations using only {\em Heegaard
embeddings}, i.e., embeddings whose complementary components $M_+$ 
and $M_-$ are
handlebodies. For each $g\geq 0$ choose a Heegaard embedding 
$i_g:\S_g\hookrightarrow S^3$ and consider the associated maps:
$(I \Tg)^m \to \M$ and $(I \Kg)^m \to \M$.
\begin{definition}
\lbl{def.h1}
Let $\FHT m$ and $\FHK m$ denote the union of the
 spans, over the chosen $i_g$, of the 
images of
 the maps, as defined above,  $(I \Tg )^m \to \M$  and  
$(I \Kg )^m \to \M$, respectively.
\end{definition}
The filtration $\FL m$ is more complicated to describe. If 
$i:\S\sub M$ is a Heegaard embedding and 
 $L\sub H=H_1 (\S_g )$ is a Lagrangian, we will say that $i, L$
are {\em compatible } if 
$L=(L\cap L_+^i )+(L\cap L_-^i )$, where we recall that 
$L_{\pm}^i=\ker\{ i_{\pm\ast}:H\to H_1 (M_{\pm})\}$ and 
$i_{\pm}:\S_g \hookrightarrow M_{\pm}$ are the inclusions.

Suppose that $h$ is any orientation-preserving diffeomorphism of $\S$.
Then it is easy to see that $L_{\pm}^{ih^{-1}}=h_{\ast}(L^i )$. Thus
$i, L$ are compatible if and only if $ih^{-1}, h_{\ast}(L)$ are 
compatible. 

\begin{proposition}
\lbl{prop.comp}
If\quad  $ i, L$ are compatible and $h\in\Lg^L$ then $M_h$ is an \ihs .
\end{proposition}

\begin{definition}
\lbl{def.h2}
For each genus $g$ choose a Heegaard embedding $i_g :\S_g\sub S^3$.
Let $\FHL m$ denote the union of the span, over all $g$ and all 
$L$ compatible with $i_g$,
 of the 
images of $(I \Lgl )^m \to \M$.
\end{definition}

\begin{proposition}
\lbl{prop.heeg}
The filtrations of $\M$ defined for Heegaard embeddings in definitions
\ref{def.h1} and \ref{def.h2} are the same as those defined for all
admissible embeddings in definition \ref{def.surf}
\end{proposition}
 These propositions will be proved in section \ref{sub.heeg}.

\begin{remark}
\lbl{rem.heg}
 We mention another set of related, and perhaps equal, filtrations.
Consider the lower central series subgroups  
$(\Tg )_m$ or 
$(\Kg )_m$ and their 
''rational closures''  $(\Tg )_{(m)}$, $(\Kg )_{(m)}$. 
We will show in section \ref{sub.iadic} that $\FT m$ 
(respectively $\FK m$) contains the subspace of $\M$ spanned by 
elements of the form $M-M_f$ for all  admissible surfaces
$f: \S \to M$  where 
$f\in (\Tg )_{(m)}$ (respectively $(\Kg )_{(m)}$). But these 
subspaces also define filtrations of $\M$ worth considering.
\end{remark}

\subsection{Filtrations on $\M$ from framed links}
\lbl{sub.fill}

In this section we recall filtrations on $\M$ from special classes
of framed links (\as\ and boundary)
in \ihs s, and we introduce yet another filtration from what we shall 
call {\em blinks}.
We begin by recalling  some definitions from 
\cite{Oh} and \cite{Ga}.

A link $L$ in a 3-manifold $M$ is called {\em algebraically split} 
 if the linking numbers between
its components vanish. A link is called 
{\em boundary} if each component bounds an oriented surface (often
called Seifert), such that these Seifert
surfaces are all disjoint from each other. Of course, the Seifert 
surfaces
are not unique. Let $|L|$ denote the number of components of a link 
$L$.
A {\em framing} of a link $L$ in a \ihs\ $M$ is a choice of 
(isotopy class) of  an essential simple closed curve 
in the boundary of a tubular neighborhood of each of its components. 
Since
  $M$ is an
\ihs\, then a framing $f$ 
on a $r$ component link can be described in terms of a sequence 
$f=(f_1, \cdots, f_r)$, where $f_i \in \BQ \cup \{ 1/0 \}$, with the
convention that $f_i= p_i/q_i$ is the isotopy class of the curve $ 
p_i
(meridian) + q_i (longitude)$. Note that a framing of a link does not
require the choice of an orientation of it.  A {\em unit framing}
$f$ of a link in an \ihs\ is one such that
$f_i  \in \{ -1,1 \}$  for all $i$. A link is called $AS$-{\em 
admissible}
(respectively, $B$-{\em admissible}) if it is \as\ (respectively,
boundary) and \uf . It is obvious that  \BA\ links are \ASA . The 
converse
is obviously false, as the Whitehead link shows.
If $(L,f)$ is a framed link in a 3-manifold $M$, we denote 
by $M_{L,f}$  the 3-manifold
obtained by doing Dehn surgery to each component of the framed link 
$L$.
Let $\M$  denote the set of of \ihs s.
For  an  \ASA\ link $(L,f)$ in a  \ihs\ $M$, let
\begin{equation}
\lbl{eq.def}
[M,L,f] = \sum_{L' \subseteq L} (-1)^{|L'|} M_{L',f'} \in \M
\end{equation}
where the sum is over all sublinks of $L$ (including the empty one), 
$f'$ is the restriction
of the framing $f$ of $L$ to $L'$ and $|L'|$ is the number of 
components
of $L'$. Note that since $(L,f)$ is
an \ASA\ link, $M_{L', f'}$ is a \ihs\ for every sublink $L'$ of $L$.

Let $\Fas m$ (respectively, $\Fb m$) denote the subspace of $\M$ 
spanned by 
$[M,L,f]$ for all \ASA\ (respectively \BA )\ $m$ component links $L$ 
in
\ihs s $M$. It is obvious that $\Fas \ast, \Fb \ast $ are decreasing 
filtrations on the vector space $\M$.
Following Ohtsuki \cite{Oh} and Garoufalidis \cite{Ga} we call  a 
map
$v: \M \to \BQ$ a $AS$-{\em type} (respectively $B$-{\em type}) 
$m$ invariant of
\ihs s if $v(\Fas {m+1})=0$ (respectively, $v(\Fb {m+1})=0$).

We now  present an important definition for the
present paper.
\begin{definition}
\lbl{def.blink}
A {\em blink} $L=L_{bl}$ in an \ihs\ $M$ is a link with the
following properties: 
\begin{itemize}
\item
     The components of $L_{bl}$ are partitioned into classes of 
two 
components each. These classes will be called the {\em pairs}  of 
$L_{bl}$. 
\item 
     $L_{bl}$ is an oriented  link.
\item
     The pairs $\{ p\}$ bound disjoint oriented
surfaces $\{\Sigma_p\}$ in 
$M$ 
(called Seifert surfaces of $L_{bl}$) so that, if $p=(l, l')$, then 
$\partial\Sigma_p =l-l'$ (using the orientations of $l$ 
and $l'$).
\end{itemize}
An example is given in figure \ref{f.unblink}.  An $r$-pair blink
$L_{bl}$ is one such that $|L_{bl}|=2r$.  
A {\em
subblink} $L_{bl}^{'}$ of a blink $L_{bl}$ is a 
sublink $L_{bl}^{'}$ of  $L_{bl}$
which is a union of some of the pairs of $L_{bl}$. Thus an 
$r$-pair blink $L_{bl}$ has $2^{ r }$ subblinks. 

We next discuss admissible framings of blinks.
Recall that every component of a blink is oriented (as a knot).
A {\em zero Seifert-framing} of a blink is the (isotopy class of a) 
parallel
of it in the Seifert surface that the blink bounds. The result is
independent of the Seifert surface chosen, and  depends only on the
orientation of the blink. A zero Seifert-framing (together with 
the orientation of the blink) defines,
for every choice of a pair of integers $(n,m)$, an {\em $(n,m)$ 
Seifert-framing}
of a $1$-pair blink.  A {\em unit} Seifert-framing of a $1$-pair blink
is a $(\epsilon, -\epsilon)$ Seifert-framing, where $\epsilon = \pm 
1$.
A {\em unit Seifert-framing} of a blink is the choice of a
unit Seifert-framing to each of its pairs. A blink is called $BL$-
{\em 
admissible} if it is unit Seifert-framed.

\begin{remark}
\lbl{rem.sei}
Every $r$-component link in an \ihs\ has a {\em zero framing}, 
which, for
a choice of integers $(f_1, \cdots, f_r)$ defines a 
$(f_1, \cdots, f_r)$-framing. For a $1$-pair blink, an $(n,m)$ 
Seifert-framing is equal to a $(n+ l_{12}, m+ l_{12})$ framing of it, 
where
$l_{12}$ is the linking number between the two components of the 
blink.
Note that a blink is not necessarily an \as\ link. The two components
 of a pair may have a non-zero linking number and 
components of different pairs may also have non-zero linking 
number.
\end{remark}

\begin{remark}
\lbl{rem.sursame}
If a \uf\ link $(L,f)$ is included (as a disjoint union of simple closed
curves) in the image of an embedded
surface $i : \S \hookrightarrow M$, 
then $M_{L,f}= M_{i, \tau(L,f)}$ where
$\tau(L,f)$ is an $f$-dependent product of Dehn twists along the simple closed
curves on $\S$ represented by $L$.
\end{remark}

For a \BLA\ $r$-pair blink $L_{bl}$ in an \ihs\ $M$, we denote 
\begin{equation}
\lbl{eq.dblink}
 [ M,L_{bl},f ] = \sum_{L_{bl}^{'}}
(-1)^{ 1/2|L_{bl}^{'}|} M_{L_{bl}^{'},f'} \in \M
\end{equation}
where the sum is over all subblinks $L_{bl}^{'}$ of $L_{bl}$
(including the empty one), 
and $f'$ is the restriction
of the framing $f$ of $L_{bl}$ to $L_{bl}^{'}$.
The above definition makes sense (i.e., each 3-manifold obtained 
by surgery on some pairs of the blink is an \ihs\ \!\!) because of
the following lemma:

\begin{lemma}
\lbl{lemma.sense}
If $(L_{bl},f)$ is a \BLA\ blink in an \ihs\ $M$ then $M_{L_{bl},f}$ 
is a \ihs\   \!\!.   
\end{lemma}

\begin{pf}
Since the order of the first homology of $M_{L_{bl},f}$ is the 
absolute
values of the determinant of the linking matrix of $L_{bl}$,
 we only need to check that the linking matrix of $L_{bl}$ is 
unimodular. We proceed by induction on the number of pairs of
$|L|$. If $L_{bl}= (L_1 ,L_2)$ is a $1$-pair blink, consisting of two 
components
$L_1$ and $L_2$ with unit Seifert-framing,
then the linking matrix of $L_{bl}$ is:
$$\pmatrix
  l_{12} + \epsilon & l_{12} \\
  l_{12} & l_{12} - \epsilon
\endpmatrix
$$
where $l_{12}$ is the linking number between $L_1$ and $L_2$ 
(which does 
{\em not} necessarily vanishes) and $\epsilon = \pm 1$. It is clear 
that 
this is a unimodular matrix.
In general, if $L_{bl}=L_{bl}^{'} \cup (L_1, L_2)$ is an $r$-pair blink
which is the union of an $r-1$ pair blink (with unimodular linking 
matrix
$A$)
and a $1$-pair blink, then the linking
matrix of $L_{bl}$ is:
$$\pmatrix
&&&a_1&a_1\\
&A&&\vdots&\vdots\\
&&&a_{m-1}&a_{m-1}\\
a_1&\hdots &a_{m-1}&k_m + \epsilon &k_m\\
a_1&\hdots &a_{m-1}&k_m &k_m -\epsilon
\endpmatrix
$$
We leave it as an exercise for the reader to show that this matrix
is 
unimodular.
\end{pf}

We denote by $\Fbl m$ the subspace of $\M$ spanned by
$ [ M,L_{bl},f ] $ for all  \BLA\ $m$-pair blinks $L_{bl}$ in
\ihs s $M$.
It is obvious that $\Fbl \ast $ is a  decreasing 
filtration on the vector space $\M$.
We call a map $v: \M \to \BQ$ a $BL$-type $m$ invariant of \ihs s if 
$v( \Fbl {m+1})=0$. 
\end{definition}

\begin{figure}[htpb]
$$\printname{unblink}
	\setlength{\unitlength}{0.03\standardunitlength}
	\begin{array}{c}  \hspace{-1.7mm}
        	\raisebox{-8pt}{\begingroup\makeatletter\ifx\SetFigFont\undefined
% extract first six characters in \fmtname
\def\x#1#2#3#4#5#6#7\relax{\def\x{#1#2#3#4#5#6}}%
\expandafter\x\fmtname xxxxxx\relax \def\y{splain}%
\ifx\x\y   % LaTeX or SliTeX?
\gdef\SetFigFont#1#2#3{%
  \ifnum #1<17\tiny\else \ifnum #1<20\small\else
  \ifnum #1<24\normalsize\else \ifnum #1<29\large\else
  \ifnum #1<34\Large\else \ifnum #1<41\LARGE\else
     \huge\fi\fi\fi\fi\fi\fi
  \csname #3\endcsname}%
\else
\gdef\SetFigFont#1#2#3{\begingroup
  \count@#1\relax \ifnum 25<\count@\count@25\fi
  \def\x{\endgroup\@setsize\SetFigFont{#2pt}}%
  \expandafter\x
    \csname \romannumeral\the\count@ pt\expandafter\endcsname
    \csname @\romannumeral\the\count@ pt\endcsname
  \csname #3\endcsname}%
\fi
\fi\endgroup
\begin{picture}(2367,1627)(0,-10)
\thicklines
\path(242,623)	(189.021,675.054)
	(144.461,721.192)
	(107.662,762.293)
	(77.964,799.235)
	(37.235,864.161)
	(17.000,923.000)

\path(17,923)	(14.825,996.418)
	(20.730,1037.726)
	(30.357,1079.859)
	(43.003,1121.134)
	(57.961,1159.870)
	(92.000,1223.000)

\path(92,1223)	(137.031,1266.762)
	(199.895,1308.301)
	(265.062,1344.690)
	(317.000,1373.000)

\path(317,1373)	(381.395,1409.281)
	(421.429,1430.996)
	(463.921,1453.351)
	(506.798,1475.037)
	(547.986,1494.746)
	(617.000,1523.000)

\path(617,1523)	(661.191,1535.343)
	(715.443,1547.426)
	(776.448,1558.954)
	(840.898,1569.635)
	(905.485,1579.174)
	(966.903,1587.277)
	(1021.844,1593.650)
	(1067.000,1598.000)

\path(1067,1598)	(1104.231,1600.388)
	(1149.712,1602.193)
	(1200.697,1603.367)
	(1254.436,1603.858)
	(1308.184,1603.616)
	(1359.192,1602.593)
	(1404.713,1600.737)
	(1442.000,1598.000)

\path(1442,1598)	(1479.841,1593.831)
	(1525.747,1587.920)
	(1576.976,1580.393)
	(1630.786,1571.379)
	(1684.435,1561.005)
	(1735.182,1549.399)
	(1780.284,1536.688)
	(1817.000,1523.000)

\path(1817,1523)	(1889.612,1480.174)
	(1931.122,1450.418)
	(1973.645,1417.760)
	(2015.337,1384.194)
	(2054.355,1351.716)
	(2117.000,1298.000)

\path(2117,1298)	(2170.783,1253.835)
	(2237.847,1197.312)
	(2300.738,1134.883)
	(2342.000,1073.000)

\path(2342,1073)	(2355.061,1022.077)
	(2359.415,960.500)
	(2355.061,898.923)
	(2342.000,848.000)

\path(2342,848)	(2300.966,785.937)
	(2238.455,723.207)
	(2171.467,666.625)
	(2117.000,623.000)

\path(2117,623)	(2070.660,585.660)
	(2010.826,540.567)
	(1947.830,499.191)
	(1892.000,473.000)

\path(1892,473)	(1865.439,469.925)
	(1817.000,473.000)

\path(242,623)	(300.090,600.359)
	(342.110,581.484)
	(392.000,548.000)

\path(392,548)	(437.379,476.457)
	(455.606,433.866)
	(467.000,398.000)

\path(467,398)	(464.330,323.735)
	(461.645,282.223)
	(467.000,248.000)

\path(467,248)	(502.197,206.998)
	(542.000,173.000)

\path(542,173)	(578.330,134.555)
	(617.000,98.000)

\path(617,98)	(689.143,55.366)
	(731.359,35.420)
	(767.000,23.000)

\path(767,23)	(816.679,14.444)
	(878.394,9.590)
	(940.662,11.441)
	(992.000,23.000)

\path(992,23)	(1023.271,45.564)
	(1067.000,98.000)

\path(1142,248)	(1190.601,309.662)
	(1229.457,352.872)
	(1292.000,398.000)

\path(1292,398)	(1326.312,403.738)
	(1367.971,401.690)
	(1442.000,398.000)

\path(1442,398)	(1492.622,412.228)
	(1554.500,435.500)
	(1616.378,458.772)
	(1667.000,473.000)

\path(1667,473)	(1719.380,476.173)
	(1760.777,475.379)
	(1817.000,473.000)

\path(2192,923)	(2147.792,954.943)
	(2109.497,982.238)
	(2048.011,1024.640)
	(2002.270,1053.722)
	(1967.000,1073.000)

\path(1967,1073)	(1918.240,1092.817)
	(1855.400,1113.916)
	(1792.110,1133.307)
	(1742.000,1148.000)

\path(1742,1148)	(1676.548,1166.810)
	(1636.069,1178.091)
	(1593.282,1189.614)
	(1550.336,1200.622)
	(1509.375,1210.360)
	(1442.000,1223.000)

\path(1442,1223)	(1404.631,1226.401)
	(1359.085,1228.683)
	(1308.092,1229.918)
	(1254.384,1230.181)
	(1200.690,1229.547)
	(1149.740,1228.089)
	(1104.267,1225.882)
	(1067.000,1223.000)

\path(1067,1223)	(1029.190,1218.067)
	(983.361,1210.305)
	(932.227,1200.495)
	(878.503,1189.415)
	(824.903,1177.846)
	(774.143,1166.567)
	(728.937,1156.359)
	(692.000,1148.000)

\path(692,1148)	(625.490,1133.905)
	(584.442,1125.252)
	(541.224,1115.765)
	(498.102,1105.620)
	(457.346,1094.993)
	(392.000,1073.000)

\path(392,1073)	(339.014,1044.231)
	(275.701,1003.141)
	(214.288,959.480)
	(167.000,923.000)

\path(167,923)	(139.959,897.879)
	(92.000,848.000)

\path(17,698)	(13.473,642.036)
	(12.298,600.725)
	(17.000,548.000)

\path(17,548)	(29.419,512.359)
	(49.363,470.143)
	(92.000,398.000)

\path(92,398)	(122.926,362.311)
	(165.114,320.653)
	(208.245,280.168)
	(242.000,248.000)

\path(242,248)	(275.998,208.197)
	(317.000,173.000)

\path(317,173)	(343.914,168.650)
	(392.000,173.000)

\path(542,323)	(599.415,348.729)
	(641.210,368.634)
	(692.000,398.000)

\path(692,398)	(724.059,438.853)
	(767.000,473.000)

\path(767,473)	(809.104,468.279)
	(850.602,447.478)
	(917.000,398.000)

\path(917,398)	(956.266,362.127)
	(992.000,323.000)

\path(992,323)	(1025.090,273.707)
	(1063.100,208.951)
	(1103.060,144.970)
	(1142.000,98.000)

\path(1142,98)	(1213.767,53.754)
	(1256.218,34.816)
	(1292.000,23.000)

\path(1292,23)	(1325.485,18.680)
	(1367.236,17.750)
	(1442.000,23.000)

\path(1442,23)	(1492.622,37.228)
	(1554.500,60.500)
	(1616.378,83.772)
	(1667.000,98.000)

\path(1667,98)	(1716.480,99.011)
	(1778.570,94.344)
	(1841.125,91.505)
	(1892.000,98.000)

\path(1892,98)	(1928.380,114.075)
	(1968.173,138.046)
	(2009.831,167.686)
	(2051.810,200.769)
	(2092.562,235.068)
	(2130.542,268.357)
	(2192.000,323.000)

\path(2192,323)	(2230.921,358.895)
	(2267.000,398.000)

\path(2267,398)	(2309.637,470.143)
	(2329.581,512.359)
	(2342.000,548.000)

\path(2342,548)	(2329.666,620.810)
	(2325.773,662.140)
	(2342.000,698.000)

\path(2342,698)	(2342.000,698.000)

\end{picture} }
        	\hspace{-1.9mm}
	\end{array}
 $$
\caption{An unblink bounding a genus $0$ surface.}\lbl{f.unblink}
\end{figure}
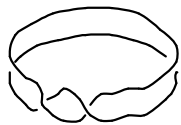

\begin{remark}
\lbl{rem.moti}
The motivation and usefulness
 of the above notion of blink comes from several
facts.
\begin{itemize}
\item
     Blinks are closely related to bounding pairs of simple closed 
curves
in an embedded oriented surface in the 3-manifold; see theorem 
\ref{thm.T2bl}. We will show that the filtration $\Fbl {\ast}$ on 
$\M$
coming from blinks is equal
to the filtration $\FT {\ast}$ coming from the $I$-adic completion 
of
the Torelli group, in much the same way that (see 
remark \ref{rem.new}) the filtration 
$\Fb {\ast}$ coming from boundary links is equal to the filtration 
$\FK {\ast}$ ,
see theorem \ref{thm.T2bl}.
\item
     Johnson \cite{Jo} proved that Dehn twists on bounding pairs 
of simple closed curves generate the 
Torelli group and  that Dehn twists on bounding closed curves 
generate 
$\Kg$. 
\end{itemize}
\end{remark}

\begin{remark}
\lbl{rem.bbl}
For later reference, let us point out that $\Fb m\subseteq \Fbl 
m$. 
Indeed, given any framed 
boundary link we can convert it to a blink by punching a small hole in 
each Seifert surface and putting the appropriate unit framing on the 
boundary of the hole. Framed surgery on this blink obviously gives 
the same result as surgery on the original boundary link.
\end{remark}

\begin{remark}
We claim  that $\Fbl m$
can be generated by $[M, L_{bl}, f]$ for \BLA\ $m$-pair blinks
such that each pair bounds a genus $1$ surface. Indeed, this follows
from the following identity (and induction on the genus)
where $L_1, L_2, K_1, K_2$
are as in  figure \ref{1gblink}:
\begin{equation}
[M, L_{bl},f]=M - M_{L_1, L_2}=[M, L_{bl}^1, f^1] +
[M_{L_{bl}^1}, L_{bl}^2, f^2]
\end{equation}
where $L_{bl}=(L_1,L_2)$, $L_{bl}^1=(L_1,K_1)$ and 
$L_{bl}^2=(K_2,L_2)$
are unit Seifert-framed $1$-pair blinks.
\end{remark}

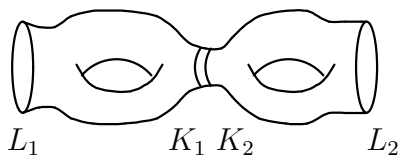
\begin{figure}[htpb]
$$\printname{1gblink}
	\setlength{\unitlength}{0.03\standardunitlength}
	\begin{array}{c}  \hspace{-1.7mm}
        	\raisebox{-8pt}{\begingroup\makeatletter\ifx\SetFigFont\undefined
% extract first six characters in \fmtname
\def\x#1#2#3#4#5#6#7\relax{\def\x{#1#2#3#4#5#6}}%
\expandafter\x\fmtname xxxxxx\relax \def\y{splain}%
\ifx\x\y   % LaTeX or SliTeX?
\gdef\SetFigFont#1#2#3{%
  \ifnum #1<17\tiny\else \ifnum #1<20\small\else
  \ifnum #1<24\normalsize\else \ifnum #1<29\large\else
  \ifnum #1<34\Large\else \ifnum #1<41\LARGE\else
     \huge\fi\fi\fi\fi\fi\fi
  \csname #3\endcsname}%
\else
\gdef\SetFigFont#1#2#3{\begingroup
  \count@#1\relax \ifnum 25<\count@\count@25\fi
  \def\x{\endgroup\@setsize\SetFigFont{#2pt}}%
  \expandafter\x
    \csname \romannumeral\the\count@ pt\expandafter\endcsname
    \csname @\romannumeral\the\count@ pt\endcsname
  \csname #3\endcsname}%
\fi
\fi\endgroup
\begin{picture}(5501,2022)(0,-10)
\thicklines
\put(1537.500,676.500){\arc{1277.204}{3.8438}{5.5809}}
\put(3937.500,676.500){\arc{1277.204}{3.8438}{5.5809}}
\put(3112.500,1201.500){\arc{977.880}{2.5749}{3.7083}}
\put(3262.500,1201.500){\arc{977.880}{2.5749}{3.7083}}
\put(225,1202){\ellipse{300}{1274}}
\put(5025,1202){\ellipse{300}{1274}}
\path(975,1239)	(1011.704,1174.189)
	(1044.639,1118.938)
	(1074.463,1072.370)
	(1101.836,1033.605)
	(1151.865,975.969)
	(1200.000,939.000)

\path(1200,939)	(1252.462,913.975)
	(1315.335,893.519)
	(1385.295,877.583)
	(1459.016,866.119)
	(1533.173,859.077)
	(1604.439,856.410)
	(1669.490,858.067)
	(1725.000,864.000)

\path(1725,864)	(1798.483,885.480)
	(1840.382,903.016)
	(1883.171,923.595)
	(1924.939,946.099)
	(1963.774,969.413)
	(2025.000,1014.000)

\path(2025,1014)	(2085.892,1088.006)
	(2125.162,1151.243)
	(2148.540,1191.727)
	(2175.000,1239.000)

\path(3375,1239)	(3411.704,1174.189)
	(3444.639,1118.938)
	(3474.463,1072.370)
	(3501.836,1033.605)
	(3551.865,975.969)
	(3600.000,939.000)

\path(3600,939)	(3652.462,913.975)
	(3715.335,893.519)
	(3785.295,877.583)
	(3859.016,866.119)
	(3933.173,859.077)
	(4004.439,856.410)
	(4069.490,858.067)
	(4125.000,864.000)

\path(4125,864)	(4198.483,885.480)
	(4240.382,903.016)
	(4283.171,923.595)
	(4324.939,946.099)
	(4363.774,969.413)
	(4425.000,1014.000)

\path(4425,1014)	(4485.892,1088.006)
	(4525.162,1151.243)
	(4548.540,1191.727)
	(4575.000,1239.000)

\path(225,1839)	(269.583,1821.823)
	(308.145,1807.314)
	(369.844,1785.424)
	(415.371,1771.572)
	(450.000,1764.000)

\path(450,1764)	(499.552,1758.122)
	(561.761,1754.678)
	(624.340,1755.894)
	(675.000,1764.000)

\path(675,1764)	(729.405,1792.120)
	(790.028,1836.713)
	(849.386,1882.448)
	(900.000,1914.000)

\path(900,1914)	(965.009,1936.936)
	(1005.470,1948.389)
	(1048.324,1959.281)
	(1091.387,1969.201)
	(1132.477,1977.735)
	(1200.000,1989.000)

\path(1200,1989)	(1244.853,1993.010)
	(1299.529,1995.577)
	(1360.745,1996.850)
	(1425.217,1996.976)
	(1489.663,1996.106)
	(1550.797,1994.387)
	(1605.338,1991.969)
	(1650.000,1989.000)

\path(1650,1989)	(1695.142,1984.948)
	(1750.048,1979.299)
	(1811.423,1972.091)
	(1875.971,1963.361)
	(1940.399,1953.149)
	(2001.409,1941.493)
	(2055.708,1928.430)
	(2100.000,1914.000)

\path(2100,1914)	(2153.847,1886.048)
	(2217.064,1844.865)
	(2277.999,1800.749)
	(2325.000,1764.000)

\path(2325,1764)	(2374.853,1714.932)
	(2434.151,1649.276)
	(2495.124,1584.733)
	(2550.000,1539.000)

\path(2550,1539)	(2598.066,1515.778)
	(2660.227,1493.505)
	(2723.526,1475.229)
	(2775.000,1464.000)

\path(2775,1464)	(2824.282,1457.161)
	(2886.041,1452.116)
	(2948.530,1453.013)
	(3000.000,1464.000)

\path(3000,1464)	(3040.928,1491.867)
	(3082.005,1534.395)
	(3119.579,1579.225)
	(3150.000,1614.000)

\path(3150,1614)	(3182.536,1647.494)
	(3223.635,1690.189)
	(3265.417,1732.290)
	(3300.000,1764.000)

\path(3300,1764)	(3347.119,1800.638)
	(3408.251,1844.569)
	(3471.508,1885.715)
	(3525.000,1914.000)

\path(3525,1914)	(3561.684,1926.878)
	(3606.772,1939.207)
	(3657.516,1950.763)
	(3711.173,1961.317)
	(3764.995,1970.645)
	(3816.237,1978.519)
	(3862.154,1984.713)
	(3900.000,1989.000)

\path(3900,1989)	(3965.947,1994.224)
	(4006.490,1996.401)
	(4049.250,1997.824)
	(4092.104,1998.143)
	(4132.928,1997.011)
	(4200.000,1989.000)

\path(4200,1989)	(4269.142,1960.306)
	(4308.838,1938.254)
	(4350.000,1914.000)
	(4391.162,1889.746)
	(4430.858,1867.694)
	(4500.000,1839.000)

\path(4500,1839)	(4550.354,1833.710)
	(4613.194,1834.298)
	(4675.687,1837.237)
	(4725.000,1839.000)

\path(4725,1839)	(4774.573,1839.000)
	(4837.500,1839.000)
	(4900.427,1839.000)
	(4950.000,1839.000)

\path(4950,1839)	(4976.014,1839.000)
	(5025.000,1839.000)

\path(225,564)	(269.583,581.177)
	(308.145,595.686)
	(369.844,617.576)
	(415.371,631.428)
	(450.000,639.000)

\path(450,639)	(499.552,644.878)
	(561.761,648.322)
	(624.340,647.106)
	(675.000,639.000)

\path(675,639)	(729.405,610.880)
	(790.028,566.288)
	(849.386,520.552)
	(900.000,489.000)

\path(900,489)	(965.009,466.064)
	(1005.470,454.611)
	(1048.324,443.719)
	(1091.387,433.799)
	(1132.477,425.265)
	(1200.000,414.000)

\path(1200,414)	(1244.853,409.990)
	(1299.529,407.423)
	(1360.745,406.150)
	(1425.217,406.024)
	(1489.663,406.894)
	(1550.797,408.613)
	(1605.338,411.031)
	(1650.000,414.000)

\path(1650,414)	(1695.142,418.052)
	(1750.048,423.701)
	(1811.423,430.909)
	(1875.971,439.639)
	(1940.399,449.851)
	(2001.409,461.507)
	(2055.708,474.570)
	(2100.000,489.000)

\path(2100,489)	(2153.847,516.952)
	(2217.064,558.135)
	(2277.999,602.251)
	(2325.000,639.000)

\path(2325,639)	(2374.853,688.068)
	(2434.151,753.724)
	(2495.124,818.267)
	(2550.000,864.000)

\path(2550,864)	(2598.066,887.222)
	(2660.227,909.495)
	(2723.526,927.771)
	(2775.000,939.000)

\path(2775,939)	(2824.282,945.839)
	(2886.041,950.884)
	(2948.530,949.987)
	(3000.000,939.000)

\path(3000,939)	(3040.928,911.133)
	(3082.005,868.605)
	(3119.579,823.775)
	(3150.000,789.000)

\path(3150,789)	(3182.536,755.506)
	(3223.635,712.811)
	(3265.417,670.710)
	(3300.000,639.000)

\path(3300,639)	(3347.119,602.362)
	(3408.251,558.431)
	(3471.508,517.285)
	(3525.000,489.000)

\path(3525,489)	(3561.684,476.122)
	(3606.772,463.793)
	(3657.516,452.237)
	(3711.173,441.683)
	(3764.995,432.355)
	(3816.237,424.481)
	(3862.154,418.287)
	(3900.000,414.000)

\path(3900,414)	(3965.947,408.776)
	(4006.490,406.599)
	(4049.250,405.176)
	(4092.104,404.857)
	(4132.928,405.989)
	(4200.000,414.000)

\path(4200,414)	(4269.142,442.694)
	(4308.838,464.746)
	(4350.000,489.000)
	(4391.162,513.254)
	(4430.858,535.306)
	(4500.000,564.000)

\path(4500,564)	(4550.354,569.290)
	(4613.194,568.702)
	(4675.687,565.763)
	(4725.000,564.000)

\path(4725,564)	(4774.573,564.000)
	(4837.500,564.000)
	(4900.427,564.000)
	(4950.000,564.000)

\path(4950,564)	(4976.014,564.000)
	(5025.000,564.000)

\put(0,39){\makebox(0,0)[lb]{$L_1$}}
\put(5025,39){\makebox(0,0)[lb]{$L_2$}}
\put(2250,39){\makebox(0,0)[lb]{$K_1$}}
\put(2925,39){\makebox(0,0)[lb]{$K_2$}}
\end{picture} }
        	\hspace{-1.9mm}
	\end{array}
 $$
\caption{Shown here ia a $1$-pair blink $(L_1, L_2)$ of genus $2$, 
and
two (parallel) knots $K_1, K_2$ on the Seifert surface that it 
bounds. 
The knots separate the surface in two genus $1$ surfaces, and are
oriented and framed in an opposite way.}\lbl{1gblink}
\end{figure}

We will  find later on  the following notation  useful.
\begin{definition}
\lbl{def.hat} 
If $\N$ is any subspace of $\M$, then $\widehat{\N}=\cap_n (\N +\Fas 
n)$.
\end{definition}

At this point we have introduced six filtrations on $\M$:
$$\FT m ,\ \FK m ,\ \FL m,\ \Fas m ,\ \Fb m ,\ \Fbl m $$ 
The purpose of the paper is 
to show, among other things, that the six associated filtrations 
$$\FT m ,\ \widehat{\FK {m}} ,\ \FL m,\ \Fas m ,\ \widehat{\Fb 
{m}},\ \Fbl 
m $$ 
are actually
{\em equal} (after renumbering); see corollary \ref{cor.equiv}.

\subsection{Statement of the results}
\lbl{sub.results}

We are now ready to formulate 
the main results of the paper which consist
of three parts. The division in three parts is for the convenience of 
the 
reader, since the methods used in each part 
are very different.

 In the first part, we compare mapping class group filtrations and
link filtrations as follows:

\begin{theorem}
\lbl{thm.T2bl}
For every non-negative integer  $m$ we have:
\begin{itemize}
\item[(a)]
     $\FT m = \Fbl m $
\item[(b)]
     $\FK m = \Fb m $
\item[(c)]
     $\FL m = \Fas m $
\end{itemize}
as subspaces of $\M$.
\end{theorem}

A part of the argument used in the proof of theorem 
\ref{thm.T2bl} will also yield the following interesting fact.

\begin{proposition}
\lbl{thm.boundary}
Every \ihs\ can be obtained by surgery on a boundary link in $S^3$.
\end{proposition}

The above proposition (which, as  the referee informs us, was
already  known to
C. Lescop using a different argument) implies the following corollary:

\begin{corollary}
\lbl{cor.conversation}
The Casson invariant $\l_C$ of an \ihs\ (which we may assume is 
diffeomorphic to $S^3_{L,f}$ for a unit-framed boundary 
$r$-component link $L$) is given by
\begin{equation}
 \l_C(M)= \sum_{i=1}^{r} f_i \phi(L_i)
\end{equation}
where $\{ L_i \}$ are the components of the link and $\phi(L_i)$
is the second derivative of the (normalized) Alexander polynomial
of the knot $L_i$, see \cite{Ho}. The point, of course, is that
the Casson invariant of an \ihs\ 
can, therefore, be calculated in terms of the associated knot 
invariant (i.e.,
the second derivative of the normalized Alexander polynomial).
\end{corollary}
 
A generalization of the above corollary appears in \cite{Ga2}.

In the second part we compare the three filtrations $\Fas {\ast}$,
$\Fb {\ast}$ and $\Fbl {\ast}$ coming from special classes of links
in \ihs s.
We have the following results:

\begin{theorem}
\lbl{thm.bl2as}
For every non-negative integer $m$ we have that:
\begin{equation}
\Fbl {2m-1} \subseteq \Fas {3m}
\end{equation}
as subspaces of $\M$.
\end{theorem}

\begin{corollary}
\lbl{cor.invariant}
 Let $v$ be an $AS$-type $3m$ invariant of \ihs s,  $M$ an
\ihs\ and $\Sigma$ an oriented embedded genus $g$ surface. Let  
$f \in (\Tg)_{(2m+1)}$ be an element of the rational closure of 
$ (\Tg)_{2m+1}$ in the Torelli
group, as in remark \ref{rem.heg}(b). Then, using that remark, we 
have that
$v(M)=v(M_f)$.
\end{corollary}
 
\begin{theorem}
\lbl{thm.as2b}
     There is an increasing function $f : \BN \to \BN$ such that 
for every non-negative integer $m$ we have that:
\begin{equation}
\Fas {f(m)} \subseteq \Fbl m
\end{equation}

\end{theorem}

\begin{remark}
\lbl{rem.f}
In fact, since $\Fas {\ast}$ and $\Fbl {\ast}$ are decreasing
filtrations, the proof of theorem \ref{thm.as2b}
shows that we can take $f(m)= c m^{13}$, for some positive integer 
constant $c$.
\end{remark}

\begin{remark}
\lbl{rem.4T}
En route to proving theorem \ref{thm.as2b} we give, in proposition
\ref{prop.tubing}, a $4$-term relation that holds in $\M$. The 
relation
is apparently new. The relation between the $4$-term and existing 
relations on $\M$ (namely the $AS$ and the $IHX$, see \cite{GO2})
is addressed in question \ref{que.relation}.
\end{remark}

\begin{theorem}
\lbl{thm.as2bl}
\begin{itemize}
\item
     With the notation of definition \ref{def.hat}, 
for every $m$ we have:
\begin{equation}
\Fas {3m} \subseteq \widehat{\Fbl {2m}}
\end{equation}
\item
     Together with theorem \ref{thm.bl2as} and theorem 
\ref{thm.as2b}
 this implies that,
for every $m$:
\begin{equation}
\Fbl {2m} = \Fbl {2m-1} = \Fas {3m}
\end{equation}
\end{itemize}
\end{theorem}

Combining theorem \ref{thm.as2bl} and theorem $2$ from [GL], 
we obtain the
following corollary:

\begin{corollary}
\lbl{cor.equiv}
For every non-negative integer $m$, the six filtrations on
$\M$  shown below are equal:

$$ 
\vbox{
\offinterlineskip
\halign{\strut#&&\vrule#&\quad\hfil#\hfil\quad\cr
\noalign{\hrule}
&& $\widehat{\FK {m}}$  && $\FT {2m}$ && $\FL {3m}$  &\cr
\noalign{\hrule}
&& $\widehat{\Fb {m}}$ && $\Fbl {2m}$  && $\Fas {3m}$  &\cr
\noalign{\hrule}
}}
$$ 
Note that in the top row the three filtrations on $\M$
 come from subgroups of the
mapping class group, and on the bottom row are their equivalent 
filtrations
that come from special classes of links (boundary, blinks and 
algebraically
split).
\end{corollary}

In the third part of the paper we discuss the relation between 
{\em blink surgery equivalence} and the {\em  Seifert 
matrix}
of a blink.
In \cite{Oh} and \cite{GL1}
 it is observed that if two \ASA\ $n$-component
links are surgery equivalent then the associated elements in 
$\Gas n $ are equal. Since it is known \cite{Le} that surgery 
equivalence is determined by the triple $\bar\mu$-invariants of 
Milnor, these numerical link invariants provide a set of generators
of $\Gas n$. We now present an analogous result 
for blinks and $\Gbl n$. This gives, as a consequence, a set of 
generators
for $\Gbl n$. Although $\Gbl {2n} = \Gas {3n}$, the
above mentioned two sets of generators seem to be {\em different}.
We hope to explore this point further some time in the near future.

\begin{definition}
Let $(M,L)$ be a blink. An {\em elementary (blink) surgery 
equivalence}
on $(M,L)$ is a surgery on a unit Seifert-framed blink 
$L^{\p\p}\subseteq M-L$ such that $L$ and $L^{\p\p}$ have 
Seifert surfaces which are disjoint from each other. 
(Thus $L\cup L^{\p\p}$ is a blink in $M$). We say 
$(M_{L^{\p\p}},L)$ is surgery equivalent to $(M,L)$.
More generally, surgery equivalence is the equivalence relation 
among blinks
generated by elementary surgery equivalence. 

\end{definition}

\begin{theorem} 
\lbl{thm.seifert}
Two blinks are surgery equivalent if and only 
if they admit equal    Seifert matrices.
\end{theorem}

\subsection{Plan of the proof}

As mentioned in the abstract and in the introduction, for the
convenience of the reader,  the proofs
of the results appear in three sections.

 In section~\ref{sec.fs} 
 we review $I$-adic and nilpotent completions
and prove the claims in remark \ref{rem.heg}, as well as theorem
\ref{thm.T2bl} and propositions \ref{prop.comp}, \ref{prop.heeg}
and  \ref{thm.boundary}.
The proofs in this section are mostly algebraic manipulations in
group algebras and a bit of cutting, pasting and tubing arguments.

In section \ref{sub.review} we review the main identities
in $\M$ in graphical and algebraic form, as well as a few facts about
blinks. There is a plethora of identities on $\M$, coming from
Kirby calculus, that is, from different ways of representing an \ihs\
by surgery on a link. Most of the identities are known, however
in sections \ref{sub.4T} and we \ref{sub.4Tmore} we introduce a new 
one 
(the $4$-term relation) that will be used 
crucially
in the present paper.

The identities of section \ref{sub.review} together with
induction  are the main tools used in the proofs of theorems 
\ref{thm.bl2as}, \ref{thm.as2b} and \ref{thm.as2bl}. We warn the 
reader
that though the statements of the above mentioned theorems seem
similar, the methods used to prove them are {\em different}. Each
theorem requires its own use of the same identities.  

In section~\ref{sec.smatrix} we discuss the notion of blink 
surgery equivalence
and prove theorem \ref{thm.seifert}.

We collect in the  appendix  some results related to the subgroup 
$\Lg$
of the 
mapping class group. These results  
are not directly  used in the present paper, but they may help  clarify
the structure of $\Lg$.

%Finally, as the reader may have observed from past work, we give
%a section of questions  (\ref{sub.que} ) and
%a section of philosophical comments (\ref{sub.philo} ) which include
%  questions  of philosophical as well as mathematical interest.
%We hope to return to them at some point in the future.

\subsection{Questions}
\lbl{sub.que}
In this section we propose a few questions that may lead to a 
better
understanding of the subject. They are addressed to different 
audiences, at least the way they are stated.

\begin{question}
\lbl{que.def} Are the subspaces of $\M$ (defined in 
Remark \ref{rem.new}) spanned by elements of the form $M-M_f$, 
for 
$(M,\Sigma )$ admissible and  $f\in  (\Lg )_{(m)}, (\Tg )_{(m)}$ 
or  
$(\Kg )_{(m)}$, respectively,  the same as $\FL{m}, \FT{m}$  or  
$\FK{m}$, respectively?
\end{question}

\begin{question}
\lbl{que.relation}
In sections \ref{sub.4T} and \ref{sub.4Tmore}
we show (three versions of) a $4$-term relation that holds on
$\M$. In the first version this
 relation really comes from the $4$-term relation on the 
space
of knots (in $S^3$) via the Dehn surgery map. On the other hand, one
knows that an antisymmetry ($AS$) and an $IHX$ relation hold on 
$\M$, see
\cite{GO2}. Is it true that the $AS$ and the $IHX$ relation are 
equivalent
to the $4$-term relation on $\M ?$ Note that the $AS$ and the $IHX$
relation are equivalent to the $4$-term relation on the space of 
knots,
see \cite{B-N1}.
\end{question}

\begin{question}
\lbl{que.el}  
The subgroup $\Lg$ of the mapping class group which is introduced 
in section~\ref{sub.fil} is related to a larger subgroup $\blg$ 
consisting of all $h\in\Gamma_{g,1}$ such that 
$h|L=\text{identity}$. In the appendix we discuss these two 
subgroups and, for example, show that $\blg\not=\Lg$. But is the 
filtration of $\M$ defined by the powers $(I\blg )^m$ different from 
that defined by $(I\Lg )^m$? 

In the appendix we also show that $(\blg )_5\sub\Kg$ but $(\blg 
)_4\not\sub\Kg$. Is $(\Lg )_4\sub\Kg$? 
\end{question}

\subsection{Philosophical comments}
\lbl{sub.philo}

In this section we propose a few questions of philosophical interest
 that may lead to a better 
understanding of the subject. They are not directly related to the 
results
of the present paper, and the questions themselves are somewhat 
vague.
Positive answers may nevertheless bring presents to a variety of 
areas.

\begin{question}
Interpret the results of the present paper in terms
of the not-yet-discovered duality of three dimensional gauge 
theories.
\end{question}

\begin{question}
Give a Hodge structure on the complexified space \mbox{$\M \otimes_{\BQ} 
\BC$}
in a way compatible with the filtration discussed in the present 
paper.
\end{question}

\subsection{Acknowledgment}
We wish to thank the {\tt Internet} for providing the required 
communications. The first author wishes to thank D. Bar-Natan for
intriguing conversations that lead to some of the results of the
present paper. 

\section{Equivalence of filtrations from links and from surfaces}
\lbl{sec.fs}

\subsection{$I$-adic and nilpotent completions}
\lbl{sub.iadic}

In this section we recall some facts about the lower central series 
of a 
group and its relation to the $I$-adic filtration of the group algebra. 
Let $G$ be a group and $\Q G$ the group $\Q$-algebra, i.e., the 
vector space over $\Q$ with $G$ as basis. Multiplication is defined 
by linearly extending the multiplication of $G$. The augmentation
ideal $IG\subseteq\Q G$ is the two-sided ideal generated by all 
elements of the form 
$g-1, g\in G$. The {\em $I$-adic filtration} of $\Q G$ is the 
sequence of powers $(IG)^n$ of $IG$. It is not difficult to show that 
if $g\in G_n$ the $n$-th lower central series subgroup of $G$
(see section \ref{sub.fil}) then $g-1\in (IG)^n$. In fact 
a bit more is true, namely let $G_{(n)}$ denote the  {\em rational
closure} of the group $G_n$ defined in remark \ref{rem.heg}. 
Then $g-1\in 
(IG)^n$ whenever $g\in G_{(n)}$. This follows from the fact that 
$g^m -1\in (IG)^n$ for some $m$, and from the 
formula: $g^m -1=\sum_{i=1}^m\binom mi (g-1)^i$.
The point of considering the rational closures is the 
theorem of Jennings (see \cite{Qu}) which says that the 
converse is true: if $g-1\in (IG)^n$, then $g\in G_{(n)}$. This 
proves Remark \ref{rem.heg} (b).

\subsection{Proof of proposition \ref{prop.comp}}
\lbl{sub.comp}

Let $h$ be as in the statement of proposition \ref{prop.comp}.
A Mayer-Vietoris argument shows that 
$$H_1 (M_h )\cong H_1 (\S )/L_- +h_{\ast}(L_+ ) $$
Thus it suffices to show that $h_{\ast}(L_+ )=L_+ \mod L_-$. First
we show that $L$ has a complementary Lagrangian $L'$ which is
compatible with $i$. Assuming this then, since $h$ is symplectic,
$h_{\ast}|L'$ is the identity $\mod L$. Thus if $\aa\in L_+ \cap
L'$ then $h_{\ast}(\aa)-\aa\in L\cap L_+ \mod L_-$. So, with
respect to the direct sum decomposition $L_+ =(L\cap L_+ )\oplus
(L'\cap L_+ )$, $h_{\ast}|L_+$ has the form $\left(\smallmatrix 
I&X\\0&I\endsmallmatrix\right)$.

To construct the complementary Lagrangian $L'$ first choose
any complementary summand $L'_+$ of $L\cap L_+$ in $L_+$.
Then let $L'_- \sub L_-$ be the annihilator, under the intersection
pairing, of $L'_+$. Since $L\cap L_-$ must pair non-singularly
with $L'_+$, then $L'_-$ must be complementary to $L\cap L_-$
in $L_-$. It is clear that $L'_+ +L'_-$ is the desired complementary
Lagrangian.

This completes the proof of proposition \ref{prop.comp}.

\subsection{Proof of proposition \ref{prop.heeg}}
\lbl{sub.heeg}

We will first prove $$\FHK m=\FK m,\ \FHT m=\FT m\ {\rm and }\ \FL m\sub\FHL m$$
 Suppose that 
$i:\S_g\subseteq M$ is an admissible surface in $M$ and $f\in\Lg^i$. 
Then 
$i(\S_g )$ separates $M$ into the two components $M_+ ,M_-$. Using 
handle decompositions of $M_+$ and $M_-$ we can find two handlebodies
 $H_+\sub M_+ , H_-\sub M_-$ such that the connected sum 
$\S_{g'} =i(\S_g )\ \sharp\ \partial H_+\ \sharp\ \partial H_-$ is a 
Heegaard embedding.Note that the complementary components of $\S_{g'}$ 
are $M_{\pm}' =(M_{\pm} -H_{\pm})\ \sharp\ H_{\mp}$, using 
boundary connected sum. Thus we  
see that $M_f =M_{f'}$, where 
$f'\in\Gamma_{g',1}$ is the image of $f$ under the canonical inclusion 
$\Lg\sub\Gamma_{g,1}\sub\Gamma_{g',1}$.  Note 
that if $f\in\Tg$ then $f'\in\T_{g',1}$ and if
$f\in\Kg$ then $f'\in\K_{g',1}$. To complete the
proof that $\FL m\sub\FHL m$ we will need to show that there is a compatible
 Lagrangian 
$L\sub H_1 (\S_{g'})$ so that 
$L\supseteq\ker\{ i_+:H_1 (\S_g )\to H_1 (M_+ )\}$.
(Note that $\ker\{ i_+:H_1 (\S_g )\to H_1 (M_+ )\}=L_+\nsubseteq L_+ '=\ker
\{ i_+:H_1 (\S_{g'} )\to H_1 ( M_+ ')\}$.) We show this below. If $g'$ is 
large enough it follows 
from \cite{Mo} that  $M=S^3_h$ for some 
$h\in\K_{g',1}$ using our chosen Heegaard surface $i_{g'}(\S_{g'})$. 
Thus $M_f =S^3_{f'h}$.

Now  a generator of $\FL m$ is, by definition, a linear 
combination 
$\sum a_i M_{f_i}$, where $\sum a_i f_i\in (I\Lg^i )^m$ and all terms 
are defined with respect to the same embedding in $M$. But 
the discussion above shows that we can rewrite this as $\sum a_i 
S^3_{f'_i h}$. Note that  $\sum a_i f'_i \in\Q\L_{g'}^L$ is the 
image of $\sum a_i f_i\in (I\Lg^i )^m$ and so belongs to 
$(I\L_{g'}^L)^m$. 
Therefore $\sum a_i f'_i h\in (I\Lgl )^m$ since $h\in\Kg$. 
The same argument works 
for $\Tg$ and $\Kg$. Note that we have completed the proof of 
proposition \ref{prop.heeg} for $\FK m$ and $\FT m$. 

To complete the proof
that  $\FL m\sub\FHL m$ we need to construct the compatible Lagrangian 
$L\sub H_1 (\S ')$. Let $A_{\pm}\sub H_1 (\S ')$ denote the subgroup 
generated by the boundary circles of the meridian disks of $H_{\pm}$.
Note that $A_{\pm}\sub L_{\mp}'$. We now define $L=L_+ +A_+ +A_-$. To 
show that this is compatible, i.e., $L=(L\cap L_+ ')+(L\cap L_- ')$, we 
first observe that $L_+ \sub L_+ ' +A_+$. This can be seen geometrically
as follows. If $\gg$ is a closed curve in $\S$ representing an element 
of $L_+$, then $i(\gg )$ bounds a (singular) surface in $M_+$. This 
surface will intersect $H_+$ generically as a union of meridian disks. 
Thus $\gg -$the corresponding element of $A_+$ 
bounds a surface in $M_+ '$ and so represents an element of $L_+ '$. 
Now suppose $\aa\in L$. Then we write $\aa =l_+ +h_+ +h_-$, where 
$l_+ \in L_+ ,h_+ \in A_+ ,h_- \in A_-$. Now we can write 
$l_+ =l_+ ' +h_+ '$, where $l_+ '\in L_+ ', h_+ '\in A_+$. So now we 
have $\aa =l_+ '+h_+ '+h_+ +h_-$. Since $h_+ '+h_+ \in L\cap L_- '$ and
$ h_- \in L\cap L_+ '$ we have $l_+ '\in L$. Thus 
$l_+ '\in L\cap L_+ '$
 and we conclude that $\aa\in (L\cap L_+ ')+(L\cap L_- ')$. This shows
that $L$ is compatible.

This completes the proof that $\FL m\sub\FHL m$. To prove the reverse
inclusion $\FHL m\sub\FL m$ we will make use of theorem \ref{thm.T2bl}
which says that $\FL m=\Fas m$. Thus we want to show $\FHL m\sub\Fas m$.
Since the ideal $(I\Lg^L )^m$ is generated by elements of the form 
$(1-h_1 )\cdots (1-h_m )$, where $h_i$ is a Dehn twist along a simple
 closed curve $l_i$ representing an
 element of $L$, it suffices to prove the following lemma.
\begin{lemma}
\lbl{lem.hl}
Let $l_1 ,\cdots ,l_m$ be simple closed curves in $\S$ representing
 elements of a Lagrangian $L$ which is compatible with an admissible 
embedding $i:\S\sub M$. Let $l_j '$ be translates of $i(l_j )$ into 
disjoint parallel copies of $i(\S )$ in $M$. Then
 the link $\{ l_1 ',\cdots ,l_m '\}$ is algebraically split in $M$. 
\end{lemma}
\begin{pf} Write $[l_i ]=\l_i^+ +\l_i^-$, where 
$\l_i^{\pm}\in L\cap L_{\pm}$. Then $\l_i^{\pm}$ bounds a surface 
$N_i^{\pm}$ in $M_{\pm}$ (translated). Suppose that $l_i '$ lies 
on the $M_+$ side of $l_j '$. Then $N_i^+$ is disjoint from 
$N_j^-$ and the linking number of $l_i '$ and $l_j '$ is just 
the intersection number of $N_i^-$ with $N_j^+$. But this is the 
same as the intersection number of $\l_i^-$ with $\l_j^+$ in $\S$.
 Since these both represent
 elements of the Lagrangian $L$, the intersection number is zero.
\end{pf}
This completes the proof of proposition \ref{prop.heeg}.

\subsection{Proof of theorem \ref{thm.T2bl}(b) and proposition 
\ref{thm.boundary}}
\lbl{sub.prfb}

We will need the following observation. Suppose that $M=M_1\cup 
M_2$, where $N=\partial M_1 =\partial M_2$. Let $f, g$ be 
diffeomorphisms of $N$ and consider the manifold $M_h$ where 
$h=fg$. We can describe $M_h$ alternatively using $f$ and $g$ 
separately by writing $M=M_1\cup (I\times N)\cup M_2$, splitting 
$M$ along two parallel copies of $N\subseteq M$. Then it is easy to 
see that $M_h =M_1\cup_f (I\times N)\cup_g M_2$.

We first prove proposition \ref{thm.boundary}. This will follow from 
the preceding observation and a theorem of Morita (\cite[proposition 
2.3]{Mo}) which says that any homology 3-sphere $M$ can be written 
in the form $S^3_f$ for some $f\in\Kg$, for some $g$. Let us write 
$g=g_1\cdots g_k$, where each $g_i$ is a Dehn twist along some 
bounding simple closed curve $\gamma_i$ in the Heegaard surface 
$\Sigma_g$ in $S^3$. According to the observation we can split 
$S^3$ 
along $k$ parallel copies of $\Sigma_g$ in $S^3$ and obtain $M$ by 
simultaneously regluing by Dehn twists along the now disjoint 
copies of $\gamma_i$ in the parallel copies of $\Sigma_g$. This is 
the same as doing simultaneous $\pm 1$ surgeries along the link $L$ 
formed by these disjoint copies of $\gamma_i$. But since each 
$\gamma_i$ is a bounding simple closed curve we see that $L$ is, in 
fact, a boundary link.  \qed

We now turn to theorem \ref{thm.T2bl}(b). Suppose $[M,L,f]$ is a 
generator of $\Fb m$. The components of $L$ bound disjoint Seifert 
surfaces $V_1 ,\ldots ,V_m$. Let $\Sigma$ be a connected sum of 
the boundaries of tubular neighborhoods of these surfaces. Thus 
$\Sigma$ is an admissible surface in $M$ and the components of $L$ 
are disjoint bounding simple closed curves $\gamma_i$ on $\Sigma$. 
Let $h_1 ,\ldots ,h_m$  be the diffeomorphisms of $\Sigma$ defined 
by Dehn twists along $\gamma_1 ,\ldots ,\gamma_m$, so that, according
to remark \ref{rem.sursame}, 
cutting $M$ along $\Sigma$ and regluing using $h_i$ is the same as 
framed surgery along $\gamma_i$ using $f|\gamma_i$. Note that the 
$h_i$ commute with each other. If $L'$ is any sublink of $L$, then, 
since the $\gamma_i$ are disjoint $M_{L',f|L'}$ is obtained from $M$ 
by cutting along $\Sigma$ and regluing using the composition of 
those $h_i$ corresponding to the components which appear in $L'$. 
Now $[M,L,f] = \sum_{L' \subseteq L} (-1)^{|L'|} M_{L',f'}$, which is, 
therefore, the image, under the map $\Kg\to\M$ of the sum 
$$\sum_{1\leq i_1\leq\cdots\leq i_k\leq m} (-1)^k h_{i_1}\cdots 
h_{i_k}= (1-h_1 )\cdots (1-h_m )$$
But this is obviously an element of $I\Kg$.

Conversely suppose $\lambda\in (I\Kg )^m$. Then $\lambda$ is a 
linear combination of elements of the form $h(1-g_1 )\cdots (1-g_m 
)$, where $h, g_i\in\Kg$.  Now each $g_i$ is a product 
$\rho_1\cdots\rho_k$, where each $\rho_i$ is a Dehn twist along a 
bounding simple closed curve. Using the identity 
$\rho_1\cdots\rho_k -1=\sum_1^k\rho_1\cdots\rho_{i-1}(\rho_i 
-1)$ and the normality of $\Kg$ again we can assume that each $g_i$ 
is a Dehn twist along a bounding simple closed curve. Now suppose 
$\Sigma$ is an admissible surface in an \ihs $M$. Let $\Sigma_0 
,\ldots ,\Sigma_m$ be parallel copies of $\Sigma$ in $M$. If $g_i$ is 
a Dehn twist along $\gamma_i\in\Sigma$, then any 
$M_{hg_{i_1}\cdots g_{i_k}}$ can be obtained by cutting $M$ along 
$\Sigma_0 ,\Sigma_{i_1},\ldots ,\Sigma_{i_k}$ and then regluing 
using $h,g_{i_1},\cdots ,g_{i_k}$ on these copies of $\Sigma$. But 
this is the same as considering the link $L$ in $M_h$ defined by the 
$\{\gamma_i\in\Sigma_i\}$ with the framing $f$ defined by the 
directions of the twists given by the $\{ g_i\}$. From this we see 
that the image of $h(1-g_1 )\cdots (1-g_m )$ in $\M$ is exactly 
$[M_h ,L,f]$. \qed

\subsection{Proof of theorem \ref{thm.T2bl}(a)}
\lbl{sub.prfa}
The proof of theorem \ref{thm.T2bl}(a) is very similar to that just 
given for (b). We need to use the result of Johnson that $\Tg$ is 
generated by what he calls BP maps if $g\geq 3$. A BP {\em map} is 
obtained by doing Dehn twists along two disjoint simple closed 
curves in $\Sigma_g$ which form a bounding pair, i.e., they are 
homologous or, equivalently, form the boundary of a subsurface of 
$\Sigma_g$. The twists are in opposite directions. Note that the 
$2$-component link associated to a BP map is nothing but a $1$-pair 
blink.
This gives the main motivation for the role of  blinks in the 
present
paper.  Suppose we have 
a product $g_1\cdots g_m$ of BP maps. If we have $m$ parallel 
copies of an admissible surface $\Sigma$ in an \ihs\ $M$ and, in the 
$i$-th copy, a pair of such curves associated to $g_i$, then the 
totality of these curves forms a blink $L_{bl}$ with a unit Seifert 
framing $f$ defined by the directions of the twists. Then, just as 
with boundary links above, we see that the manifold obtained by 
cutting $M$ along $\Sigma$ and regluing by $g_1\cdots g_m$ is 
homeomorphic to  $M_{L_{bl},f}$. The proof of theorem 
\ref{thm.T2bl}(a) now is identical to the proof of theorem 
\ref{thm.T2bl}(b) with the substitution of blinks for boundary links 
and $\Tg$ for $\Kg$.

\subsection{Proof of theorem \ref{thm.T2bl}(c)}
\lbl{sub.prfc}
 Suppose $[M,L,f]$ is a generator of $\Fas m$. Consider a tubular 
neighborhood of $L$ and define a connected submanifold 
$M_1\sub M$ by connecting these components by solid tubes 
in $M$. For each component $l_i$ of $L$ there is a canonical 
longitude $\lambda_i\subseteq\partial M_1$ defined by the 
requirement that $\lambda_i$ be null-homologous in the 
complement of $l_i$. Since $L$ is algebraically split, 
$\{\lambda_i\}$ span Ker$\{ H_1 (\partial M_1 )\to H_1 
(\overline{M-
M_1})\}$ which is a Lagrangian in $H_1 (\partial M_1 )$. The result 
of 
doing $\pm 1$-surgery on any sublink of $L$ is the same as cutting 
$M$
along $\partial M_1$ and regluing 
by simultaneous Dehn twists on the corresponding 
$\lambda_i\subseteq\partial M_1$. If we choose an identification 
$\S_g\cong\partial M_1$, these Dehn twists define elements 
$h_i\in\Lg^i$, where $i$ is the composition $\S_g\cong\partial 
M_1\subseteq M$. As in section~\ref{sub.prfb} we see that $[M,L,f]$ 
is the image, under the map $\Lg\to\M$ defined by $i$, of $(1-h_1 
)\cdots (1-h_m )$.

We now prove the converse statement.
Suppose $\lambda\in(I\Lg )^m$. Then $\lambda$ is a linear 
combination of elements of the form $h(1-g_1 )\cdots (1-g_m )$, 
where $h, g_i\in\Lg$. Now writing each $g_i$ as a product of $L$-
twists
 and using repeatedly the identity $1-
gh=g(1-h)+1-g$, we can write this as a sum of elements of the form 
$h(1-\rho_1 )\cdots (1-\rho_m )$, where each $\rho_i$ is an $L$-
twist. 

We now proceed as in Sections~\ref{sub.prfb} and \ref{sub.prfa}. We 
see that the image of \hfil\break $h(1-\rho_1 )\cdots (1-\rho_m )$ 
in $\M$ is 
an element $[M,K,f]$ which can be described as follows. We begin 
with an embedded surface $\Sigma$ in some \ihs\ $M'$ and let 
$M_1$ 
be the closure of one of the components of $M'-\Sigma$, so that 
$L=\ker\{H_1 (\Sigma )\to H_1 (M_1 )\}$. Let $\Sigma 
,\Sigma_1 ,\cdots ,\Sigma_m$ be parallel copies of $\Sigma$ 
appearing in that order as we move away from $M_1$. Then $\rho_i$ 
is a Dehn twist on a curve $\lambda_i\sub\Sigma_i$ and 
$\lambda_i$ is null-homologous in $M_1$. Then $M$ is obtained by 
cutting $M'$ along $\Sigma$ and reattaching with $h$, $K$ is the link 
consisting of the $\{\lambda_i\}$ and $f$ is the framing defined by 
the signs of the Dehn twists. To complete the proof we need to see 
that $K$ is algebraically split. But since $\lambda_i$ bounds a chain 
in $M_1$ and $h$ does not alter this, it follows that the linking 
number of $\lambda_i$ with any $\lambda_j$, when $j>i$, is zero.

\section{Equivalence of $AS$, $B$ and $BL$ filtrations}
\lbl{sec.equiv}

In this section we prove theorems \ref{thm.bl2as}, \ref{thm.as2b} 
and
\ref{thm.as2bl}. 

%%% \lbl{sub.blinks}

\subsection{A review of relations on $\M$ in graphical and algebraic 
form}
\lbl{sub.review}

In this section we review the forms of notation used in this paper
as well as some important identities among the elements
 of $\M$. We follow the conventions of \cite{Oh}, \cite{GL1}, 
\cite{GO1}.
 For the convenience of the uninitiated reader, we review them again 
here.
The notation comes in two forms:  {\em algebraic} and
{\em graphical}. By algebraic notation we mean 
$[M,L,f]$. Note that the various filtrations on $\M$ have been
written in algebraic notation.
The rules of the graphical notation are summarized in figure
\ref{f.bands}. Since we are talking about 
identities of (linear combinations) of
\ihs s obtained by framed links in other \ihs s, it is almost 
unavoidable 
to use graphical notation to represent framed links, and surgeries 
on them.

\begin{remark}
\lbl{rem.conventions}
We mention 
once and for all that any link (whether \as\ or boundary 
or blink) drawn in a figure corresponds to a {\em linear combination}
of \ihs s, and therefore represents an element of $\M$. The figures 
represent identities of these
elements.
We cannot stress too strongly the fact that in papers prior to
the ones talking about \fti s 
figures corresponded to links or 3-manifolds, but never to {\em 
linear
combinations} of them. Nevertheless, this point of view is
very fundamental in the world of \fti s.
\end{remark}

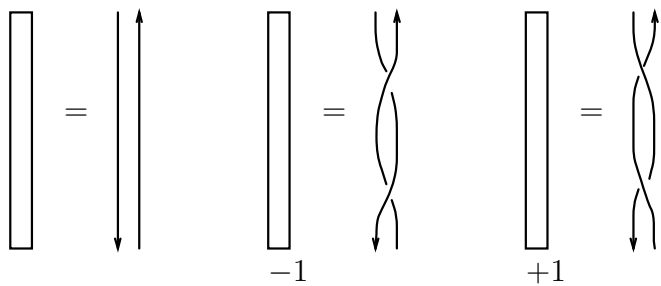
\begin{figure}[htpb]
$$\printname{bands}
	\setlength{\unitlength}{0.03\standardunitlength}
	\begin{array}{c}  \hspace{-1.7mm}
        	\raisebox{-8pt}{\begingroup\makeatletter\ifx\SetFigFont\undefined
% extract first six characters in \fmtname
\def\x#1#2#3#4#5#6#7\relax{\def\x{#1#2#3#4#5#6}}%
\expandafter\x\fmtname xxxxxx\relax \def\y{splain}%
\ifx\x\y   % LaTeX or SliTeX?
\gdef\SetFigFont#1#2#3{%
  \ifnum #1<17\tiny\else \ifnum #1<20\small\else
  \ifnum #1<24\normalsize\else \ifnum #1<29\large\else
  \ifnum #1<34\Large\else \ifnum #1<41\LARGE\else
     \huge\fi\fi\fi\fi\fi\fi
  \csname #3\endcsname}%
\else
\gdef\SetFigFont#1#2#3{\begingroup
  \count@#1\relax \ifnum 25<\count@\count@25\fi
  \def\x{\endgroup\@setsize\SetFigFont{#2pt}}%
  \expandafter\x
    \csname \romannumeral\the\count@ pt\expandafter\endcsname
    \csname @\romannumeral\the\count@ pt\endcsname
  \csname #3\endcsname}%
\fi
\fi\endgroup
\begin{picture}(9020,3810)(0,-10)
\thicklines
\path(3612,3783)(3612,483)(3912,483)
	(3912,3783)(3612,3783)
\path(12,3783)(12,483)(312,483)
	(312,3783)(12,3783)
\path(1849.500,3648.000)(1812.000,3783.000)(1774.500,3648.000)
\path(1812,3783)(1812,483)
\path(1512,3783)(1512,483)
\path(1474.500,618.000)(1512.000,483.000)(1549.500,618.000)
\path(7212,3783)(7212,483)(7512,483)
	(7512,3783)(7212,3783)
\spline(5112,3783)
(5112,3408)(5187,3108)(5262,2958)
\spline(5337,1158)
(5412,933)(5412,633)(5412,483)
\spline(5337,2658)
(5412,2358)(5412,1908)
	(5412,1608)(5337,1308)(5187,1008)
	(5112,858)(5112,483)
\path(5074.500,618.000)(5112.000,483.000)(5149.500,618.000)
\path(5449.500,3648.000)(5412.000,3783.000)(5374.500,3648.000)
\spline(5412,3783)
(5412,3333)(5412,3108)
	(5262,2808)(5187,2583)(5112,2283)
	(5112,1758)(5262,1383)
\spline(8712,3783)
(8712,3408)(8787,3108)
	(8937,2733)(9012,2433)(9012,2058)
	(9012,1683)(8937,1458)
\spline(8787,2883)
(8712,2658)(8712,2358)
	(8712,2058)(8712,1758)(8787,1533)
	(8862,1308)(8937,1083)(9012,933)
	(9012,633)(9012,483)
\spline(8787,1308)
(8712,1083)(8712,783)(8712,483)
\path(8674.500,618.000)(8712.000,483.000)(8749.500,618.000)
\path(9049.500,3648.000)(9012.000,3783.000)(8974.500,3648.000)
\spline(9012,3783)
(9012,3408)(8862,3033)
\put(762,2283){\makebox(0,0)[lb]{$=$}}
\put(4362,2283){\makebox(0,0)[lb]{$=$}}
\put(7962,2283){\makebox(0,0)[lb]{$=$}}
\put(7212,33){\makebox(0,0)[lb]{$+1$}}
\put(3612,33){\makebox(0,0)[lb]{$-1$}}
\end{picture} }
        	\hspace{-1.9mm}
	\end{array}
 $$
 \caption{Some drawing conventions for bands.
Shown here
are  ribbon parts of $AS$-admissible links that represent 
(linear combinations of) \ihs s.
The numbers in the bottom of each band indicate the 
number of twists that we put in the band.}
 \lbl{f.bands}
\end{figure}

We are now ready to review identities of elements of $\M$. 
We begin with the following fundamental identity, (in algebraic
notation) for an \ASA\ link $(L,f)$
in an \ihs\ $M$:

Let $L_{bl} \cup L$ be  the union of an \ASA\ link $L$ and a 
\BLA\ blink $L_{bl}$ (in the complement of $L$). Let $l$
denote either a component of $L$, or a $1$-pair of $L_{bl}$.
Define 
$$L' =\cases L & \text{ if } l\notin L \\
            L-l & \text{ if } l\in L \endcases
$$
and 
$$L_{bl}' =\cases L_{bl} & \text{ if } l\notin L_{bl} \\
            L_{bl}-l & \text{ if } l\in L_{bl} \endcases
$$

Let  $f'$ (respectively, $f|_{l}$) denote
 the restriction of the framing $f$ of $L_{bl} \cup L$ to 
$L_{bl}' \cup L'$ (respectively $l$).
Then we have the following fundamental relation:
\begin{equation}
\lbl{eq.fundamental}
[M,L_{bl} \cup L,f]=[M,L_{bl}' \cup L',f']-
[M_{(l,f|_{l})}L_{bl}' \cup L', f']
\end{equation}
The proof of the above equation follows by definition of the 
symbol $[M, L_{bl} \cup L,f]$ and the following exercise, left to the
reader:

\begin{exercise}
Show that the unit-Seifert framing $f|L_{bl}$ of $L_{bl}$ in $M$
is the same as the one of $L_{bl}$ in $M_{(l,f|_{l})}$.
\end{exercise}

\begin{remark}
\lbl{rem.oldfundamental}
This extends the notation of previous papers
 \cite{Ga}, \cite{GL1}, \cite{GL2},
\cite{GO1}, \cite{GO2}, where $L_{bl}$ is the empty blink.
If $l$ is a knot that
 bounds a disk $D$ in $M$, then $M_{(l,f|_{l})}$ is diffeomorphic to
$M$ and we may construct 
the link in $M$ corresponding to $L_{bl}' \cup L'$ in $M_{(l,f|_{l})}$ 
from
 $L_{bl}' \cup L'$ by 
just giving 
the bundle of strands of $L_{bl}' \cup L'$ which pass through $D$ a 
full
 clockwise twist if $f=+1$ or counterclockwise twist if $f=-1$.   
\end{remark}

Examples of equation \eqref{eq.fundamental} in graphical notation
are given in figures \ref{f.bands1},
\ref{f.onecross} and \ref{f.twocross}. 

\begin{figure}[htpb]
$$\printname{bands1}
	\setlength{\unitlength}{0.03\standardunitlength}
	\begin{array}{c}  \hspace{-1.7mm}
        	\raisebox{-8pt}{\begingroup\makeatletter\ifx\SetFigFont\undefined
% extract first six characters in \fmtname
\def\x#1#2#3#4#5#6#7\relax{\def\x{#1#2#3#4#5#6}}%
\expandafter\x\fmtname xxxxxx\relax \def\y{splain}%
\ifx\x\y   % LaTeX or SliTeX?
\gdef\SetFigFont#1#2#3{%
  \ifnum #1<17\tiny\else \ifnum #1<20\small\else
  \ifnum #1<24\normalsize\else \ifnum #1<29\large\else
  \ifnum #1<34\Large\else \ifnum #1<41\LARGE\else
     \huge\fi\fi\fi\fi\fi\fi
  \csname #3\endcsname}%
\else
\gdef\SetFigFont#1#2#3{\begingroup
  \count@#1\relax \ifnum 25<\count@\count@25\fi
  \def\x{\endgroup\@setsize\SetFigFont{#2pt}}%
  \expandafter\x
    \csname \romannumeral\the\count@ pt\expandafter\endcsname
    \csname @\romannumeral\the\count@ pt\endcsname
  \csname #3\endcsname}%
\fi
\fi\endgroup
\begin{picture}(10434,2685)(0,-10)
\thicklines
\path(2850,558)(2850,2658)(3150,2658)
	(3150,558)(2850,558)
\path(4350,558)(4350,2658)(4650,2658)
	(4650,558)(4350,558)
\path(8550,558)(8550,2658)(8850,2658)
	(8850,558)(8550,558)
\path(10050,558)(10050,2658)(10350,2658)
	(10350,558)(10050,558)
\path(6450,1758)(6450,2658)(6750,2658)(6750,1758)
\path(6450,1533)(6450,558)(6750,558)(6750,1533)
\path(750,1758)(750,2658)(1050,2658)(1050,1758)
\path(750,1533)(750,558)(1050,558)(1050,1533)
\spline(6375,1908)
(6150,1908)(6000,1758)
	(6150,1608)(7050,1608)(7200,1758)
	(7050,1908)(6825,1908)
\spline(675,1908)
(450,1908)(300,1758)
	(450,1608)(1350,1608)(1500,1758)
	(1350,1908)(1125,1908)
\put(1950,1683){\makebox(0,0)[lb]{$=$}}
\put(3600,1683){\makebox(0,0)[lb]{$-$}}
\put(7725,1683){\makebox(0,0)[lb]{$=$}}
\put(9300,1683){\makebox(0,0)[lb]{$-$}}
\put(4350,33){\makebox(0,0)[lb]{$-1$}}
\put(10050,33){\makebox(0,0)[lb]{$+1$}}
\put(5700,1383){\makebox(0,0)[lb]{$-1$}}
\put(0,1383){\makebox(0,0)[lb]{$+1$}}
\end{picture} }
        	\hspace{-1.9mm}
	\end{array}
 $$
 \caption{A special case of equation \eqref{eq.fundamental} in a 
graphical way.} 
 \lbl{f.bands1}
\end{figure}

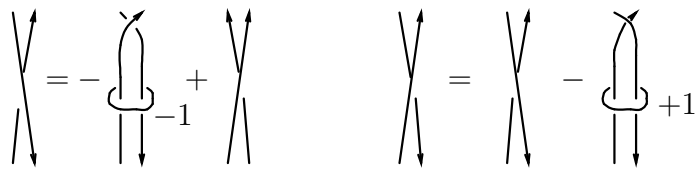
\begin{figure}[htpb]
$$\printname{onecross}
	\setlength{\unitlength}{0.03\standardunitlength}
	\begin{array}{c}  \hspace{-1.7mm}
        	\raisebox{-8pt}{\begingroup\makeatletter\ifx\SetFigFont\undefined
% extract first six characters in \fmtname
\def\x#1#2#3#4#5#6#7\relax{\def\x{#1#2#3#4#5#6}}%
\expandafter\x\fmtname xxxxxx\relax \def\y{splain}%
\ifx\x\y   % LaTeX or SliTeX?
\gdef\SetFigFont#1#2#3{%
  \ifnum #1<17\tiny\else \ifnum #1<20\small\else
  \ifnum #1<24\normalsize\else \ifnum #1<29\large\else
  \ifnum #1<34\Large\else \ifnum #1<41\LARGE\else
     \huge\fi\fi\fi\fi\fi\fi
  \csname #3\endcsname}%
\else
\gdef\SetFigFont#1#2#3{\begingroup
  \count@#1\relax \ifnum 25<\count@\count@25\fi
  \def\x{\endgroup\@setsize\SetFigFont{#2pt}}%
  \expandafter\x
    \csname \romannumeral\the\count@ pt\expandafter\endcsname
    \csname @\romannumeral\the\count@ pt\endcsname
  \csname #3\endcsname}%
\fi
\fi\endgroup
\begin{picture}(9396,2139)(0,-10)
\thicklines
\path(12,12)(87,762)
\path(1512,2112)(1587,2037)
\path(1512,1212)(1512,912)
\path(1512,687)(1512,12)
\path(12,2112)(312,12)
\path(265.331,126.551)(312.000,12.000)(324.728,135.037)
\path(1812,687)(1812,12)
\path(1782.000,132.000)(1812.000,12.000)(1842.000,132.000)
\path(3012,12)(3312,2112)
\path(3324.728,1988.963)(3312.000,2112.000)(3265.331,1997.449)
\path(3062.982,1999.302)(3012.000,2112.000)(3003.950,1988.569)
\path(3012,2112)(3162,1287)
\path(320.050,1988.569)(312.000,2112.000)(261.018,1999.302)
\path(312,2112)(162,1287)
\path(3312,12)(3237,912)
\path(5412,2112)(5562,1212)
\path(6987,912)(6912,12)
\path(8712,1512)(8712,912)
\path(8412,1362)(8412,912)
\path(8412,687)(8412,12)
\path(8712,687)(8712,12)
\path(8682.000,132.000)(8712.000,12.000)(8742.000,132.000)
\path(5412,12)(5712,2112)
\path(5724.728,1988.963)(5712.000,2112.000)(5665.331,1997.449)
\path(6912,2112)(7212,12)
\path(7165.331,126.551)(7212.000,12.000)(7224.728,135.037)
\path(7220.050,1988.569)(7212.000,2112.000)(7161.018,1999.302)
\path(7212,2112)(7062,1287)
\path(5637,912)(5712,12)
\path(5672.138,129.094)(5712.000,12.000)(5731.931,134.077)
\path(8637,2037)(8712,2112)
\path(8648.360,2005.934)(8712.000,2112.000)(8605.934,2048.360)
\path(1737,1887)	(1768.506,1832.802)
	(1790.336,1792.080)
	(1812.000,1737.000)

\path(1812,1737)	(1818.001,1698.979)
	(1820.817,1653.080)
	(1821.185,1601.977)
	(1819.838,1548.345)
	(1817.511,1494.857)
	(1814.939,1444.188)
	(1812.857,1399.011)
	(1812.000,1362.000)

\path(1812,1362)	(1812.000,1295.898)
	(1812.000,1254.543)
	(1812.000,1205.895)
	(1812.000,1148.635)
	(1812.000,1081.445)
	(1812.000,1043.714)
	(1812.000,1003.006)
	(1812.000,959.156)
	(1812.000,912.000)

\path(1437,1062)	(1384.564,1018.271)
	(1362.000,987.000)

\path(1362,987)	(1347.030,912.817)
	(1349.512,871.867)
	(1362.000,837.000)

\path(1362,837)	(1432.406,786.420)
	(1475.708,771.441)
	(1512.000,762.000)

\path(1512,762)	(1562.354,756.710)
	(1625.194,757.298)
	(1687.687,760.237)
	(1737.000,762.000)

\path(1737,762)	(1810.800,756.195)
	(1852.602,755.469)
	(1887.000,762.000)

\path(1887,762)	(1930.305,793.695)
	(1962.000,837.000)

\path(1962,837)	(1964.902,875.100)
	(1962.000,912.000)

\path(1962,912)	(1964.902,948.900)
	(1962.000,987.000)

\path(1962,987)	(1939.436,1018.271)
	(1887.000,1062.000)

\path(1740.816,2010.843)(1812.000,2112.000)(1701.587,2056.242)
\path(1812,2112)	(1763.867,2070.036)
	(1722.769,2033.344)
	(1659.041,1973.141)
	(1615.543,1926.118)
	(1587.000,1887.000)

\path(1587,1887)	(1563.778,1838.934)
	(1541.505,1776.773)
	(1523.229,1713.474)
	(1512.000,1662.000)

\path(1512,1662)	(1504.856,1595.127)
	(1503.070,1553.580)
	(1502.475,1504.868)
	(1503.070,1447.672)
	(1504.856,1380.674)
	(1506.196,1343.088)
	(1507.833,1302.556)
	(1509.768,1258.915)
	(1512.000,1212.000)

\path(8412,2112)	(8459.713,2086.879)
	(8500.511,2064.457)
	(8563.999,2025.952)
	(8637.000,1962.000)

\path(8637,1962)	(8663.108,1914.095)
	(8685.060,1851.915)
	(8701.733,1788.528)
	(8712.000,1737.000)

\path(8712,1737)	(8715.572,1703.564)
	(8716.763,1658.434)
	(8715.572,1596.337)
	(8714.084,1557.278)
	(8712.000,1512.000)

\path(8562,1962)	(8532.755,1906.081)
	(8511.679,1864.785)
	(8487.000,1812.000)

\path(8487,1812)	(8467.761,1763.390)
	(8445.045,1701.289)
	(8424.557,1638.293)
	(8412.000,1587.000)

\path(8412,1587)	(8408.428,1553.564)
	(8407.237,1508.434)
	(8408.428,1446.337)
	(8409.916,1407.278)
	(8412.000,1362.000)

\path(8337,1062)	(8284.564,1018.271)
	(8262.000,987.000)

\path(8262,987)	(8247.030,912.817)
	(8249.512,871.867)
	(8262.000,837.000)

\path(8262,837)	(8332.406,786.420)
	(8375.708,771.441)
	(8412.000,762.000)

\path(8412,762)	(8462.354,756.710)
	(8525.194,757.298)
	(8587.687,760.237)
	(8637.000,762.000)

\path(8637,762)	(8710.800,756.195)
	(8752.602,755.469)
	(8787.000,762.000)

\path(8787,762)	(8830.305,793.695)
	(8862.000,837.000)

\path(8862,837)	(8864.903,875.100)
	(8862.000,912.000)

\path(8862,912)	(8864.903,948.900)
	(8862.000,987.000)

\path(8862,987)	(8839.436,1018.271)
	(8787.000,1062.000)

\put(1962,537){\makebox(0,0)[lb]{$-1$}}
\put(2412,1062){\makebox(0,0)[lb]{$+$}}
\put(912,1062){\makebox(0,0)[lb]{$-$}}
\put(462,1062){\makebox(0,0)[lb]{$=$}}
\put(6087,1062){\makebox(0,0)[lb]{$=$}}
\put(7662,1062){\makebox(0,0)[lb]{$-$}}
\put(9012,687){\makebox(0,0)[lb]{$+1$}}
\end{picture} }
        	\hspace{-1.9mm}
	\end{array}
 $$
 \caption{ Another special case of equation \eqref{eq.fundamental}.
The figure represents an identity in $\cal M$. There are
{\em two} interpretations of the above figure. Either each of the 
crossings
shown belong to the same component (of an algebraically split link),
or each of the crossings shown are part of a ribbon of a piece of a 
Seifert
surface of a $1$-pair blink. It will be clear each time we use the 
identity
shown in the figure which interpretation we have in mind.}
 \lbl{f.onecross}
\end{figure}

\begin{figure}[htpb]
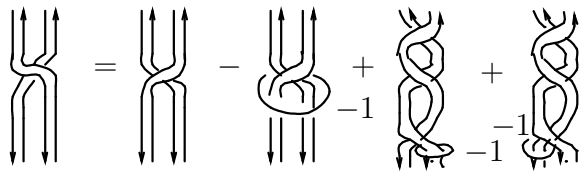

$$\printname{twocross}
	\setlength{\unitlength}{0.03\standardunitlength}
	\begin{array}{c}  \hspace{-1.7mm}
        	\raisebox{-8pt}{\input draws/twocross.tex }
        	\hspace{-1.9mm}
	\end{array}
 $$
 \caption{Another special case of equation \eqref{eq.fundamental}.
The figure represents an identity in $\cal M$.  Crossings are in the 
same
component of the link.}
 \lbl{f.twocross}
\end{figure}

\begin{figure}[htpb]
$$\printname{3bands}
	\setlength{\unitlength}{0.03\standardunitlength}
	\begin{array}{c}  \hspace{-1.7mm}
        	\raisebox{-8pt}{\begingroup\makeatletter\ifx\SetFigFont\undefined
% extract first six characters in \fmtname
\def\x#1#2#3#4#5#6#7\relax{\def\x{#1#2#3#4#5#6}}%
\expandafter\x\fmtname xxxxxx\relax \def\y{splain}%
\ifx\x\y   % LaTeX or SliTeX?
\gdef\SetFigFont#1#2#3{%
  \ifnum #1<17\tiny\else \ifnum #1<20\small\else
  \ifnum #1<24\normalsize\else \ifnum #1<29\large\else
  \ifnum #1<34\Large\else \ifnum #1<41\LARGE\else
     \huge\fi\fi\fi\fi\fi\fi
  \csname #3\endcsname}%
\else
\gdef\SetFigFont#1#2#3{\begingroup
  \count@#1\relax \ifnum 25<\count@\count@25\fi
  \def\x{\endgroup\@setsize\SetFigFont{#2pt}}%
  \expandafter\x
    \csname \romannumeral\the\count@ pt\expandafter\endcsname
    \csname @\romannumeral\the\count@ pt\endcsname
  \csname #3\endcsname}%
\fi
\fi\endgroup
\begin{picture}(10461,12810)(0,-10)
\thicklines
\path(2930,7353)(2930,8283)(3080,8283)(3080,7353)
\path(3530,7323)(3530,8253)(3680,8253)(3680,7323)
\path(4130,6783)(4130,8283)(4280,8283)(4280,6783)
\path(2900,7083)(2900,6663)
\path(3080,7083)(3080,6663)
\path(2930,6453)(2930,4953)(3080,4953)(3080,6483)
\path(3530,6453)(3530,4953)(3680,4953)(3680,6483)
\path(3170,6963)(4040,6963)
\path(3530,7113)(3530,6993)
\path(3680,7113)(3680,7023)
\path(3530,6873)(3530,6663)
\path(3680,6873)(3680,6693)
\path(4130,6453)(4130,4953)(4280,4953)(4280,6483)
\path(5600,7353)(5600,8283)(5750,8283)(5750,7353)
\path(6200,7353)(6200,8283)(6350,8283)(6350,7353)
\path(5630,6933)(5630,4983)(5780,4983)(5780,6933)
\path(6230,6963)(6230,5973)
\path(6380,6993)(6380,5973)
\path(6230,5733)(6230,4983)(6380,4983)(6380,5733)
\path(5870,7473)(6080,7473)
\path(6800,6063)(6800,8283)(6980,8283)(6980,6033)
\path(6800,5733)(6800,4983)(6980,5013)(6980,5763)
\path(8630,7353)(8630,8283)(8780,8283)(8780,7353)
\path(9830,7353)(9830,8283)(9980,8283)(9980,7353)
\path(8900,7563)(9680,7563)
\path(9230,7683)(9230,8283)(9380,8283)(9380,7683)
\path(9230,5733)(9230,4983)(9380,4983)(9380,5733)
\path(9800,5703)(9800,4953)(9980,4983)(9980,5733)
\path(8630,7053)(8630,4983)(8780,4983)(8780,7053)
\path(9230,7053)(9230,5973)
\path(9380,7053)(9380,5943)
\path(9830,7053)(9830,5943)
\path(9980,7053)(9980,5913)
\path(9230,7473)(9230,7263)(9260,7263)
\path(9380,7503)(9380,7263)
\path(3170,7503)(3470,7503)
\path(6440,6123)(6740,6123)
\path(9470,6093)(9770,6093)
\spline(2720,6993)
(2480,6993)(2390,6843)
	(2420,6663)(2510,6603)(4490,6603)
	(4580,6693)(4610,6813)(4580,6933)
	(4460,6993)(4370,6993)
\spline(8570,7563)
(8330,7563)(8240,7413)
	(8270,7233)(8360,7173)(10340,7173)
	(10430,7263)(10460,7383)(10430,7503)
	(10310,7563)(10220,7563)
\spline(6200,6153)
(6020,6153)(5930,6003)
	(6020,5853)(7220,5853)(7310,6003)
	(7220,6153)(7040,6153)
\spline(9200,6123)
(9020,6123)(8930,5973)
	(9020,5823)(10220,5823)(10310,5973)
	(10220,6123)(10040,6123)
\spline(2900,7503)
(2720,7503)(2630,7353)
	(2720,7203)(3920,7203)(4010,7353)
	(3920,7503)(3740,7503)
\spline(5540,7443)
(5360,7443)(5270,7293)
	(5360,7143)(6560,7143)(6650,7293)
	(6560,7443)(6380,7443)
\put(1790,6753){\makebox(0,0)[lb]{$-$}}
\put(5030,6723){\makebox(0,0)[lb]{$-$}}
\put(7910,6753){\makebox(0,0)[lb]{$-$}}
\put(2840,4533){\makebox(0,0)[lb]{$+1$}}
\put(6140,4533){\makebox(0,0)[lb]{$+1$}}
\put(6770,4533){\makebox(0,0)[lb]{$+1$}}
\put(5540,4533){\makebox(0,0)[lb]{$+1$}}
\put(8540,4533){\makebox(0,0)[lb]{$+1$}}
\put(9140,4533){\makebox(0,0)[lb]{$+1$}}
\put(9740,4533){\makebox(0,0)[lb]{$+1$}}
\put(4040,4533){\makebox(0,0)[lb]{$+1$}}
\put(3440,4533){\makebox(0,0)[lb]{$+1$}}
\path(4130,3783)(4130,453)(4280,453)
	(4280,3783)(4130,3783)
\path(2930,2283)(2930,3783)(3080,3783)(3080,2283)
\path(2930,1983)(2930,483)(3080,483)(3080,1983)
\path(3530,2283)(3530,3783)(3680,3783)(3680,2283)
\path(3530,1953)(3530,453)(3680,453)(3680,1953)
\path(5630,1983)(5630,483)(5780,483)(5780,1983)
\path(5630,2253)(5630,3753)(5780,3753)(5780,2253)
\path(6830,2283)(6830,3783)(6980,3783)(6980,2283)
\path(6830,1953)(6830,453)(6980,453)(6980,1953)
\path(5900,2493)(6680,2493)
\path(6230,1953)(6230,453)(6380,453)(6380,1953)
\path(6230,2613)(6230,3753)(6380,3753)(6380,2613)
\path(8630,3783)(8630,453)(8780,453)
	(8780,3783)(8630,3783)
\path(9830,2283)(9830,3783)(9980,3783)(9980,2283)
\path(9230,2283)(9230,3783)(9380,3783)(9380,2283)
\path(9230,1983)(9230,483)(9380,483)(9380,1983)
\path(9830,1983)(9830,483)(9980,483)(9980,1983)
\path(3140,2403)(3440,2403)
\path(9440,2403)(9770,2403)
\path(6215,2433)(6215,2208)
\path(6365,2433)(6365,2208)
\spline(5570,2493)
(5330,2493)(5240,2343)
	(5270,2163)(5360,2103)(7340,2103)
	(7430,2193)(7460,2313)(7430,2433)
	(7310,2493)(7220,2493)
\spline(2870,2403)
(2690,2403)(2600,2253)
	(2690,2103)(3890,2103)(3980,2253)
	(3890,2403)(3710,2403)
\spline(9200,2433)
(9020,2433)(8930,2283)
	(9020,2133)(10220,2133)(10310,2283)
	(10220,2433)(10040,2433)
\put(1970,2283){\makebox(0,0)[lb]{$+$}}
\put(4700,2253){\makebox(0,0)[lb]{$+$}}
\put(8030,2253){\makebox(0,0)[lb]{$+$}}
\put(2840,33){\makebox(0,0)[lb]{$+1$}}
\put(3440,33){\makebox(0,0)[lb]{$+1$}}
\put(4070,33){\makebox(0,0)[lb]{$+1$}}
\put(8570,33){\makebox(0,0)[lb]{$+1$}}
\put(9170,33){\makebox(0,0)[lb]{$+1$}}
\put(9770,33){\makebox(0,0)[lb]{$+1$}}
\put(5570,33){\makebox(0,0)[lb]{$+1$}}
\put(6170,33){\makebox(0,0)[lb]{$+1$}}
\put(6770,33){\makebox(0,0)[lb]{$+1$}}
\path(5630,11853)(5630,12783)(5780,12783)(5780,11853)
\path(5630,11583)(5630,11193)
\path(5780,11553)(5780,11223)
\path(6230,11853)(6230,12783)(6380,12783)(6380,11853)
\path(6830,11283)(6830,12783)(6980,12783)(6980,11283)
\path(6230,10953)(6230,10203)
\path(6380,10983)(6380,10233)
\path(6830,10983)(6830,10203)(6830,10203)
\path(6980,10983)(6980,10173)
\path(5630,10953)(5630,9453)(5780,9453)(5780,10983)
\path(6230,9963)(6230,9453)(6380,9453)(6380,9963)
\path(6830,9963)(6830,9453)(6980,9453)(6980,9963)
\path(4130,12753)(4130,9483)(4280,9483)
	(4280,12783)(4130,12783)
\path(3530,12753)(3530,9483)(3680,9483)
	(3680,12783)(3530,12783)
\path(2930,12753)(2930,9483)(3080,9483)
	(3080,12783)(2930,12783)
\path(8630,12723)(8630,9453)(8780,9453)
	(8780,12753)(8630,12753)
\path(9200,12723)(9200,9453)(9350,9453)
	(9350,12753)(9200,12753)
\path(9830,12723)(9830,9453)(9980,9453)
	(9980,12753)(9830,12753)
\path(5870,12033)(6170,12033)
\spline(5540,11523)
(5300,11523)(5210,11373)
	(5240,11193)(5330,11133)(7310,11133)
	(7400,11223)(7430,11343)(7400,11463)
	(7280,11523)(7190,11523)
\path(6470,10473)(6740,10473)
\path(5900,11493)(6710,11493)
\path(6230,11673)(6230,11553)
\path(6380,11673)(6380,11553)
\path(6230,11403)(6230,11253)
\path(6380,11403)(6380,11253)
\spline(6170,10383)
(5990,10383)(5900,10233)
	(5990,10083)(7190,10083)(7280,10233)
	(7190,10383)(7010,10383)
\spline(5600,12033)
(5420,12033)(5330,11883)
	(5420,11733)(6620,11733)(6710,11883)
	(6620,12033)(6440,12033)
\put(4760,11013){\makebox(0,0)[lb]{$-$}}
\put(8000,11013){\makebox(0,0)[lb]{$-$}}
\put(2870,9033){\makebox(0,0)[lb]{$0$}}
\put(3500,9033){\makebox(0,0)[lb]{$0$}}
\put(4100,9033){\makebox(0,0)[lb]{$0$}}
\put(5570,9063){\makebox(0,0)[lb]{$+1$}}
\put(6140,9063){\makebox(0,0)[lb]{$+1$}}
\put(6740,9063){\makebox(0,0)[lb]{$+1$}}
\put(8570,9033){\makebox(0,0)[lb]{$+1$}}
\put(9170,9033){\makebox(0,0)[lb]{$+1$}}
\put(9770,9033){\makebox(0,0)[lb]{$+1$}}
\put(2015,11058){\makebox(0,0)[lb]{$=$}}
\path(215,11283)(215,12783)(365,12783)(365,11283)
\path(1415,11283)(1415,12783)(1565,12783)(1565,11283)
\path(815,11283)(815,12783)(965,12783)(965,11283)
\path(215,10908)(215,9483)(365,9483)(365,10908)
\path(815,10908)(815,9483)(965,9483)(965,10908)
\path(1415,10908)(1415,9483)(1565,9483)(1565,10908)
\path(440,11358)(740,11358)
\path(1040,11358)(1340,11358)
\spline(215,11358)
(65,11358)(-10,11208)
	(65,10983)(1790,10983)(1865,11208)
	(1790,11358)(1640,11358)
\put(1790,10608){\makebox(0,0)[lb]{$+1$}}
\put(140,9033){\makebox(0,0)[lb]{$0$}}
\put(740,9033){\makebox(0,0)[lb]{$0$}}
\put(1340,9033){\makebox(0,0)[lb]{$0$}}
\end{picture} }
        	\hspace{-1.9mm}
	\end{array}
 $$   
 \caption{A figure showing an identity in $\cal M$. The links shown
are algebraically split. 
For a proof, see [GL1]. The numbers in the bottom of each 
band represent twists, with the conventions of figure \ref{f.bands}.
The framings in all the horizontal components is $+1$. }\lbl{3bands}
\end{figure}
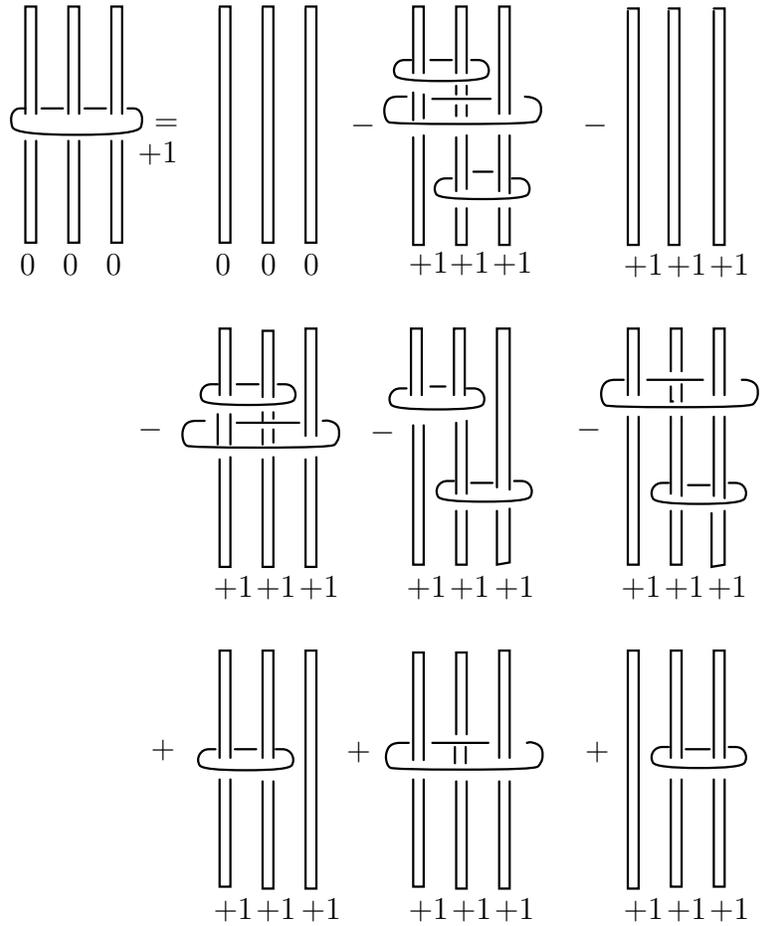

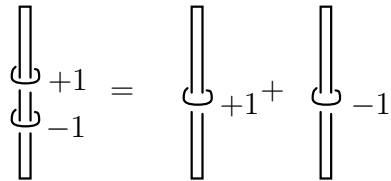
\begin{figure}[htpb]
$$\printname{plusminus}
	\setlength{\unitlength}{0.03\standardunitlength}
	\begin{array}{c}  \hspace{-1.7mm}
        	\raisebox{-8pt}{\begingroup\makeatletter\ifx\SetFigFont\undefined
% extract first six characters in \fmtname
\def\x#1#2#3#4#5#6#7\relax{\def\x{#1#2#3#4#5#6}}%
\expandafter\x\fmtname xxxxxx\relax \def\y{splain}%
\ifx\x\y   % LaTeX or SliTeX?
\gdef\SetFigFont#1#2#3{%
  \ifnum #1<17\tiny\else \ifnum #1<20\small\else
  \ifnum #1<24\normalsize\else \ifnum #1<29\large\else
  \ifnum #1<34\Large\else \ifnum #1<41\LARGE\else
     \huge\fi\fi\fi\fi\fi\fi
  \csname #3\endcsname}%
\else
\gdef\SetFigFont#1#2#3{\begingroup
  \count@#1\relax \ifnum 25<\count@\count@25\fi
  \def\x{\endgroup\@setsize\SetFigFont{#2pt}}%
  \expandafter\x
    \csname \romannumeral\the\count@ pt\expandafter\endcsname
    \csname @\romannumeral\the\count@ pt\endcsname
  \csname #3\endcsname}%
\fi
\fi\endgroup
\begin{picture}(5132,2439)(0,-10)
\thicklines
\path(128,1422)(128,2412)(278,2412)(278,1422)
\path(128,642)(128,12)(278,12)(278,642)
\path(128,1242)(128,822)
\path(278,1212)(278,822)
\path(2528,1212)(2528,2412)(2678,2412)(2678,1212)
\path(2528,942)(2528,12)(2678,12)(2678,912)
\path(4328,1212)(4328,2412)(4478,2412)(4478,1212)
\path(4328,942)(4328,12)(4478,12)(4478,912)
\spline(68,1572)
(8,1512)(8,1332)
	(398,1332)(398,1512)(338,1542)
\spline(68,972)
(8,912)(8,732)
	(398,732)(398,912)(338,942)
\spline(2468,1272)
(2408,1212)(2408,1032)
	(2798,1032)(2798,1212)(2738,1242)
\spline(4268,1272)
(4208,1212)(4208,1032)
	(4598,1032)(4598,1212)(4538,1242)
\put(518,1242){\makebox(0,0)[lb]{$+1$}}
\put(488,582){\makebox(0,0)[lb]{$-1$}}
\put(4748,912){\makebox(0,0)[lb]{$-1$}}
\put(2918,912){\makebox(0,0)[lb]{$+1$}}
\put(3488,1122){\makebox(0,0)[lb]{$+$}}
\put(1388,1122){\makebox(0,0)[lb]{$=$}}
\end{picture} }
        	\hspace{-1.9mm}
	\end{array}
 $$
 \caption{A figure showing an identity in $\cal M$. The links shown
are algebraically split. 
For a proof, see [GL1]. This figure is used to change
the $+1$ framing of an unknot to a $-1$ an vice versa. 
}\lbl{plusminus}
\end{figure}

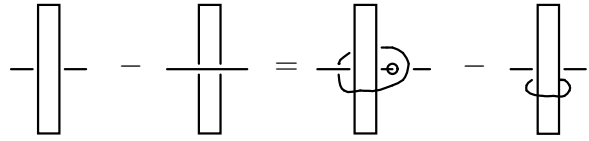
\begin{figure}[htpb]
$$\printname{Kauf}
	\setlength{\unitlength}{0.03\standardunitlength}
	\begin{array}{c}  \hspace{-1.7mm}
        	\raisebox{-8pt}{\begingroup\makeatletter\ifx\SetFigFont\undefined
% extract first six characters in \fmtname
\def\x#1#2#3#4#5#6#7\relax{\def\x{#1#2#3#4#5#6}}%
\expandafter\x\fmtname xxxxxx\relax \def\y{splain}%
\ifx\x\y   % LaTeX or SliTeX?
\gdef\SetFigFont#1#2#3{%
  \ifnum #1<17\tiny\else \ifnum #1<20\small\else
  \ifnum #1<24\normalsize\else \ifnum #1<29\large\else
  \ifnum #1<34\Large\else \ifnum #1<41\LARGE\else
     \huge\fi\fi\fi\fi\fi\fi
  \csname #3\endcsname}%
\else
\gdef\SetFigFont#1#2#3{\begingroup
  \count@#1\relax \ifnum 25<\count@\count@25\fi
  \def\x{\endgroup\@setsize\SetFigFont{#2pt}}%
  \expandafter\x
    \csname \romannumeral\the\count@ pt\expandafter\endcsname
    \csname @\romannumeral\the\count@ pt\endcsname
  \csname #3\endcsname}%
\fi
\fi\endgroup
\begin{picture}(8050,1840)(0,-10)
\thicklines
\put(5337,913){\ellipse{150}{150}}
\path(2187,913)(3312,913)
\path(2637,988)(2637,1813)(2937,1813)(2937,988)
\path(2637,838)(2637,13)(2937,13)(2937,838)
\path(4812,1813)(4812,13)(5112,13)
	(5112,1813)(4812,1813)
\path(5187,913)(5337,913)
\path(5637,913)(5862,913)
\path(4287,913)(4737,913)
\path(7363,1812)(7363,12)(7663,12)
	(7663,1812)(7363,1812)
\path(6988,912)(7288,912)
\path(7738,912)(8038,912)
\path(387,1813)(387,13)(687,13)
	(687,1813)(387,1813)
\path(12,913)(312,913)
\path(762,913)(1062,913)
\path(4737,1138)	(4694.052,1109.407)
	(4658.896,1084.505)
	(4610.205,1044.017)
	(4587.000,988.000)

\path(4587,988)	(4587.000,988.000)

\path(4587,838)	(4590.672,791.199)
	(4595.534,751.053)
	(4609.706,688.086)
	(4631.276,643.826)
	(4662.000,613.000)

\path(4662,613)	(4715.484,594.594)
	(4777.976,596.639)
	(4838.730,606.865)
	(4887.000,613.000)

\path(4887,613)	(4936.445,611.809)
	(4999.159,609.824)
	(5062.043,609.427)
	(5112.000,613.000)

\path(5112,613)	(5163.293,625.557)
	(5226.289,646.045)
	(5288.390,668.761)
	(5337.000,688.000)

\path(5337,688)	(5371.252,701.119)
	(5414.205,718.202)
	(5487.000,763.000)

\path(5487,763)	(5517.037,809.569)
	(5541.915,869.695)
	(5558.086,932.723)
	(5562.000,988.000)

\path(5562,988)	(5553.196,1026.199)
	(5534.681,1068.531)
	(5487.000,1138.000)

\path(5487,1138)	(5415.458,1183.382)
	(5372.866,1201.610)
	(5337.000,1213.000)

\path(5337,1213)	(5299.271,1214.568)
	(5262.000,1213.000)

\path(5262,1213)	(5235.986,1213.000)
	(5187.000,1213.000)

\path(7287,763)	(7234.564,719.271)
	(7212.000,688.000)

\path(7212,688)	(7203.645,651.059)
	(7212.000,613.000)

\path(7212,613)	(7258.824,575.616)
	(7321.121,554.916)
	(7385.108,544.509)
	(7437.000,538.000)

\path(7437,538)	(7470.305,535.621)
	(7512.229,535.885)
	(7587.000,538.000)

\path(7587,538)	(7660.800,532.195)
	(7702.602,531.469)
	(7737.000,538.000)

\path(7737,538)	(7780.305,569.695)
	(7812.000,613.000)

\path(7812,613)	(7817.805,650.500)
	(7812.000,688.000)

\path(7812,688)	(7780.305,731.305)
	(7737.000,763.000)

\path(7737,763)	(7692.491,765.902)
	(7662.000,763.000)

\path(7662,763)	(7678.406,763.000)
	(7737.000,763.000)

\put(1512,838){\makebox(0,0)[lb]{$-$}}
\put(3687,838){\makebox(0,0)[lb]{$=$}}
\put(6312,838){\makebox(0,0)[lb]{$-$}}
\end{picture} }
        	\hspace{-1.9mm}
	\end{array}
 $$
\caption{Another special case of equation \eqref{eq.fundamental}. 
We
assume that all the arcs in the band lie in the same link component. 
Then
the link shown in the third picture from the left contains
 a unit Seifert-framed $1$-pair blink
(indeed, it bounds a Seifert surface obtained by tubing the 
obvious
disc along the strands contained in the band). The fourth picture 
from 
the left contains a 
$1$-component boundary link in the complement of the rest of the 
link, and 
thus (using remark \ref{rem.bbl}) a $1$-pair blink.}\lbl{Kauf}
\end{figure}

\subsection{A $4$-term relation on $\M$ }
\lbl{sub.4T}

In this and in the  two section 
 we introduce three versions of  a $4$-term relation on $\M$, which 
will be
used crucially in the proof of theorem \ref{thm.as2b}.

We begin by noting that there is a map from $+1$-framed knots  in 
$S^3$
to $\M$, defined by $K \to [K, S^3, +1]= S^3 -S^3_{K,+1}$.

originally introduced   in \cite{Ga}. Dually
this map induces a map from invariants of \ihs s to invariants of 
knots
in $S^3$. Furthermore it is trivial to show that \fti s of \ihs s map
to \fti s of knots. Actually much more is known, namely that 
$AS$-type
$3m$ of \ihs s map to type $2m$ invariants of knots, see \cite{Ha} 
and
\cite{GL2}.
The vector space on the set of knots satisfies a basic relation, the
$4$-term relation. This relation   plays an important role in
studying \fti s of knots. According to the above defined map, we 
have a $4$-term relation satisfied on the image of this map in $\M$. 
Strangely enough, the above 
relation has not been explicitly introduced or noticed before. In the 
present
paper we describe the $4$-term relation on $\M$ and use it in a 
crucial
way in the proof of theorem \ref{thm.as2b}.
Furthermore, in the section \ref{sub.que} we pose the question of a 
possible
relation between the $4$-term relation and the more well known 
$AS$ and
$IHX$ relations on $\M$.
With these preliminaries and motivation in mind, we begin to 
describe the 
$4$-term relation on $\M$.

Let $C \cup K$ be a two component sublink of an \ASA\ link $L$ in 
$S^3$,
such that:
\begin{itemize}
\item
     The intersection of $K$ with a $3$-ball $B$ in $S^3$ consists of 
$3$
arcs shown in figure \ref{3arcs}.
\item
     $C$ is an unknot that bounds a disc $D$ which  lies in the interior 
of the ball $B$. The disc $D$ intersects $K$ in two points, $b,c$.
See figure \ref{3arcs}.
\end{itemize}

\begin{figure}[htpb]
$$\printname{3arcs}
	\setlength{\unitlength}{0.04\standardunitlength}
	\begin{array}{c}  \hspace{-1.7mm}
        	\raisebox{-8pt}{\begingroup\makeatletter\ifx\SetFigFont\undefined
% extract first six characters in \fmtname
\def\x#1#2#3#4#5#6#7\relax{\def\x{#1#2#3#4#5#6}}%
\expandafter\x\fmtname xxxxxx\relax \def\y{splain}%
\ifx\x\y   % LaTeX or SliTeX?
\gdef\SetFigFont#1#2#3{%
  \ifnum #1<17\tiny\else \ifnum #1<20\small\else
  \ifnum #1<24\normalsize\else \ifnum #1<29\large\else
  \ifnum #1<34\Large\else \ifnum #1<41\LARGE\else
     \huge\fi\fi\fi\fi\fi\fi
  \csname #3\endcsname}%
\else
\gdef\SetFigFont#1#2#3{\begingroup
  \count@#1\relax \ifnum 25<\count@\count@25\fi
  \def\x{\endgroup\@setsize\SetFigFont{#2pt}}%
  \expandafter\x
    \csname \romannumeral\the\count@ pt\expandafter\endcsname
    \csname @\romannumeral\the\count@ pt\endcsname
  \csname #3\endcsname}%
\fi
\fi\endgroup
\begin{picture}(2562,2139)(0,-10)
\thicklines
\put(1050,1212){\whiten\ellipse{76}{76}}
\put(1050,1212){\ellipse{76}{76}}
\put(1350,912){\whiten\ellipse{76}{76}}
\put(1350,912){\ellipse{76}{76}}
\put(375,912){\whiten\ellipse{76}{76}}
\put(375,912){\ellipse{76}{76}}
\path(150,1212)(825,1212)
\path(1350,12)(1350,687)
\path(1050,1212)(1200,1212)
\path(1500,1212)(2550,1212)
\path(2430.000,1182.000)(2550.000,1212.000)(2430.000,1242.000)
\path(1350,912)(1350,2112)
\path(1380.000,1992.000)(1350.000,2112.000)(1320.000,1992.000)
\path(1800,1287)	(1800.000,1343.482)
	(1800.000,1384.965)
	(1800.000,1437.000)

\path(1800,1437)	(1804.582,1473.495)
	(1800.000,1512.000)

\path(1800,1512)	(1728.382,1561.039)
	(1684.929,1575.825)
	(1650.000,1587.000)

\path(1650,1587)	(1585.471,1609.177)
	(1545.538,1621.625)
	(1503.188,1633.796)
	(1460.438,1644.806)
	(1419.310,1653.767)
	(1350.000,1662.000)

\path(1350,1662)	(1280.678,1656.915)
	(1239.541,1650.707)
	(1196.783,1642.189)
	(1154.426,1631.459)
	(1114.495,1618.619)
	(1050.000,1587.000)

\path(1050,1587)	(1005.757,1542.390)
	(963.416,1480.226)
	(926.867,1415.200)
	(900.000,1362.000)

\path(900,1362)	(879.309,1313.664)
	(856.755,1251.533)
	(837.073,1188.384)
	(825.000,1137.000)

\path(825,1137)	(818.161,1087.718)
	(813.116,1025.959)
	(814.013,963.470)
	(825.000,912.000)

\path(825,912)	(851.519,870.070)
	(891.799,827.325)
	(936.180,789.417)
	(975.000,762.000)

\path(975,762)	(1022.260,737.524)
	(1083.079,713.160)
	(1146.108,694.466)
	(1200.000,687.000)

\path(1200,687)	(1253.802,695.309)
	(1316.681,715.406)
	(1377.470,740.051)
	(1425.000,762.000)

\path(1425,762)	(1477.223,791.542)
	(1540.620,832.549)
	(1602.458,875.781)
	(1650.000,912.000)

\path(1650,912)	(1689.127,947.734)
	(1725.000,987.000)

\path(1725,987)	(1754.366,1037.790)
	(1774.271,1079.585)
	(1800.000,1137.000)

\put(1725,612){\makebox(0,0)[lb]{$c$}}
\put(600,1512){\makebox(0,0)[lb]{$b$}}
\put(0,537){\makebox(0,0)[lb]{$a$}}
\end{picture} }
        	\hspace{-1.9mm}
	\end{array}
 $$
\caption{We show the intersection of a $3$-ball $B$ with a knot 
$K$,
which consists of $3$ arcs (two of them are shown, the third is 
perpendicular
to the page pointing towards you). Shown also are $3$ points $a,b,c$ 
on the
knot $K$, as well as an unknot that bounds a disc which intersects 
$K$ in
$b,c$. }\lbl{3arcs}
\end{figure}
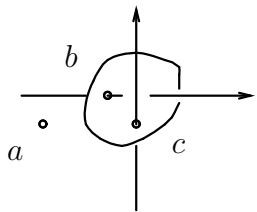

Choose $4$ disjoint discs $D_j$ (for $j=s,n,e,w$) ($s,n,e,w$ stands 
for
south, north, east and west) in the ball
 that intersect $K$ in the $4$ tuples of points $(j,a)$, as shown in 
figure
\ref{snew} with the abbreviations of figure \ref{abrev}. Let $C_j$ 
(for
$j=s,n,e,w$) denote the boundary of $D_j$ with framing $-1$.
Note that $C_j \cup L$ is an \ASA\ link in $S^3$. Let us momentarily
abbreviate the elements $[S^3, C_j \cup L, -1 \cup f]$ by $[C_j \cup 
L]$,
where $f$ is the framing of the \ASA\ link $L$. We can now state 
the 
following proposition:

\begin{proposition}
\lbl{prop.4T1}
With the above notation, we have the following relation on $\M$:
\begin{equation}
\lbl{2TM}
[C_s \cup L] - [C_n \cup L]=[C_e \cup L] -[C_w \cup L]
\end{equation}
We call the above relation the (first version of the)
$4$-term relation on $\M$.
\end{proposition}

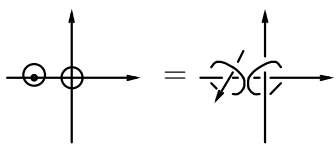
\begin{figure}[htpb]
$$\printname{abrev}
	\setlength{\unitlength}{0.03\standardunitlength}
	\begin{array}{c}  \hspace{-1.7mm}
        	\raisebox{-8pt}{\begingroup\makeatletter\ifx\SetFigFont\undefined
% extract first six characters in \fmtname
\def\x#1#2#3#4#5#6#7\relax{\def\x{#1#2#3#4#5#6}}%
\expandafter\x\fmtname xxxxxx\relax \def\y{splain}%
\ifx\x\y   % LaTeX or SliTeX?
\gdef\SetFigFont#1#2#3{%
  \ifnum #1<17\tiny\else \ifnum #1<20\small\else
  \ifnum #1<24\normalsize\else \ifnum #1<29\large\else
  \ifnum #1<34\Large\else \ifnum #1<41\LARGE\else
     \huge\fi\fi\fi\fi\fi\fi
  \csname #3\endcsname}%
\else
\gdef\SetFigFont#1#2#3{\begingroup
  \count@#1\relax \ifnum 25<\count@\count@25\fi
  \def\x{\endgroup\@setsize\SetFigFont{#2pt}}%
  \expandafter\x
    \csname \romannumeral\the\count@ pt\expandafter\endcsname
    \csname @\romannumeral\the\count@ pt\endcsname
  \csname #3\endcsname}%
\fi
\fi\endgroup
\begin{picture}(4524,1839)(0,-10)
\thicklines
\put(912,912){\ellipse{300}{300}}
\put(387,950){\ellipse{310}{310}}
\put(387,912){\whiten\ellipse{76}{76}}
\put(387,912){\ellipse{76}{76}}
\path(912,12)(912,1812)
\path(942.000,1692.000)(912.000,1812.000)(882.000,1692.000)
\path(12,912)(1812,912)
\path(1692.000,882.000)(1812.000,912.000)(1692.000,942.000)
\path(2712,912)(3012,912)
\path(3312,1287)(3237,1137)
\path(3162,912)(3237,912)
\path(3462,912)(3537,912)
\path(3837,837)(3687,687)
\path(2862,837)(2937,762)
\path(3612,12)(3612,987)
\path(3687,912)(4512,912)
\path(4392.000,882.000)(4512.000,912.000)(4392.000,942.000)
\path(3612,1212)(3612,1812)
\path(3642.000,1692.000)(3612.000,1812.000)(3582.000,1692.000)
\path(3162,987)(2937,612)
\path(2973.015,730.334)(2937.000,612.000)(3024.464,699.464)
\path(3837,1062)	(3796.406,1116.349)
	(3762.000,1137.000)

\path(3762,1137)	(3703.539,1136.659)
	(3641.017,1117.174)
	(3582.737,1088.851)
	(3537.000,1062.000)

\path(3537,1062)	(3498.180,1034.583)
	(3453.799,996.675)
	(3413.519,953.930)
	(3387.000,912.000)

\path(3387,912)	(3378.293,878.052)
	(3375.390,837.000)
	(3378.293,795.948)
	(3387.000,762.000)

\path(3387,762)	(3418.695,718.695)
	(3462.000,687.000)

\path(3462,687)	(3488.914,682.650)
	(3537.000,687.000)

\path(2862,1062)	(2902.594,1116.349)
	(2937.000,1137.000)

\path(2937,1137)	(2995.461,1136.659)
	(3057.983,1117.174)
	(3116.263,1088.851)
	(3162.000,1062.000)

\path(3162,1062)	(3200.820,1034.583)
	(3245.201,996.675)
	(3285.481,953.930)
	(3312.000,912.000)

\path(3312,912)	(3320.707,878.052)
	(3323.610,837.000)
	(3320.707,795.948)
	(3312.000,762.000)

\path(3312,762)	(3280.305,718.695)
	(3237.000,687.000)

\path(3237,687)	(3210.086,682.650)
	(3162.000,687.000)

\put(2187,837){\makebox(0,0)[lb]{$=$}}
\end{picture} }
        	\hspace{-1.9mm}
	\end{array}
 $$
\caption{Some abbreviation conventions for drawing the next 
figure.}\lbl{abrev}
\end{figure}

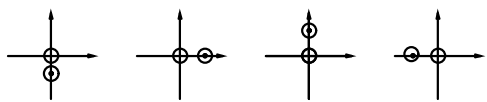
\begin{figure}[htpb]
$$\printname{snew}
	\setlength{\unitlength}{0.02\standardunitlength}
	\begin{array}{c}  \hspace{-1.7mm}
        	\raisebox{-8pt}{\begingroup\makeatletter\ifx\SetFigFont\undefined
% extract first six characters in \fmtname
\def\x#1#2#3#4#5#6#7\relax{\def\x{#1#2#3#4#5#6}}%
\expandafter\x\fmtname xxxxxx\relax \def\y{splain}%
\ifx\x\y   % LaTeX or SliTeX?
\gdef\SetFigFont#1#2#3{%
  \ifnum #1<17\tiny\else \ifnum #1<20\small\else
  \ifnum #1<24\normalsize\else \ifnum #1<29\large\else
  \ifnum #1<34\Large\else \ifnum #1<41\LARGE\else
     \huge\fi\fi\fi\fi\fi\fi
  \csname #3\endcsname}%
\else
\gdef\SetFigFont#1#2#3{\begingroup
  \count@#1\relax \ifnum 25<\count@\count@25\fi
  \def\x{\endgroup\@setsize\SetFigFont{#2pt}}%
  \expandafter\x
    \csname \romannumeral\the\count@ pt\expandafter\endcsname
    \csname @\romannumeral\the\count@ pt\endcsname
  \csname #3\endcsname}%
\fi
\fi\endgroup
\begin{picture}(9925,1840)(0,-10)
\thicklines
\put(912,538){\whiten\ellipse{76}{76}}
\put(912,538){\ellipse{76}{76}}
\put(912,913){\ellipse{300}{300}}
\put(912,538){\ellipse{310}{310}}
\put(3612,913){\ellipse{300}{300}}
\put(9013,912){\ellipse{300}{300}}
\put(4137,913){\ellipse{300}{300}}
\put(6312,1438){\ellipse{300}{300}}
\put(8456,944){\ellipse{300}{300}}
\put(4137,913){\whiten\ellipse{76}{76}}
\put(4137,913){\ellipse{76}{76}}
\put(6312,1438){\whiten\ellipse{76}{76}}
\put(6312,1438){\ellipse{76}{76}}
\put(8487,913){\whiten\ellipse{76}{76}}
\put(8487,913){\ellipse{76}{76}}
\put(6313,912){\ellipse{300}{300}}
\put(6313,912){\ellipse{300}{300}}
\path(912,13)(912,1813)
\path(942.000,1693.000)(912.000,1813.000)(882.000,1693.000)
\path(12,913)(1812,913)
\path(1692.000,883.000)(1812.000,913.000)(1692.000,943.000)
\path(2712,913)(4512,913)
\path(4392.000,883.000)(4512.000,913.000)(4392.000,943.000)
\path(3612,13)(3612,1813)
\path(3642.000,1693.000)(3612.000,1813.000)(3582.000,1693.000)
\path(9013,12)(9013,1812)
\path(9043.000,1692.000)(9013.000,1812.000)(8983.000,1692.000)
\path(8113,912)(9913,912)
\path(9793.000,882.000)(9913.000,912.000)(9793.000,942.000)
\path(5413,912)(7213,912)
\path(6313,12)(6313,1812)
\path(6313,12)(6313,1812)
\path(6343.000,1692.000)(6313.000,1812.000)(6283.000,1692.000)
\path(5413,912)(7213,912)
\path(7093.000,882.000)(7213.000,912.000)(7093.000,942.000)
\end{picture} }
        	\hspace{-1.9mm}
	\end{array}
 $$
\caption{The terms $s,e,n,w$ appearing in the $4$-term relation on 
$\cal M$
with the drawing conventions of the previous figure. }\lbl{snew}
\end{figure}

\begin{pf}
The proof of the $4$-term relation on $\M$ is the same as the proof 
of the
$4$-term relation on the space of knots, see \cite{B-N1}. In both 
proofs, we move the arc $a$ of figure \ref{3arcs} from the $SW$ 
quarter,
to the $NE$ quarter in two ways: either by passing through the $NW$
quarter, or by passing through the $SE$ quarter. Using figure 
\ref{crosschange},  the difference in the
first (respectively, second) way equals to the left (respectively,
right) hand side of equation \eqref{2TM}, thus proving the 
proposition.
\end{pf}

\begin{figure}[htpb]
$$\printname{crosschange}
	\setlength{\unitlength}{0.02\standardunitlength}
	\begin{array}{c}  \hspace{-1.7mm}
        	\raisebox{-8pt}{\begingroup\makeatletter\ifx\SetFigFont\undefined
% extract first six characters in \fmtname
\def\x#1#2#3#4#5#6#7\relax{\def\x{#1#2#3#4#5#6}}%
\expandafter\x\fmtname xxxxxx\relax \def\y{splain}%
\ifx\x\y   % LaTeX or SliTeX?
\gdef\SetFigFont#1#2#3{%
  \ifnum #1<17\tiny\else \ifnum #1<20\small\else
  \ifnum #1<24\normalsize\else \ifnum #1<29\large\else
  \ifnum #1<34\Large\else \ifnum #1<41\LARGE\else
     \huge\fi\fi\fi\fi\fi\fi
  \csname #3\endcsname}%
\else
\gdef\SetFigFont#1#2#3{\begingroup
  \count@#1\relax \ifnum 25<\count@\count@25\fi
  \def\x{\endgroup\@setsize\SetFigFont{#2pt}}%
  \expandafter\x
    \csname \romannumeral\the\count@ pt\expandafter\endcsname
    \csname @\romannumeral\the\count@ pt\endcsname
  \csname #3\endcsname}%
\fi
\fi\endgroup
\begin{picture}(6324,939)(0,-10)
\thicklines
\put(462,462){\whiten\ellipse{76}{76}}
\put(462,462){\ellipse{76}{76}}
\path(12,12)(912,912)
\path(12,12)(912,912)
\path(848.360,805.934)(912.000,912.000)(805.934,848.360)
\path(912,12)(12,912)
\path(912,12)(12,912)
\path(118.066,848.360)(12.000,912.000)(75.640,805.934)
\path(2712,12)(2337,387)
\path(2712,12)(2337,387)
\path(3612,12)(3987,387)
\path(3612,12)(3987,387)
\path(5412,12)(5562,162)
\path(6312,12)(5862,387)
\path(5787,462)(5712,537)
\path(2187,537)(1812,912)
\path(1918.066,848.360)(1812.000,912.000)(1875.640,805.934)
\path(1812,12)(2712,912)
\path(2648.360,805.934)(2712.000,912.000)(2605.934,848.360)
\path(4512,12)(3612,912)
\path(3718.066,848.360)(3612.000,912.000)(3675.640,805.934)
\path(4062,537)(4512,912)
\path(4439.019,812.131)(4512.000,912.000)(4400.608,858.225)
\path(5787,387)(6312,912)
\path(6248.360,805.934)(6312.000,912.000)(6205.934,848.360)
\path(5637,687)(5412,912)
\path(5518.066,848.360)(5412.000,912.000)(5475.640,805.934)
\path(6012,462)(6012,387)
\path(5937,612)	(5893.271,664.436)
	(5862.000,687.000)

\path(5862,687)	(5787.818,701.970)
	(5746.868,699.488)
	(5712.000,687.000)

\path(5712,687)	(5661.420,616.594)
	(5646.441,573.292)
	(5637.000,537.000)

\path(5637,537)	(5625.597,470.762)
	(5621.493,430.667)
	(5619.124,388.474)
	(5618.989,346.096)
	(5621.589,305.448)
	(5637.000,237.000)

\path(5637,237)	(5668.695,193.695)
	(5712.000,162.000)

\path(5712,162)	(5746.398,155.469)
	(5788.200,156.195)
	(5862.000,162.000)

\path(5862,162)	(5898.176,154.268)
	(5937.000,162.000)

\path(5937,162)	(5960.205,218.018)
	(5954.404,258.505)
	(5937.000,312.000)

\put(1287,387){\makebox(0,0)[lb]{$=$}}
\put(3087,462){\makebox(0,0)[lb]{$-$}}
\put(4812,462){\makebox(0,0)[lb]{$=$}}
\end{picture} }
        	\hspace{-1.9mm}
	\end{array}
 $$
\caption{An identity useful for the $4$-term 
relation.}\lbl{crosschange}
\end{figure}
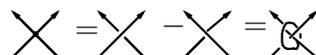

\begin{remark}
\lbl{rem.order}
The points $a,b,c,s,w,e,n$ (on the knot $K$) are displayed in the 
order
shown in figure \ref{order}.  
\end{remark}

\begin{figure}[htpb]
$$\printname{order}
	\setlength{\unitlength}{0.02\standardunitlength}
	\begin{array}{c}  \hspace{-1.7mm}
        	\raisebox{-8pt}{\begingroup\makeatletter\ifx\SetFigFont\undefined
% extract first six characters in \fmtname
\def\x#1#2#3#4#5#6#7\relax{\def\x{#1#2#3#4#5#6}}%
\expandafter\x\fmtname xxxxxx\relax \def\y{splain}%
\ifx\x\y   % LaTeX or SliTeX?
\gdef\SetFigFont#1#2#3{%
  \ifnum #1<17\tiny\else \ifnum #1<20\small\else
  \ifnum #1<24\normalsize\else \ifnum #1<29\large\else
  \ifnum #1<34\Large\else \ifnum #1<41\LARGE\else
     \huge\fi\fi\fi\fi\fi\fi
  \csname #3\endcsname}%
\else
\gdef\SetFigFont#1#2#3{\begingroup
  \count@#1\relax \ifnum 25<\count@\count@25\fi
  \def\x{\endgroup\@setsize\SetFigFont{#2pt}}%
  \expandafter\x
    \csname \romannumeral\the\count@ pt\expandafter\endcsname
    \csname @\romannumeral\the\count@ pt\endcsname
  \csname #3\endcsname}%
\fi
\fi\endgroup
\begin{picture}(2843,3030)(0,-10)
\thicklines
\put(1217,1533){\ellipse{2418}{2418}}
\put(842,2658){\whiten\ellipse{76}{76}}
\put(842,2658){\ellipse{76}{76}}
\put(392,2433){\whiten\ellipse{76}{76}}
\put(392,2433){\ellipse{76}{76}}
\put(1217,2733){\whiten\ellipse{76}{76}}
\put(1217,2733){\ellipse{76}{76}}
\put(1067,333){\whiten\ellipse{76}{76}}
\put(1067,333){\ellipse{76}{76}}
\put(1517,333){\whiten\ellipse{76}{76}}
\put(1517,333){\ellipse{76}{76}}
\put(2417,1608){\whiten\ellipse{76}{76}}
\put(2417,1608){\ellipse{76}{76}}
\put(542,558){\whiten\ellipse{76}{76}}
\put(542,558){\ellipse{76}{76}}
\put(2567,1533){\makebox(0,0)[lb]{$a$}}
\put(692,2808){\makebox(0,0)[lb]{$b$}}
\put(992,33){\makebox(0,0)[lb]{$c$}}
\put(242,2583){\makebox(0,0)[lb]{$e$}}
\put(1217,2883){\makebox(0,0)[lb]{$w$}}
\put(392,258){\makebox(0,0)[lb]{$s$}}
\put(1517,33){\makebox(0,0)[lb]{$n$}}
\end{picture} }
        	\hspace{-1.9mm}
	\end{array}
 $$
\caption{The order various points on a knot $K$.}\lbl{order}
\end{figure}
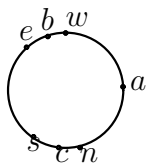

\begin{remark}
\lbl{rem.4T1}
If we represent the knot $K$ by a  circle and the knots $C, C_j$
(for $j=s,n,e,w$) by chords (that intersect the above mentioned 
circle
in two points each, namely the points of intersection $K \cap D_j$) 
then the $4$-term relation reads as in figure 
\ref{4Tv1}. The $4$-term relation will be
used, in the form of figure \ref{4Tv1}, in the proof of theorem 
\ref{thm.as2b}.
\end{remark}

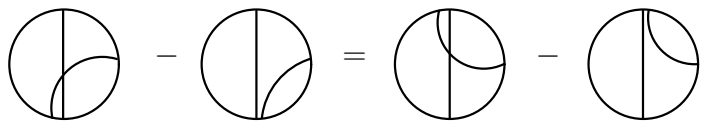
\begin{figure}[htpb]
$$\printname{4Tv1}
	\setlength{\unitlength}{0.03\standardunitlength}
	\begin{array}{c}  \hspace{-1.7mm}
        	\raisebox{-8pt}{\begingroup\makeatletter\ifx\SetFigFont\undefined
% extract first six characters in \fmtname
\def\x#1#2#3#4#5#6#7\relax{\def\x{#1#2#3#4#5#6}}%
\expandafter\x\fmtname xxxxxx\relax \def\y{splain}%
\ifx\x\y   % LaTeX or SliTeX?
\gdef\SetFigFont#1#2#3{%
  \ifnum #1<17\tiny\else \ifnum #1<20\small\else
  \ifnum #1<24\normalsize\else \ifnum #1<29\large\else
  \ifnum #1<34\Large\else \ifnum #1<41\LARGE\else
     \huge\fi\fi\fi\fi\fi\fi
  \csname #3\endcsname}%
\else
\gdef\SetFigFont#1#2#3{\begingroup
  \count@#1\relax \ifnum 25<\count@\count@25\fi
  \def\x{\endgroup\@setsize\SetFigFont{#2pt}}%
  \expandafter\x
    \csname \romannumeral\the\count@ pt\expandafter\endcsname
    \csname @\romannumeral\the\count@ pt\endcsname
  \csname #3\endcsname}%
\fi
\fi\endgroup
\begin{picture}(9635,1562)(0,-10)
\thicklines
\put(1309.500,166.500){\arc{1423.815}{2.9401}{5.0008}}
\put(4507.000,-81.000){\arc{1951.333}{3.2494}{4.4052}}
\put(6635.077,1356.692){\arc{1290.298}{1.1271}{3.4039}}
\put(9579.732,1423.709){\arc{1301.019}{1.5212}{3.2964}}
\put(773,773){\ellipse{1530}{1530}}
\put(3462,774){\ellipse{1530}{1530}}
\put(6162,774){\ellipse{1530}{1530}}
\put(8862,774){\ellipse{1530}{1530}}
\path(762,1524)(762,24)
\path(3462,1524)(3462,24)
\path(6162,1524)(6162,24)
\path(8862,1524)(8862,24)
\put(2037,774){\makebox(0,0)[lb]{$-$}}
\put(4662,774){\makebox(0,0)[lb]{$=$}}
\put(7362,774){\makebox(0,0)[lb]{$-$}}
\end{picture} }
        	\hspace{-1.9mm}
	\end{array}
 $$
\caption{The first version of the $4$-term relation on $\cal M$, 
with
the notation of remark \ref{rem.4T1}.}\lbl{4Tv1}
\end{figure}

\subsection{Two more versions of the $4$-term relation on $\M$}
\lbl{sub.4Tmore}

In this section we will 
introduce two more refined versions of the $4$-term
relation on $\M$.

Let $C \cup K_1 \cup K_2$ be a three
 component sublink of an \ASA\ link $L$ in 
$S^3$,
such that there is a three ball $B$ in $S^3$ with the properties:
\begin{itemize}
\item
     The intersection of $K_1$ with $B$  consists of $3$
arcs shown in figure \ref{3arcs2}.
\item
     The intersection of $K_2$ with $B$  consists of $2$
arcs shown in figure \ref{3arcs2}.
\item
     $C$ is an unknot that bounds a disc $D$ which  lies in the interior 
of the ball $B$. The disc $D$ intersects $K_1$ in $2$ points 
and intersects $K_2$ in two points. Furthermore, $D$ intersects no 
other 
component of $L$. See figure \ref{3arcs2}.
\end{itemize}

Let $D_j$ (for $j=s,n,e,w$) be $4$ discs as in section \ref{sub.4T}
(with figure \ref{3arcs2} replacing figure \ref{3arcs}). Let $C_j$
be the boundary of the disc $D_j$.

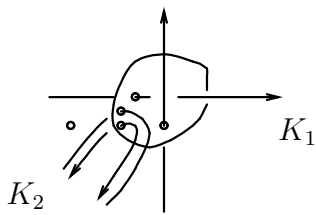
\begin{figure}[htpb]
$$\printname{3arcs2}
	\setlength{\unitlength}{0.04\standardunitlength}
	\begin{array}{c}  \hspace{-1.7mm}
        	\raisebox{-8pt}{\begingroup\makeatletter\ifx\SetFigFont\undefined
% extract first six characters in \fmtname
\def\x#1#2#3#4#5#6#7\relax{\def\x{#1#2#3#4#5#6}}%
\expandafter\x\fmtname xxxxxx\relax \def\y{splain}%
\ifx\x\y   % LaTeX or SliTeX?
\gdef\SetFigFont#1#2#3{%
  \ifnum #1<17\tiny\else \ifnum #1<20\small\else
  \ifnum #1<24\normalsize\else \ifnum #1<29\large\else
  \ifnum #1<34\Large\else \ifnum #1<41\LARGE\else
     \huge\fi\fi\fi\fi\fi\fi
  \csname #3\endcsname}%
\else
\gdef\SetFigFont#1#2#3{\begingroup
  \count@#1\relax \ifnum 25<\count@\count@25\fi
  \def\x{\endgroup\@setsize\SetFigFont{#2pt}}%
  \expandafter\x
    \csname \romannumeral\the\count@ pt\expandafter\endcsname
    \csname @\romannumeral\the\count@ pt\endcsname
  \csname #3\endcsname}%
\fi
\fi\endgroup
\begin{picture}(3357,2139)(0,-10)
\thicklines
\put(1350,1212){\whiten\ellipse{76}{76}}
\put(1350,1212){\ellipse{76}{76}}
\put(1650,912){\whiten\ellipse{76}{76}}
\put(1650,912){\ellipse{76}{76}}
\put(675,912){\whiten\ellipse{76}{76}}
\put(675,912){\ellipse{76}{76}}
\put(1200,1062){\whiten\ellipse{76}{76}}
\put(1200,1062){\ellipse{76}{76}}
\put(1200,912){\whiten\ellipse{76}{76}}
\put(1200,912){\ellipse{76}{76}}
\path(450,1212)(1125,1212)
\path(1650,12)(1650,687)
\path(1350,1212)(1500,1212)
\path(1800,1212)(2850,1212)
\path(2730.000,1182.000)(2850.000,1212.000)(2730.000,1242.000)
\path(1650,912)(1650,2112)
\path(1680.000,1992.000)(1650.000,2112.000)(1620.000,1992.000)
\path(2100,1287)	(2100.000,1343.482)
	(2100.000,1384.965)
	(2100.000,1437.000)

\path(2100,1437)	(2104.582,1473.495)
	(2100.000,1512.000)

\path(2100,1512)	(2028.382,1561.039)
	(1984.929,1575.825)
	(1950.000,1587.000)

\path(1950,1587)	(1885.471,1609.177)
	(1845.538,1621.625)
	(1803.188,1633.796)
	(1760.438,1644.806)
	(1719.310,1653.767)
	(1650.000,1662.000)

\path(1650,1662)	(1580.678,1656.915)
	(1539.541,1650.707)
	(1496.783,1642.189)
	(1454.426,1631.459)
	(1414.495,1618.619)
	(1350.000,1587.000)

\path(1350,1587)	(1305.757,1542.390)
	(1263.416,1480.226)
	(1226.867,1415.200)
	(1200.000,1362.000)

\path(1200,1362)	(1179.309,1313.664)
	(1156.755,1251.533)
	(1137.073,1188.384)
	(1125.000,1137.000)

\path(1125,1137)	(1118.161,1087.718)
	(1113.116,1025.959)
	(1114.013,963.470)
	(1125.000,912.000)

\path(1125,912)	(1151.519,870.070)
	(1191.799,827.325)
	(1236.180,789.417)
	(1275.000,762.000)

\path(1275,762)	(1322.260,737.524)
	(1383.079,713.160)
	(1446.108,694.466)
	(1500.000,687.000)

\path(1500,687)	(1553.802,695.309)
	(1616.681,715.406)
	(1677.470,740.051)
	(1725.000,762.000)

\path(1725,762)	(1777.223,791.542)
	(1840.620,832.549)
	(1902.458,875.781)
	(1950.000,912.000)

\path(1950,912)	(1989.127,947.734)
	(2025.000,987.000)

\path(2025,987)	(2054.366,1037.790)
	(2074.271,1079.585)
	(2100.000,1137.000)

\path(600,462)	(630.794,501.772)
	(659.470,538.726)
	(710.904,604.725)
	(755.183,661.098)
	(793.185,708.941)
	(853.873,783.436)
	(900.000,837.000)

\path(900,837)	(951.007,890.145)
	(992.748,931.351)
	(1050.000,987.000)

\path(1200,1062)	(1248.634,1064.348)
	(1275.000,1062.000)

\path(1275,1062)	(1310.866,1050.610)
	(1353.457,1032.382)
	(1425.000,987.000)

\path(1425,987)	(1473.165,917.963)
	(1491.741,875.684)
	(1500.000,837.000)

\path(1500,837)	(1496.667,792.283)
	(1484.622,743.005)
	(1466.019,691.138)
	(1443.011,638.655)
	(1417.752,587.527)
	(1392.395,539.727)
	(1369.093,497.227)
	(1350.000,462.000)

\path(1350,462)	(1319.863,410.322)
	(1278.922,346.864)
	(1236.021,284.724)
	(1200.000,237.000)

\path(1200,237)	(1149.750,182.917)
	(1107.820,141.946)
	(1050.000,87.000)

\path(1200,912)	(1253.495,929.404)
	(1293.983,935.205)
	(1350.000,912.000)

\path(1350,912)	(1368.797,876.197)
	(1369.189,834.964)
	(1350.000,762.000)

\path(1350,762)	(1334.215,719.362)
	(1310.140,671.200)
	(1280.156,619.672)
	(1246.646,566.932)
	(1211.993,515.139)
	(1178.580,466.448)
	(1148.788,423.017)
	(1125.000,387.000)

\path(1125,387)	(1103.478,353.597)
	(1073.648,308.479)
	(1031.993,246.371)
	(1005.634,207.298)
	(975.000,162.000)

\path(1017.423,278.191)(975.000,162.000)(1067.110,244.558)
\path(724.691,500.273)(675.000,387.000)(772.150,463.563)
\path(675,387)	(705.794,426.772)
	(734.470,463.726)
	(785.904,529.725)
	(830.183,586.098)
	(868.185,633.941)
	(928.873,708.436)
	(975.000,762.000)

\path(975,762)	(1000.504,788.569)
	(1050.000,837.000)

\put(2850,762){\makebox(0,0)[lb]{$K_1$}}
\put(0,87){\makebox(0,0)[lb]{$K_2$}}
\end{picture} }
        	\hspace{-1.9mm}
	\end{array}
 $$
\caption{Shown here are the $3$ arcs along which a ball $B$ 
intersects the 
knot $K_1$
and the $2$ arcs along which it intersects the knot $K_2$. 
Shown also with bold %$\circle$ 
are points of $K_1$ and $K_2$ where they intersect 
$D$.}\lbl{3arcs2}
\end{figure}

\begin{proposition}
\lbl{prop.4T2}
With the above notation, we have the following relation on $\M$:
\begin{equation}
\lbl{eq.4T2}
[C_s \cup L] - [C_n \cup L]=[C_e \cup L] -[C_w \cup L] + [E_1^{'}] - 
[E_2^{'}]
\end{equation}
where $E_1, E_2$ are links shown in figure \ref{E1E2}.
In the rest of the paper, the above 
relation will be called  the (second version of the)
$4$-term relation on $\M$.
\end{proposition}

\begin{figure}[htpb]
$$\printname{E1E2}
	\setlength{\unitlength}{0.03\standardunitlength}
	\begin{array}{c}  \hspace{-1.7mm}
        	\raisebox{-8pt}{\begingroup\makeatletter\ifx\SetFigFont\undefined
% extract first six characters in \fmtname
\def\x#1#2#3#4#5#6#7\relax{\def\x{#1#2#3#4#5#6}}%
\expandafter\x\fmtname xxxxxx\relax \def\y{splain}%
\ifx\x\y   % LaTeX or SliTeX?
\gdef\SetFigFont#1#2#3{%
  \ifnum #1<17\tiny\else \ifnum #1<20\small\else
  \ifnum #1<24\normalsize\else \ifnum #1<29\large\else
  \ifnum #1<34\Large\else \ifnum #1<41\LARGE\else
     \huge\fi\fi\fi\fi\fi\fi
  \csname #3\endcsname}%
\else
\gdef\SetFigFont#1#2#3{\begingroup
  \count@#1\relax \ifnum 25<\count@\count@25\fi
  \def\x{\endgroup\@setsize\SetFigFont{#2pt}}%
  \expandafter\x
    \csname \romannumeral\the\count@ pt\expandafter\endcsname
    \csname @\romannumeral\the\count@ pt\endcsname
  \csname #3\endcsname}%
\fi
\fi\endgroup
\begin{picture}(7708,2140)(0,-10)
\thicklines
\put(1350,1213){\whiten\ellipse{76}{76}}
\put(1350,1213){\ellipse{76}{76}}
\put(1650,913){\whiten\ellipse{76}{76}}
\put(1650,913){\ellipse{76}{76}}
\put(1200,1063){\whiten\ellipse{76}{76}}
\put(1200,1063){\ellipse{76}{76}}
\put(1200,913){\whiten\ellipse{76}{76}}
\put(1200,913){\ellipse{76}{76}}
\put(600,838){\whiten\ellipse{76}{76}}
\put(600,838){\ellipse{76}{76}}
\put(600,876){\ellipse{300}{226}}
\put(5701,1212){\whiten\ellipse{76}{76}}
\put(5701,1212){\ellipse{76}{76}}
\put(6001,912){\whiten\ellipse{76}{76}}
\put(6001,912){\ellipse{76}{76}}
\put(5551,1062){\whiten\ellipse{76}{76}}
\put(5551,1062){\ellipse{76}{76}}
\put(5551,912){\whiten\ellipse{76}{76}}
\put(5551,912){\ellipse{76}{76}}
\put(4951,837){\whiten\ellipse{76}{76}}
\put(4951,837){\ellipse{76}{76}}
\path(450,1213)(1125,1213)
\path(1650,13)(1650,688)
\path(1350,1213)(1500,1213)
\path(1800,1213)(2850,1213)
\path(2730.000,1183.000)(2850.000,1213.000)(2730.000,1243.000)
\path(1650,913)(1650,2113)
\path(1680.000,1993.000)(1650.000,2113.000)(1620.000,1993.000)
\path(750,388)(900,238)
\path(4801,1212)(5476,1212)
\path(6001,12)(6001,687)
\path(5701,1212)(5851,1212)
\path(6151,1212)(7201,1212)
\path(7081.000,1182.000)(7201.000,1212.000)(7081.000,1242.000)
\path(6001,912)(6001,2112)
\path(6031.000,1992.000)(6001.000,2112.000)(5971.000,1992.000)
\path(5101,387)(5251,237)
\path(2100,1288)	(2100.000,1344.482)
	(2100.000,1385.965)
	(2100.000,1438.000)

\path(2100,1438)	(2104.582,1474.495)
	(2100.000,1513.000)

\path(2100,1513)	(2028.382,1562.039)
	(1984.929,1576.825)
	(1950.000,1588.000)

\path(1950,1588)	(1885.471,1610.177)
	(1845.538,1622.625)
	(1803.188,1634.796)
	(1760.438,1645.806)
	(1719.310,1654.767)
	(1650.000,1663.000)

\path(1650,1663)	(1580.678,1657.915)
	(1539.541,1651.707)
	(1496.783,1643.189)
	(1454.426,1632.459)
	(1414.495,1619.619)
	(1350.000,1588.000)

\path(1350,1588)	(1305.757,1543.390)
	(1263.416,1481.226)
	(1226.867,1416.200)
	(1200.000,1363.000)

\path(1200,1363)	(1179.309,1314.664)
	(1156.755,1252.533)
	(1137.073,1189.384)
	(1125.000,1138.000)

\path(1125,1138)	(1118.161,1088.718)
	(1113.116,1026.959)
	(1114.013,964.470)
	(1125.000,913.000)

\path(1125,913)	(1151.519,871.070)
	(1191.799,828.325)
	(1236.180,790.417)
	(1275.000,763.000)

\path(1275,763)	(1322.260,738.524)
	(1383.079,714.160)
	(1446.108,695.466)
	(1500.000,688.000)

\path(1500,688)	(1553.802,696.309)
	(1616.681,716.406)
	(1677.470,741.051)
	(1725.000,763.000)

\path(1725,763)	(1777.223,792.542)
	(1840.620,833.549)
	(1902.458,876.781)
	(1950.000,913.000)

\path(1950,913)	(1989.127,948.734)
	(2025.000,988.000)

\path(2025,988)	(2054.366,1038.790)
	(2074.271,1080.585)
	(2100.000,1138.000)

\path(600,463)	(630.794,502.772)
	(659.470,539.726)
	(710.904,605.725)
	(755.183,662.098)
	(793.185,709.941)
	(853.873,784.436)
	(900.000,838.000)

\path(900,838)	(951.007,891.145)
	(992.748,932.351)
	(1050.000,988.000)

\path(1200,1063)	(1248.634,1065.348)
	(1275.000,1063.000)

\path(1275,1063)	(1310.866,1051.610)
	(1353.457,1033.382)
	(1425.000,988.000)

\path(1425,988)	(1473.165,918.963)
	(1491.741,876.684)
	(1500.000,838.000)

\path(1500,838)	(1496.667,793.283)
	(1484.622,744.005)
	(1466.019,692.138)
	(1443.011,639.655)
	(1417.752,588.527)
	(1392.395,540.727)
	(1369.093,498.227)
	(1350.000,463.000)

\path(1350,463)	(1319.863,411.322)
	(1278.922,347.864)
	(1236.021,285.724)
	(1200.000,238.000)

\path(1200,238)	(1149.750,183.917)
	(1107.820,142.946)
	(1050.000,88.000)

\path(1200,913)	(1253.495,930.404)
	(1293.983,936.205)
	(1350.000,913.000)

\path(1350,913)	(1368.797,877.197)
	(1369.189,835.964)
	(1350.000,763.000)

\path(1350,763)	(1334.215,720.362)
	(1310.140,672.200)
	(1280.156,620.672)
	(1246.646,567.932)
	(1211.993,516.139)
	(1178.580,467.448)
	(1148.788,424.017)
	(1125.000,388.000)

\path(1125,388)	(1103.478,354.597)
	(1073.648,309.479)
	(1031.993,247.371)
	(1005.634,208.298)
	(975.000,163.000)

\path(1017.423,279.191)(975.000,163.000)(1067.110,245.558)
\path(724.691,501.273)(675.000,388.000)(772.150,464.563)
\path(675,388)	(705.794,427.772)
	(734.470,464.726)
	(785.904,530.725)
	(830.183,587.098)
	(868.185,634.941)
	(928.873,709.436)
	(975.000,763.000)

\path(975,763)	(1000.504,789.569)
	(1050.000,838.000)

\path(525,1063)	(584.096,1070.619)
	(635.263,1076.061)
	(679.382,1079.327)
	(717.330,1080.415)
	(778.232,1076.061)
	(825.000,1063.000)

\path(825,1063)	(886.394,1021.198)
	(948.007,957.408)
	(1004.367,890.163)
	(1050.000,838.000)

\path(1050,838)	(1118.423,775.157)
	(1161.299,737.164)
	(1206.206,696.674)
	(1250.338,655.099)
	(1290.886,613.854)
	(1325.043,574.350)
	(1350.000,538.000)

\path(1350,538)	(1376.413,473.008)
	(1388.010,432.548)
	(1398.371,389.691)
	(1407.374,346.623)
	(1414.894,305.529)
	(1425.000,238.000)

\path(1425,238)	(1442.062,165.520)
	(1443.000,123.984)
	(1425.000,88.000)

\path(1425,88)	(1374.964,53.198)
	(1312.973,41.598)
	(1275.296,44.498)
	(1231.995,53.198)
	(1182.189,67.699)
	(1125.000,88.000)

\path(525,1063)	(463.338,1014.399)
	(420.128,975.543)
	(375.000,913.000)

\path(375,913)	(366.293,879.052)
	(363.390,838.000)
	(366.293,796.948)
	(375.000,763.000)

\path(375,763)	(402.867,722.072)
	(445.395,680.995)
	(490.225,643.421)
	(525.000,613.000)

\path(525,613)	(551.014,586.986)
	(600.000,538.000)

\path(6451,1287)	(6451.000,1343.482)
	(6451.000,1384.965)
	(6451.000,1437.000)

\path(6451,1437)	(6455.583,1473.495)
	(6451.000,1512.000)

\path(6451,1512)	(6379.382,1561.039)
	(6335.929,1575.825)
	(6301.000,1587.000)

\path(6301,1587)	(6236.471,1609.177)
	(6196.538,1621.625)
	(6154.188,1633.796)
	(6111.438,1644.806)
	(6070.310,1653.767)
	(6001.000,1662.000)

\path(6001,1662)	(5931.678,1656.915)
	(5890.541,1650.707)
	(5847.782,1642.189)
	(5805.426,1631.459)
	(5765.495,1618.619)
	(5701.000,1587.000)

\path(5701,1587)	(5656.757,1542.390)
	(5614.416,1480.226)
	(5577.867,1415.200)
	(5551.000,1362.000)

\path(5551,1362)	(5530.309,1313.664)
	(5507.755,1251.533)
	(5488.073,1188.384)
	(5476.000,1137.000)

\path(5476,1137)	(5469.161,1087.718)
	(5464.116,1025.959)
	(5465.013,963.470)
	(5476.000,912.000)

\path(5476,912)	(5502.519,870.070)
	(5542.799,827.325)
	(5587.180,789.417)
	(5626.000,762.000)

\path(5626,762)	(5673.260,737.524)
	(5734.079,713.160)
	(5797.108,694.466)
	(5851.000,687.000)

\path(5851,687)	(5904.802,695.309)
	(5967.681,715.406)
	(6028.470,740.051)
	(6076.000,762.000)

\path(6076,762)	(6128.222,791.542)
	(6191.620,832.549)
	(6253.457,875.781)
	(6301.000,912.000)

\path(6301,912)	(6340.128,947.734)
	(6376.000,987.000)

\path(6376,987)	(6405.366,1037.790)
	(6425.271,1079.585)
	(6451.000,1137.000)

\path(4951,462)	(4981.794,501.772)
	(5010.470,538.726)
	(5061.904,604.725)
	(5106.183,661.098)
	(5144.185,708.941)
	(5204.873,783.436)
	(5251.000,837.000)

\path(5251,837)	(5302.007,890.145)
	(5343.748,931.351)
	(5401.000,987.000)

\path(5551,1062)	(5599.634,1064.348)
	(5626.000,1062.000)

\path(5626,1062)	(5661.866,1050.610)
	(5704.458,1032.382)
	(5776.000,987.000)

\path(5776,987)	(5824.165,917.963)
	(5842.741,875.684)
	(5851.000,837.000)

\path(5851,837)	(5847.667,792.283)
	(5835.622,743.005)
	(5817.019,691.138)
	(5794.011,638.655)
	(5768.752,587.527)
	(5743.395,539.727)
	(5720.093,497.227)
	(5701.000,462.000)

\path(5701,462)	(5670.863,410.322)
	(5629.922,346.864)
	(5587.021,284.724)
	(5551.000,237.000)

\path(5551,237)	(5500.750,182.917)
	(5458.820,141.946)
	(5401.000,87.000)

\path(5551,912)	(5604.495,929.404)
	(5644.983,935.205)
	(5701.000,912.000)

\path(5701,912)	(5719.797,876.197)
	(5720.189,834.964)
	(5701.000,762.000)

\path(5701,762)	(5685.215,719.362)
	(5661.140,671.200)
	(5631.156,619.672)
	(5597.646,566.932)
	(5562.993,515.139)
	(5529.580,466.448)
	(5499.788,423.017)
	(5476.000,387.000)

\path(5476,387)	(5454.478,353.597)
	(5424.648,308.479)
	(5382.993,246.371)
	(5356.634,207.298)
	(5326.000,162.000)

\path(5368.423,278.191)(5326.000,162.000)(5418.110,244.558)
\path(5075.691,500.273)(5026.000,387.000)(5123.150,463.563)
\path(5026,387)	(5056.794,426.772)
	(5085.470,463.726)
	(5136.904,529.725)
	(5181.183,586.098)
	(5219.185,633.941)
	(5279.873,708.436)
	(5326.000,762.000)

\path(5326,762)	(5351.504,788.569)
	(5401.000,837.000)

\path(4950,538)	(4943.469,593.132)
	(4941.293,634.165)
	(4950.000,688.000)

\path(4950,688)	(4981.695,731.305)
	(5025.000,763.000)

\path(5025,763)	(5075.923,776.061)
	(5137.500,780.415)
	(5199.077,776.061)
	(5250.000,763.000)

\path(5250,763)	(5290.417,734.687)
	(5330.640,691.416)
	(5368.043,646.437)
	(5400.000,613.000)

\path(5400,613)	(5449.293,579.914)
	(5514.049,541.904)
	(5578.030,501.941)
	(5625.000,463.000)

\path(5625,463)	(5669.243,391.232)
	(5688.180,348.782)
	(5700.000,313.000)

\path(5700,313)	(5710.446,263.780)
	(5718.450,202.829)
	(5717.229,140.713)
	(5700.000,88.000)

\path(5700,88)	(5664.348,60.527)
	(5612.925,41.654)
	(5550.000,13.000)

\path(5550,13)	(5550.000,13.000)

\path(5550,13)	(5518.909,13.000)
	(5475.000,13.000)

\path(5475,13)	(5448.986,13.000)
	(5400.000,13.000)

\put(2850,763){\makebox(0,0)[lb]{$K_1$}}
\put(0,88){\makebox(0,0)[lb]{$K_2$}}
\put(7201,762){\makebox(0,0)[lb]{$K_1$}}
\put(4351,87){\makebox(0,0)[lb]{$K_2$}}
\end{picture} }
        	\hspace{-1.9mm}
	\end{array}
 $$
\caption{Shown here are the two links $E_1^{'}$ and $E_2^{'}$ 
mentioned in
equation \eqref{eq.4T2}. Note that $E_1^{'}$ contains a $1$-pair 
blink and $E_2^{'}$ contains a $1$-component boundary 
link.}\lbl{E1E2}
\end{figure}
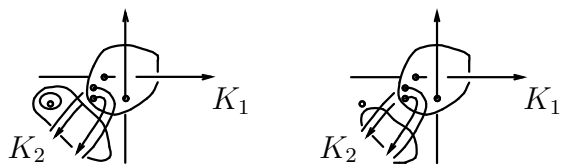

\begin{pf}
Moving the arc of $K$ (of figure \ref{3arcs2}) from the $SW$ quarter 
 to the $NE$ quarter as in proposition \ref{prop.4T1}, it follows
that:
\begin{equation}
\lbl{eq.well}
[C_s \cup L] - [C_n \cup L]- [C_e \cup L] +[C_w \cup L] = [E_1] - 
[E_2]
\end{equation}
where $E_1, E_2$ are as in figure \ref{E11E22}.
Using figure \ref{Kauf} the result follows.
\end{pf}

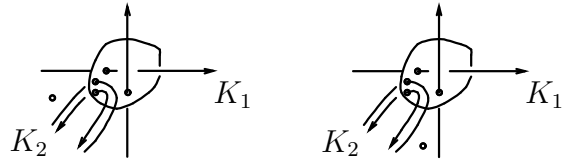
\begin{figure}[htpb]
$$\printname{E11E22}
	\setlength{\unitlength}{0.03\standardunitlength}
	\begin{array}{c}  \hspace{-1.7mm}
        	\raisebox{-8pt}{\begingroup\makeatletter\ifx\SetFigFont\undefined
% extract first six characters in \fmtname
\def\x#1#2#3#4#5#6#7\relax{\def\x{#1#2#3#4#5#6}}%
\expandafter\x\fmtname xxxxxx\relax \def\y{splain}%
\ifx\x\y   % LaTeX or SliTeX?
\gdef\SetFigFont#1#2#3{%
  \ifnum #1<17\tiny\else \ifnum #1<20\small\else
  \ifnum #1<24\normalsize\else \ifnum #1<29\large\else
  \ifnum #1<34\Large\else \ifnum #1<41\LARGE\else
     \huge\fi\fi\fi\fi\fi\fi
  \csname #3\endcsname}%
\else
\gdef\SetFigFont#1#2#3{\begingroup
  \count@#1\relax \ifnum 25<\count@\count@25\fi
  \def\x{\endgroup\@setsize\SetFigFont{#2pt}}%
  \expandafter\x
    \csname \romannumeral\the\count@ pt\expandafter\endcsname
    \csname @\romannumeral\the\count@ pt\endcsname
  \csname #3\endcsname}%
\fi
\fi\endgroup
\begin{picture}(7708,2140)(0,-10)
\thicklines
\put(1350,1213){\whiten\ellipse{76}{76}}
\put(1350,1213){\ellipse{76}{76}}
\put(1650,913){\whiten\ellipse{76}{76}}
\put(1650,913){\ellipse{76}{76}}
\put(1200,1063){\whiten\ellipse{76}{76}}
\put(1200,1063){\ellipse{76}{76}}
\put(1200,913){\whiten\ellipse{76}{76}}
\put(1200,913){\ellipse{76}{76}}
\put(600,838){\whiten\ellipse{76}{76}}
\put(600,838){\ellipse{76}{76}}
\put(5701,1212){\whiten\ellipse{76}{76}}
\put(5701,1212){\ellipse{76}{76}}
\put(6001,912){\whiten\ellipse{76}{76}}
\put(6001,912){\ellipse{76}{76}}
\put(5551,1062){\whiten\ellipse{76}{76}}
\put(5551,1062){\ellipse{76}{76}}
\put(5551,912){\whiten\ellipse{76}{76}}
\put(5551,912){\ellipse{76}{76}}
\put(5775,163){\whiten\ellipse{76}{76}}
\put(5775,163){\ellipse{76}{76}}
\path(450,1213)(1125,1213)
\path(1650,13)(1650,688)
\path(1350,1213)(1500,1213)
\path(1800,1213)(2850,1213)
\path(2730.000,1183.000)(2850.000,1213.000)(2730.000,1243.000)
\path(1650,913)(1650,2113)
\path(1680.000,1993.000)(1650.000,2113.000)(1620.000,1993.000)
\path(4801,1212)(5476,1212)
\path(6001,12)(6001,687)
\path(5701,1212)(5851,1212)
\path(6151,1212)(7201,1212)
\path(7081.000,1182.000)(7201.000,1212.000)(7081.000,1242.000)
\path(6001,912)(6001,2112)
\path(6031.000,1992.000)(6001.000,2112.000)(5971.000,1992.000)
\path(2100,1288)	(2100.000,1344.482)
	(2100.000,1385.965)
	(2100.000,1438.000)

\path(2100,1438)	(2104.582,1474.495)
	(2100.000,1513.000)

\path(2100,1513)	(2028.382,1562.039)
	(1984.929,1576.825)
	(1950.000,1588.000)

\path(1950,1588)	(1885.471,1610.177)
	(1845.538,1622.625)
	(1803.188,1634.796)
	(1760.438,1645.806)
	(1719.310,1654.767)
	(1650.000,1663.000)

\path(1650,1663)	(1580.678,1657.915)
	(1539.541,1651.707)
	(1496.783,1643.189)
	(1454.426,1632.459)
	(1414.495,1619.619)
	(1350.000,1588.000)

\path(1350,1588)	(1305.757,1543.390)
	(1263.416,1481.226)
	(1226.867,1416.200)
	(1200.000,1363.000)

\path(1200,1363)	(1179.309,1314.664)
	(1156.755,1252.533)
	(1137.073,1189.384)
	(1125.000,1138.000)

\path(1125,1138)	(1118.161,1088.718)
	(1113.116,1026.959)
	(1114.013,964.470)
	(1125.000,913.000)

\path(1125,913)	(1151.519,871.070)
	(1191.799,828.325)
	(1236.180,790.417)
	(1275.000,763.000)

\path(1275,763)	(1322.260,738.524)
	(1383.079,714.160)
	(1446.108,695.466)
	(1500.000,688.000)

\path(1500,688)	(1553.802,696.309)
	(1616.681,716.406)
	(1677.470,741.051)
	(1725.000,763.000)

\path(1725,763)	(1777.223,792.542)
	(1840.620,833.549)
	(1902.458,876.781)
	(1950.000,913.000)

\path(1950,913)	(1989.127,948.734)
	(2025.000,988.000)

\path(2025,988)	(2054.366,1038.790)
	(2074.271,1080.585)
	(2100.000,1138.000)

\path(600,463)	(630.794,502.772)
	(659.470,539.726)
	(710.904,605.725)
	(755.183,662.098)
	(793.185,709.941)
	(853.873,784.436)
	(900.000,838.000)

\path(900,838)	(951.007,891.145)
	(992.748,932.351)
	(1050.000,988.000)

\path(1200,1063)	(1248.634,1065.348)
	(1275.000,1063.000)

\path(1275,1063)	(1310.866,1051.610)
	(1353.457,1033.382)
	(1425.000,988.000)

\path(1425,988)	(1473.165,918.963)
	(1491.741,876.684)
	(1500.000,838.000)

\path(1500,838)	(1496.667,793.283)
	(1484.622,744.005)
	(1466.019,692.138)
	(1443.011,639.655)
	(1417.752,588.527)
	(1392.395,540.727)
	(1369.093,498.227)
	(1350.000,463.000)

\path(1350,463)	(1319.863,411.322)
	(1278.922,347.864)
	(1236.021,285.724)
	(1200.000,238.000)

\path(1200,238)	(1149.750,183.917)
	(1107.820,142.946)
	(1050.000,88.000)

\path(1200,913)	(1253.495,930.404)
	(1293.983,936.205)
	(1350.000,913.000)

\path(1350,913)	(1368.797,877.197)
	(1369.189,835.964)
	(1350.000,763.000)

\path(1350,763)	(1334.215,720.362)
	(1310.140,672.200)
	(1280.156,620.672)
	(1246.646,567.932)
	(1211.993,516.139)
	(1178.580,467.448)
	(1148.788,424.017)
	(1125.000,388.000)

\path(1125,388)	(1103.478,354.597)
	(1073.648,309.479)
	(1031.993,247.371)
	(1005.634,208.298)
	(975.000,163.000)

\path(1017.423,279.191)(975.000,163.000)(1067.110,245.558)
\path(724.691,501.273)(675.000,388.000)(772.150,464.563)
\path(675,388)	(705.794,427.772)
	(734.470,464.726)
	(785.904,530.725)
	(830.183,587.098)
	(868.185,634.941)
	(928.873,709.436)
	(975.000,763.000)

\path(975,763)	(1000.504,789.569)
	(1050.000,838.000)

\path(6451,1287)	(6451.000,1343.482)
	(6451.000,1384.965)
	(6451.000,1437.000)

\path(6451,1437)	(6455.583,1473.495)
	(6451.000,1512.000)

\path(6451,1512)	(6379.382,1561.039)
	(6335.929,1575.825)
	(6301.000,1587.000)

\path(6301,1587)	(6236.471,1609.177)
	(6196.538,1621.625)
	(6154.188,1633.796)
	(6111.438,1644.806)
	(6070.310,1653.767)
	(6001.000,1662.000)

\path(6001,1662)	(5931.678,1656.915)
	(5890.541,1650.707)
	(5847.782,1642.189)
	(5805.426,1631.459)
	(5765.495,1618.619)
	(5701.000,1587.000)

\path(5701,1587)	(5656.757,1542.390)
	(5614.416,1480.226)
	(5577.867,1415.200)
	(5551.000,1362.000)

\path(5551,1362)	(5530.309,1313.664)
	(5507.755,1251.533)
	(5488.073,1188.384)
	(5476.000,1137.000)

\path(5476,1137)	(5469.161,1087.718)
	(5464.116,1025.959)
	(5465.013,963.470)
	(5476.000,912.000)

\path(5476,912)	(5502.519,870.070)
	(5542.799,827.325)
	(5587.180,789.417)
	(5626.000,762.000)

\path(5626,762)	(5673.260,737.524)
	(5734.079,713.160)
	(5797.108,694.466)
	(5851.000,687.000)

\path(5851,687)	(5904.802,695.309)
	(5967.681,715.406)
	(6028.470,740.051)
	(6076.000,762.000)

\path(6076,762)	(6128.222,791.542)
	(6191.620,832.549)
	(6253.457,875.781)
	(6301.000,912.000)

\path(6301,912)	(6340.128,947.734)
	(6376.000,987.000)

\path(6376,987)	(6405.366,1037.790)
	(6425.271,1079.585)
	(6451.000,1137.000)

\path(4951,462)	(4981.794,501.772)
	(5010.470,538.726)
	(5061.904,604.725)
	(5106.183,661.098)
	(5144.185,708.941)
	(5204.873,783.436)
	(5251.000,837.000)

\path(5251,837)	(5302.007,890.145)
	(5343.748,931.351)
	(5401.000,987.000)

\path(5551,1062)	(5599.634,1064.348)
	(5626.000,1062.000)

\path(5626,1062)	(5661.866,1050.610)
	(5704.458,1032.382)
	(5776.000,987.000)

\path(5776,987)	(5824.165,917.963)
	(5842.741,875.684)
	(5851.000,837.000)

\path(5851,837)	(5847.667,792.283)
	(5835.622,743.005)
	(5817.019,691.138)
	(5794.011,638.655)
	(5768.752,587.527)
	(5743.395,539.727)
	(5720.093,497.227)
	(5701.000,462.000)

\path(5701,462)	(5670.863,410.322)
	(5629.922,346.864)
	(5587.021,284.724)
	(5551.000,237.000)

\path(5551,237)	(5500.750,182.917)
	(5458.820,141.946)
	(5401.000,87.000)

\path(5551,912)	(5604.495,929.404)
	(5644.983,935.205)
	(5701.000,912.000)

\path(5701,912)	(5719.797,876.197)
	(5720.189,834.964)
	(5701.000,762.000)

\path(5701,762)	(5685.215,719.362)
	(5661.140,671.200)
	(5631.156,619.672)
	(5597.646,566.932)
	(5562.993,515.139)
	(5529.580,466.448)
	(5499.788,423.017)
	(5476.000,387.000)

\path(5476,387)	(5454.478,353.597)
	(5424.648,308.479)
	(5382.993,246.371)
	(5356.634,207.298)
	(5326.000,162.000)

\path(5368.423,278.191)(5326.000,162.000)(5418.110,244.558)
\path(5075.691,500.273)(5026.000,387.000)(5123.150,463.563)
\path(5026,387)	(5056.794,426.772)
	(5085.470,463.726)
	(5136.904,529.725)
	(5181.183,586.098)
	(5219.185,633.941)
	(5279.873,708.436)
	(5326.000,762.000)

\path(5326,762)	(5351.504,788.569)
	(5401.000,837.000)

\put(2850,763){\makebox(0,0)[lb]{$K_1$}}
\put(0,88){\makebox(0,0)[lb]{$K_2$}}
\put(7201,762){\makebox(0,0)[lb]{$K_1$}}
\put(4351,87){\makebox(0,0)[lb]{$K_2$}}
\end{picture} }
        	\hspace{-1.9mm}
	\end{array}
 $$
\caption{Shown here are the links $E_1$ and $E_2$ of equation
\eqref{eq.well}.}\lbl{E11E22}
\end{figure}

A few remarks are in order:

\begin{remark}
\lbl{rem.4T2a}
In analogy with remark \ref{rem.4T1}, 
if we represent the knots $K_1, K_2$ by two circles and the knots 
$C, C_j$
(for $j=s,n,e,w$) by chords (that intersect the above mentioned 
circles
in two points each, namely the points of intersection $K_1 \cap D_j$ 
and
$K_2 \cap D_j$) then the $4$-term relation reads as in figure 
\ref{4Tv2}. Note also that the error terms contain $1$-pair blinks in
the complement of the link $L-C$. The $4$-term relation will be
used, in the form of figure \ref{4Tv2}, in the proof of theorem 
\ref{thm.as2b}. 
\end{remark}

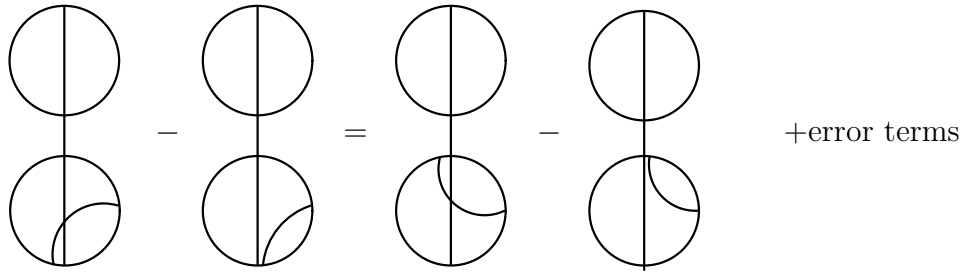
\begin{figure}[htpb]
$$\printname{4Tv2}
	\setlength{\unitlength}{0.03\standardunitlength}
	\begin{array}{c}  \hspace{-1.7mm}
        	\raisebox{-8pt}{\begingroup\makeatletter\ifx\SetFigFont\undefined
% extract first six characters in \fmtname
\def\x#1#2#3#4#5#6#7\relax{\def\x{#1#2#3#4#5#6}}%
\expandafter\x\fmtname xxxxxx\relax \def\y{splain}%
\ifx\x\y   % LaTeX or SliTeX?
\gdef\SetFigFont#1#2#3{%
  \ifnum #1<17\tiny\else \ifnum #1<20\small\else
  \ifnum #1<24\normalsize\else \ifnum #1<29\large\else
  \ifnum #1<34\Large\else \ifnum #1<41\LARGE\else
     \huge\fi\fi\fi\fi\fi\fi
  \csname #3\endcsname}%
\else
\gdef\SetFigFont#1#2#3{\begingroup
  \count@#1\relax \ifnum 25<\count@\count@25\fi
  \def\x{\endgroup\@setsize\SetFigFont{#2pt}}%
  \expandafter\x
    \csname \romannumeral\the\count@ pt\expandafter\endcsname
    \csname @\romannumeral\the\count@ pt\endcsname
  \csname #3\endcsname}%
\fi
\fi\endgroup
\begin{picture}(12561,3725)(0,-10)
\thicklines
\put(1320.500,229.500){\arc{1423.815}{2.9401}{5.0008}}
\put(4518.000,-18.000){\arc{1951.333}{3.2494}{4.4052}}
\put(6646.077,1419.692){\arc{1290.298}{1.1271}{3.4039}}
\put(9590.732,1486.709){\arc{1301.019}{1.5212}{3.2964}}
\put(784,836){\ellipse{1530}{1530}}
\put(3473,837){\ellipse{1530}{1530}}
\put(6173,837){\ellipse{1530}{1530}}
\put(8873,837){\ellipse{1530}{1530}}
\put(773,2937){\ellipse{1530}{1530}}
\put(3473,2937){\ellipse{1530}{1530}}
\put(6173,2937){\ellipse{1530}{1530}}
\put(8873,2862){\ellipse{1530}{1530}}
\path(773,3687)(773,87)
\path(3473,3687)(3473,87)
\path(6173,3687)(6173,87)
\path(8873,3612)(8873,12)
\put(2048,1812){\makebox(0,0)[lb]{$-$}}
\put(4673,1812){\makebox(0,0)[lb]{$=$}}
\put(7373,1812){\makebox(0,0)[lb]{$-$}}
\put(10823,1812){\makebox(0,0)[lb]{$+\text{error terms}$}}
\end{picture} }
        	\hspace{-1.9mm}
	\end{array}
 $$
\caption{The second  version of the $4$-term relation on $\cal M$, 
with
the notation of remark \ref{rem.4T2a}.}\lbl{4Tv2}
\end{figure}

\begin{remark}
\lbl{rem.4T2b}
Proposition \ref{prop.4T2} implies proposition \ref{prop.4T1}. 
Indeed,
consider the figure \ref{E1alt}.
Notice that in this case, the error terms vanish.
The reason that we introduced proposition \ref{prop.4T1} at all was 
as a warm-up exercise to make the proof of 
proposition
\ref{prop.4T2} more accessible, in light of the similarity between
proposition \ref{prop.4T1} and the $4$-term relation in the theory
of finite type knot invariants, see \cite{B-N1}.
\end{remark}

\begin{figure}[htpb]
$$\printname{E1alt}
	\setlength{\unitlength}{0.03\standardunitlength}
	\begin{array}{c}  \hspace{-1.7mm}
        	\raisebox{-8pt}{\begingroup\makeatletter\ifx\SetFigFont\undefined
% extract first six characters in \fmtname
\def\x#1#2#3#4#5#6#7\relax{\def\x{#1#2#3#4#5#6}}%
\expandafter\x\fmtname xxxxxx\relax \def\y{splain}%
\ifx\x\y   % LaTeX or SliTeX?
\gdef\SetFigFont#1#2#3{%
  \ifnum #1<17\tiny\else \ifnum #1<20\small\else
  \ifnum #1<24\normalsize\else \ifnum #1<29\large\else
  \ifnum #1<34\Large\else \ifnum #1<41\LARGE\else
     \huge\fi\fi\fi\fi\fi\fi
  \csname #3\endcsname}%
\else
\gdef\SetFigFont#1#2#3{\begingroup
  \count@#1\relax \ifnum 25<\count@\count@25\fi
  \def\x{\endgroup\@setsize\SetFigFont{#2pt}}%
  \expandafter\x
    \csname \romannumeral\the\count@ pt\expandafter\endcsname
    \csname @\romannumeral\the\count@ pt\endcsname
  \csname #3\endcsname}%
\fi
\fi\endgroup
\begin{picture}(3357,2303)(0,-10)
\thicklines
\put(1350,1376){\whiten\ellipse{76}{76}}
\put(1350,1376){\ellipse{76}{76}}
\put(1650,1076){\whiten\ellipse{76}{76}}
\put(1650,1076){\ellipse{76}{76}}
\put(675,1076){\whiten\ellipse{76}{76}}
\put(675,1076){\ellipse{76}{76}}
\put(1200,1151){\whiten\ellipse{76}{76}}
\put(1200,1151){\ellipse{76}{76}}
\put(1275,1076){\whiten\ellipse{76}{76}}
\put(1275,1076){\ellipse{76}{76}}
\path(450,1376)(1125,1376)
\path(1650,176)(1650,851)
\path(1350,1376)(1500,1376)
\path(1800,1376)(2850,1376)
\path(2730.000,1346.000)(2850.000,1376.000)(2730.000,1406.000)
\path(1650,1076)(1650,2276)
\path(1680.000,2156.000)(1650.000,2276.000)(1620.000,2156.000)
\path(2100,1451)	(2100.000,1507.482)
	(2100.000,1548.965)
	(2100.000,1601.000)

\path(2100,1601)	(2104.582,1637.495)
	(2100.000,1676.000)

\path(2100,1676)	(2028.382,1725.039)
	(1984.929,1739.825)
	(1950.000,1751.000)

\path(1950,1751)	(1885.471,1773.177)
	(1845.538,1785.625)
	(1803.188,1797.796)
	(1760.438,1808.806)
	(1719.310,1817.767)
	(1650.000,1826.000)

\path(1650,1826)	(1580.678,1820.915)
	(1539.541,1814.707)
	(1496.783,1806.189)
	(1454.426,1795.459)
	(1414.495,1782.619)
	(1350.000,1751.000)

\path(1350,1751)	(1305.757,1706.390)
	(1263.416,1644.226)
	(1226.867,1579.200)
	(1200.000,1526.000)

\path(1200,1526)	(1179.309,1477.664)
	(1156.755,1415.533)
	(1137.073,1352.384)
	(1125.000,1301.000)

\path(1125,1301)	(1118.161,1251.718)
	(1113.116,1189.959)
	(1114.013,1127.470)
	(1125.000,1076.000)

\path(1125,1076)	(1151.519,1034.070)
	(1191.799,991.325)
	(1236.180,953.417)
	(1275.000,926.000)

\path(1275,926)	(1322.260,901.524)
	(1383.079,877.160)
	(1446.108,858.466)
	(1500.000,851.000)

\path(1500,851)	(1553.802,859.309)
	(1616.681,879.406)
	(1677.470,904.051)
	(1725.000,926.000)

\path(1725,926)	(1777.223,955.542)
	(1840.620,996.549)
	(1902.458,1039.781)
	(1950.000,1076.000)

\path(1950,1076)	(1989.127,1111.734)
	(2025.000,1151.000)

\path(2025,1151)	(2054.366,1201.790)
	(2074.271,1243.585)
	(2100.000,1301.000)

\path(600,626)	(630.794,665.772)
	(659.470,702.726)
	(710.904,768.725)
	(755.183,825.098)
	(793.185,872.941)
	(853.873,947.436)
	(900.000,1001.000)

\path(900,1001)	(951.007,1054.145)
	(992.748,1095.351)
	(1050.000,1151.000)

\path(1200,1151)	(1255.964,1154.527)
	(1297.275,1155.702)
	(1350.000,1151.000)

\path(1350,1151)	(1386.866,1143.236)
	(1431.124,1131.050)
	(1500.000,1076.000)

\path(1500,1076)	(1505.082,1014.413)
	(1482.634,952.561)
	(1450.119,896.179)
	(1425.000,851.000)

\path(1425,851)	(1394.569,783.501)
	(1374.770,742.376)
	(1353.416,699.268)
	(1331.636,656.376)
	(1310.555,615.900)
	(1275.000,551.000)

\path(1275,551)	(1244.015,500.039)
	(1203.169,436.501)
	(1160.738,373.963)
	(1125.000,326.000)

\path(1125,326)	(1099.879,298.959)
	(1050.000,251.000)

\path(975,326)	(917.585,300.271)
	(875.790,280.366)
	(825.000,251.000)

\path(825,251)	(782.033,221.056)
	(750.000,176.000)

\path(750,176)	(751.741,133.245)
	(767.576,90.661)
	(825.000,26.000)

\path(825,26)	(859.561,12.252)
	(900.000,7.670)
	(940.439,12.252)
	(975.000,26.000)

\path(975,26)	(1025.580,96.406)
	(1040.559,139.708)
	(1050.000,176.000)

\path(1050,176)	(1052.348,202.366)
	(1050.000,251.000)

\path(675,551)	(705.794,590.772)
	(734.470,627.726)
	(785.904,693.725)
	(830.183,750.098)
	(868.185,797.941)
	(928.873,872.436)
	(975.000,926.000)

\path(975,926)	(1000.504,952.569)
	(1050.000,1001.000)

\path(1200,1076)	(1253.495,1093.404)
	(1293.983,1099.205)
	(1350.000,1076.000)

\path(1350,1076)	(1368.797,1040.197)
	(1369.189,998.964)
	(1350.000,926.000)

\path(1350,926)	(1334.215,883.362)
	(1310.140,835.200)
	(1280.156,783.672)
	(1246.646,730.932)
	(1211.993,679.139)
	(1178.580,630.448)
	(1148.788,587.017)
	(1125.000,551.000)

\path(1125,551)	(1103.478,517.597)
	(1073.648,472.479)
	(1031.993,410.371)
	(1005.634,371.298)
	(975.000,326.000)

%\put(2025,776){\makebox(0,0)[lb]{\smash{{{\SetFigFont{12}{14.4}{rm}\$c\$}}}}}
%\put(900,1676){\makebox(0,0)[lb]{\smash{{{\SetFigFont{12}{14.4}{rm}\$b\$}}}}}
%\put(225,1001){\makebox(0,0)[lb]{\smash{{{\SetFigFont{12}{14.4}{rm}\$a\$}}}}}
\put(2850,926){\makebox(0,0)[lb]{$K_1$}}
\put(0,251){\makebox(0,0)[lb]{$K_2$}}
\end{picture} }
        	\hspace{-1.9mm}
	\end{array}
 $$
\caption{A special link that reduces proposition \ref{prop.4T2}
to proposition \ref{prop.4T1}.}\lbl{E1alt}
\end{figure}
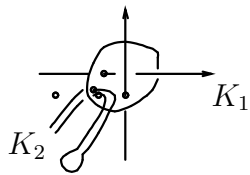

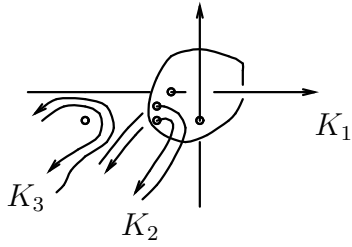
\begin{figure}[htpb]
$$\printname{3arcs3}
	\setlength{\unitlength}{0.04\standardunitlength}
	\begin{array}{c}  \hspace{-1.7mm}
        	\raisebox{-8pt}{\begingroup\makeatletter\ifx\SetFigFont\undefined
% extract first six characters in \fmtname
\def\x#1#2#3#4#5#6#7\relax{\def\x{#1#2#3#4#5#6}}%
\expandafter\x\fmtname xxxxxx\relax \def\y{splain}%
\ifx\x\y   % LaTeX or SliTeX?
\gdef\SetFigFont#1#2#3{%
  \ifnum #1<17\tiny\else \ifnum #1<20\small\else
  \ifnum #1<24\normalsize\else \ifnum #1<29\large\else
  \ifnum #1<34\Large\else \ifnum #1<41\LARGE\else
     \huge\fi\fi\fi\fi\fi\fi
  \csname #3\endcsname}%
\else
\gdef\SetFigFont#1#2#3{\begingroup
  \count@#1\relax \ifnum 25<\count@\count@25\fi
  \def\x{\endgroup\@setsize\SetFigFont{#2pt}}%
  \expandafter\x
    \csname \romannumeral\the\count@ pt\expandafter\endcsname
    \csname @\romannumeral\the\count@ pt\endcsname
  \csname #3\endcsname}%
\fi
\fi\endgroup
\begin{picture}(3732,2466)(0,-10)
\thicklines
\put(1725,1539){\whiten\ellipse{76}{76}}
\put(1725,1539){\ellipse{76}{76}}
\put(2025,1239){\whiten\ellipse{76}{76}}
\put(2025,1239){\ellipse{76}{76}}
\put(1575,1389){\whiten\ellipse{76}{76}}
\put(1575,1389){\ellipse{76}{76}}
\put(1575,1239){\whiten\ellipse{76}{76}}
\put(1575,1239){\ellipse{76}{76}}
\put(825,1239){\whiten\ellipse{76}{76}}
\put(825,1239){\ellipse{76}{76}}
\path(825,1539)(1500,1539)
\path(2025,339)(2025,1014)
\path(1725,1539)(1875,1539)
\path(2175,1539)(3225,1539)
\path(3105.000,1509.000)(3225.000,1539.000)(3105.000,1569.000)
\path(2025,1239)(2025,2439)
\path(2055.000,2319.000)(2025.000,2439.000)(1995.000,2319.000)
\path(225,1539)(825,1539)
\path(2475,1614)	(2475.000,1670.482)
	(2475.000,1711.965)
	(2475.000,1764.000)

\path(2475,1764)	(2479.582,1800.495)
	(2475.000,1839.000)

\path(2475,1839)	(2403.383,1888.039)
	(2359.929,1902.825)
	(2325.000,1914.000)

\path(2325,1914)	(2260.471,1936.177)
	(2220.538,1948.625)
	(2178.188,1960.796)
	(2135.438,1971.806)
	(2094.310,1980.767)
	(2025.000,1989.000)

\path(2025,1989)	(1955.678,1983.915)
	(1914.541,1977.707)
	(1871.783,1969.189)
	(1829.426,1958.459)
	(1789.495,1945.619)
	(1725.000,1914.000)

\path(1725,1914)	(1680.757,1869.390)
	(1638.416,1807.226)
	(1601.867,1742.200)
	(1575.000,1689.000)

\path(1575,1689)	(1554.309,1640.664)
	(1531.755,1578.533)
	(1512.073,1515.384)
	(1500.000,1464.000)

\path(1500,1464)	(1493.161,1414.718)
	(1488.116,1352.959)
	(1489.013,1290.470)
	(1500.000,1239.000)

\path(1500,1239)	(1526.519,1197.070)
	(1566.799,1154.325)
	(1611.180,1116.417)
	(1650.000,1089.000)

\path(1650,1089)	(1697.260,1064.524)
	(1758.079,1040.160)
	(1821.108,1021.466)
	(1875.000,1014.000)

\path(1875,1014)	(1928.802,1022.309)
	(1991.681,1042.406)
	(2052.470,1067.051)
	(2100.000,1089.000)

\path(2100,1089)	(2152.222,1118.542)
	(2215.620,1159.549)
	(2277.457,1202.781)
	(2325.000,1239.000)

\path(2325,1239)	(2364.128,1274.734)
	(2400.000,1314.000)

\path(2400,1314)	(2429.366,1364.790)
	(2449.271,1406.585)
	(2475.000,1464.000)

\path(975,789)	(1005.794,828.772)
	(1034.470,865.726)
	(1085.904,931.725)
	(1130.183,988.098)
	(1168.185,1035.941)
	(1228.873,1110.436)
	(1275.000,1164.000)

\path(1275,1164)	(1326.007,1217.145)
	(1367.748,1258.351)
	(1425.000,1314.000)

\path(1575,1389)	(1623.634,1391.348)
	(1650.000,1389.000)

\path(1650,1389)	(1685.866,1377.610)
	(1728.457,1359.382)
	(1800.000,1314.000)

\path(1800,1314)	(1848.165,1244.963)
	(1866.741,1202.684)
	(1875.000,1164.000)

\path(1875,1164)	(1871.667,1119.283)
	(1859.622,1070.005)
	(1841.019,1018.138)
	(1818.011,965.655)
	(1792.752,914.527)
	(1767.395,866.727)
	(1744.093,824.227)
	(1725.000,789.000)

\path(1725,789)	(1694.863,737.322)
	(1653.922,673.864)
	(1611.021,611.724)
	(1575.000,564.000)

\path(1575,564)	(1524.750,509.917)
	(1482.820,468.946)
	(1425.000,414.000)

\path(1575,1239)	(1628.495,1256.404)
	(1668.983,1262.205)
	(1725.000,1239.000)

\path(1725,1239)	(1743.797,1203.197)
	(1744.189,1161.964)
	(1725.000,1089.000)

\path(1725,1089)	(1709.215,1046.362)
	(1685.140,998.200)
	(1655.156,946.672)
	(1621.646,893.932)
	(1586.993,842.139)
	(1553.580,793.448)
	(1523.788,750.017)
	(1500.000,714.000)

\path(1500,714)	(1478.478,680.597)
	(1448.648,635.479)
	(1406.993,573.371)
	(1380.634,534.298)
	(1350.000,489.000)

\path(1392.423,605.191)(1350.000,489.000)(1442.110,571.558)
\path(1099.691,827.273)(1050.000,714.000)(1147.150,790.563)
\path(1050,714)	(1080.794,753.772)
	(1109.470,790.726)
	(1160.904,856.725)
	(1205.183,913.098)
	(1243.185,960.941)
	(1303.873,1035.436)
	(1350.000,1089.000)

\path(1350,1089)	(1375.504,1115.569)
	(1425.000,1164.000)

\path(375,1239)	(417.105,1272.895)
	(453.898,1301.585)
	(514.181,1345.102)
	(561.124,1373.069)
	(600.000,1389.000)

\path(600,1389)	(650.528,1397.678)
	(712.886,1399.849)
	(775.051,1396.595)
	(825.000,1389.000)

\path(825,1389)	(861.292,1379.559)
	(904.594,1364.580)
	(975.000,1314.000)

\path(975,1314)	(994.838,1261.988)
	(1000.406,1201.043)
	(993.271,1140.325)
	(975.000,1089.000)

\path(975,1089)	(931.287,1040.908)
	(867.709,1000.628)
	(801.526,967.033)
	(750.000,939.000)

\path(750,939)	(699.997,906.705)
	(636.352,866.006)
	(573.282,824.305)
	(525.000,789.000)

\path(525,789)	(497.959,763.879)
	(450.000,714.000)

\path(510.481,821.898)(450.000,714.000)(554.142,780.742)
\path(374.086,1413.051)(300.000,1314.000)(411.983,1366.534)
\path(300,1314)	(342.105,1347.895)
	(378.898,1376.585)
	(439.181,1420.102)
	(486.124,1448.069)
	(525.000,1464.000)

\path(525,1464)	(592.548,1474.810)
	(633.305,1476.583)
	(675.986,1476.431)
	(718.544,1474.725)
	(758.931,1471.837)
	(825.000,1464.000)

\path(825,1464)	(876.614,1455.260)
	(940.144,1441.132)
	(1002.352,1419.689)
	(1050.000,1389.000)

\path(1050,1389)	(1082.580,1342.591)
	(1107.176,1282.447)
	(1121.934,1219.330)
	(1125.000,1164.000)

\path(1125,1164)	(1115.770,1125.885)
	(1096.545,1083.694)
	(1050.000,1014.000)

\path(1050,1014)	(1003.030,975.060)
	(939.049,935.100)
	(874.293,897.090)
	(825.000,864.000)

\path(825,864)	(789.709,833.170)
	(746.749,792.536)
	(705.414,750.135)
	(675.000,714.000)

\path(675,714)	(656.804,681.574)
	(637.500,639.000)
	(618.196,596.426)
	(600.000,564.000)

\path(600,564)	(575.396,536.224)
	(525.000,489.000)

\put(3225,1089){\makebox(0,0)[lb]{$K_1$}}
\put(1200,39){\makebox(0,0)[lb]{$K_2$}}
\put(0,339){\makebox(0,0)[lb]{$K_3$}}
\end{picture} }
        	\hspace{-1.9mm}
	\end{array}
 $$
\caption{ Shown here on the left hand side 
 are the $3$ arcs along which a ball $B$ 
intersects the 
knot $K_1$,
the $2$ arcs along which it intersects $K_2$ and the $2$ arcs that it 
intersects 
the knot $K_3$. Shown also is an unknotted component $C$ which bounds a 
disc $D$, and with bold dots are the 
points where $D$ intersects $K_1$ and $K_2$. On the right hand side of the
figure is shown one of the links ($C_w$)
 that appears in equation \eqref{eq.4T3}.}\lbl{3arcs3}
\end{figure}

We close this section with the third (and most refined) version of
the $4$-term relation on $\M$. 
 
\begin{proposition}
\lbl{prop.4T3}
Let $C \cup K_1 \cup K_2 \cup K_3$ be a four component sublink of 
an
\ASA\ link $L$ in $S^3$. Assume that there is a ball $B$ such that 
the
intersection of $L$ with $B$ is as in figure \ref{3arcs3}.
Let $D$ and $D_j$ (for $j=s,n,e,w$) be discs as in section 
\ref{sub.4T}
(with figure \ref{3arcs3} replacing figure \ref{3arcs}). Let $C_j$
be the boundary of the disc $D_j$ and $C$ the boundary of $D$. Then, 
with the abbreviations before
the statement of proposition \ref{prop.4T1}, we have the following:
\begin{equation}
\lbl{eq.4T3}
[C_s \cup L] - [C_n \cup L]=[C_e \cup L] -[C_w \cup L] + \text{error 
terms}
\end{equation}
where the error terms include $1$-pair blinks in the complement of
$L -( C \cup_j C_j)$.
In the rest of the paper, the above 
relation will be called  the (third version of the)
$4$-term relation on $\M$.
\end{proposition}

\begin{pf}
The proof is similar to that of proposition \ref{prop.4T2} and we 
briefly
sketch it here.
Using the equation of figure \ref{dotcross}, after  
moving the arc of $K$ (of figure \ref{3arcs3}) from the $SW$ quarter 
 to the $NE$ quarter, we get the following equality involving the
terms $[C_j \cup L]$, as well as two kinds of error terms: $ET_1$
(respectively, $ET_2$) that come from 
from moving the
$K_1$ arc around  the $K_3$ arcs (respectively,
around  the $K_2$ arcs):
\begin{equation}
[C_s \cup L] - [C_n \cup L]= [C_e \cup L] -[C_w \cup L] +
ET_1 + ET_2
\end{equation}
By  definition, we have:
\begin{equation}
ET_1= 
\sum_{i=1,2}([E_{i,s}]-[E_{i,n}]) - ([E_{i,e}]-[E_{i,w}])
\end{equation}
where $E_{i,j}$ (for $i=1,2$ and $j=s,n,e,w$) are $1$-pair blinks 
that
come from figure \ref{dotcross}.
Each of the differences in the parentheses is a sum (with signs)
of four terms (recall that for a $1$-pair blink $E$, $[E]$ is a sum 
(with
signs) of two terms), two of which cancel. The   remaining two 
(with signs) can be combined as $[e]$ for some $1$-pair blink
$e$. These $1$-pair blinks $\{ e \}$
%
%surgeries on $1$-pair blinks. But these two terms can be rewritten 
%as a single term containing an extra blink (see Figure ? ({\bf I think 
%another picture is needed here.})), whose 
bound Seifert surfaces in the
interior of the ball $B$ (of figure \ref{3arcs3}), and these surfaces
intersect the disc that $C$ bounds, where $C$ is as in figure 
\ref{3arcs3}.
By the same 
argument as in proposition \ref{prop.4T2}, it follows that
$ET_2$ can be written as a linear combination of links that contain 
$1$-pair
blinks that bound surfaces in the interior of the ball $B$.
This proves proposition \ref{prop.4T3}.
\end{pf}

\begin{figure}[htpb]
$$\printname{dotcross}
	\setlength{\unitlength}{0.03\standardunitlength}
	\begin{array}{c}  \hspace{-1.7mm}
        	\raisebox{-8pt}{\begingroup\makeatletter\ifx\SetFigFont\undefined
% extract first six characters in \fmtname
\def\x#1#2#3#4#5#6#7\relax{\def\x{#1#2#3#4#5#6}}%
\expandafter\x\fmtname xxxxxx\relax \def\y{splain}%
\ifx\x\y   % LaTeX or SliTeX?
\gdef\SetFigFont#1#2#3{%
  \ifnum #1<17\tiny\else \ifnum #1<20\small\else
  \ifnum #1<24\normalsize\else \ifnum #1<29\large\else
  \ifnum #1<34\Large\else \ifnum #1<41\LARGE\else
     \huge\fi\fi\fi\fi\fi\fi
  \csname #3\endcsname}%
\else
\gdef\SetFigFont#1#2#3{\begingroup
  \count@#1\relax \ifnum 25<\count@\count@25\fi
  \def\x{\endgroup\@setsize\SetFigFont{#2pt}}%
  \expandafter\x
    \csname \romannumeral\the\count@ pt\expandafter\endcsname
    \csname @\romannumeral\the\count@ pt\endcsname
  \csname #3\endcsname}%
\fi
\fi\endgroup
\begin{picture}(8424,1426)(0,-10)
\thicklines
\put(612,91){\whiten\ellipse{76}{76}}
\put(612,91){\ellipse{76}{76}}
\put(5412,166){\whiten\ellipse{76}{76}}
\put(5412,166){\ellipse{76}{76}}
\put(3012,916){\whiten\ellipse{76}{76}}
\put(3012,916){\ellipse{76}{76}}
\put(7812,1366){\whiten\ellipse{76}{76}}
\put(7812,1366){\ellipse{76}{76}}
\path(12,616)(1212,616)
\path(2412,616)(3612,616)
\path(4812,616)(6012,616)
\path(7212,616)(8412,616)
\path(4962,466)(5037,466)
\path(5262,316)(5712,316)
\path(7362,466)(7512,466)
\path(7737,316)(8112,316)
\path(5262,466)(5712,466)
\path(5592.000,436.000)(5712.000,466.000)(5592.000,496.000)
\path(7662,466)(8112,466)
\path(7992.000,436.000)(8112.000,466.000)(7992.000,496.000)
\path(7482.000,346.000)(7362.000,316.000)(7482.000,286.000)
\path(7362,316)(7512,316)
\path(5082.000,346.000)(4962.000,316.000)(5082.000,286.000)
\path(4962,316)(5037,316)
\path(5562,691)	(5566.613,746.672)
	(5568.150,787.885)
	(5562.000,841.000)

\path(5562,841)	(5542.507,899.403)
	(5512.294,966.689)
	(5469.434,1027.380)
	(5412.000,1066.000)

\path(5412,1066)	(5351.337,1072.089)
	(5288.812,1056.314)
	(5231.632,1026.632)
	(5187.000,991.000)

\path(5187,991)	(5156.311,943.352)
	(5134.868,881.144)
	(5120.740,817.614)
	(5112.000,766.000)

\path(5112,766)	(5106.042,699.491)
	(5104.797,658.719)
	(5104.579,615.812)
	(5105.288,572.929)
	(5106.826,532.228)
	(5112.000,466.000)

\path(5112,466)	(5122.106,398.471)
	(5129.626,357.377)
	(5138.629,314.309)
	(5148.990,271.452)
	(5160.587,230.992)
	(5187.000,166.000)

\path(5187,166)	(5212.365,125.391)
	(5247.716,79.514)
	(5290.209,39.380)
	(5337.000,16.000)

\path(5337,16)	(5395.870,11.620)
	(5459.081,23.237)
	(5517.502,49.986)
	(5562.000,91.000)

\path(5562,91)	(5580.406,144.484)
	(5578.361,206.976)
	(5568.135,267.730)
	(5562.000,316.000)

\path(5562,316)	(5562.000,368.035)
	(5562.000,409.518)
	(5562.000,466.000)

\path(8037,766)	(8041.613,821.672)
	(8043.150,862.885)
	(8037.000,916.000)

\path(8037,916)	(8017.507,974.403)
	(7987.294,1041.689)
	(7944.434,1102.380)
	(7887.000,1141.000)

\path(7887,1141)	(7826.337,1147.089)
	(7763.812,1131.314)
	(7706.632,1101.632)
	(7662.000,1066.000)

\path(7662,1066)	(7631.311,1018.352)
	(7609.868,956.144)
	(7595.740,892.614)
	(7587.000,841.000)

\path(7587,841)	(7581.042,774.491)
	(7579.797,733.719)
	(7579.579,690.812)
	(7580.288,647.929)
	(7581.826,607.228)
	(7587.000,541.000)

\path(7587,541)	(7597.106,473.471)
	(7604.626,432.377)
	(7613.629,389.309)
	(7623.990,346.452)
	(7635.587,305.992)
	(7662.000,241.000)

\path(7662,241)	(7687.365,200.391)
	(7722.716,154.514)
	(7765.209,114.380)
	(7812.000,91.000)

\path(7812,91)	(7870.870,86.620)
	(7934.081,98.237)
	(7992.502,124.986)
	(8037.000,166.000)

\path(8037,166)	(8055.406,219.484)
	(8053.361,281.976)
	(8043.135,342.730)
	(8037.000,391.000)

\path(8037,391)	(8037.000,443.035)
	(8037.000,484.518)
	(8037.000,541.000)

\path(7887,691)	(7851.611,749.594)
	(7812.000,766.000)

\path(7812,766)	(7761.142,741.520)
	(7737.000,691.000)

\path(7737,691)	(7738.741,648.245)
	(7754.576,605.661)
	(7812.000,541.000)

\path(7812,541)	(7850.505,536.418)
	(7887.000,541.000)

\path(7887,541)	(7887.000,541.000)

\path(312,466)(912,466)
\path(432.000,346.000)(312.000,316.000)(432.000,286.000)
\path(312,316)(912,316)
\path(2832.000,346.000)(2712.000,316.000)(2832.000,286.000)
\path(2712,316)(3312,316)
\path(2712,466)(3312,466)
\put(1587,541){\makebox(0,0)[lb]{$-$}}
\put(4062,541){\makebox(0,0)[lb]{$=$}}
\put(6462,541){\makebox(0,0)[lb]{$-$}}
\end{picture} }
        	\hspace{-1.9mm}
	\end{array}
 $$
\caption{Moving an arc perpendicular to the page past three others. 
On the
left hand side, there are four arcs shown, one perpendicular to the 
plane,
pointing towards your eyes. Of these four arcs, the top horizontal 
one and
the one represented by a dot, belong to the same link component, and 
so do
the two others.
Note that only two components of the link are involved, and that
this figure is a redrawing of figure \ref{Kauf}. The two blinks on the
right hand side are denoted by $E_1$ and $E_2$ respectively. 
}\lbl{dotcross}
\end{figure}
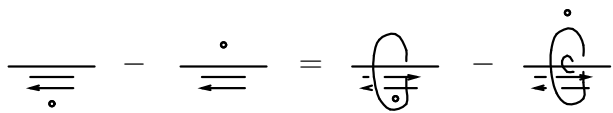

\begin{remark}
\lbl{rem.4T3}
In analogy with remark \ref{rem.4T2a}, 
if we represent the knots $K_1, K_2, K_3$ by three circles and the 
knots 
$C, C_j$
(for $j=s,n,e,w$) by chords (that intersect the above mentioned 
circles
in two points each, namely the points of intersection $K_i \cap D_j$ 
(for $i=1,2,3$) then the $4$-term relation reads as in figure 
\ref{4Tv3}. Note  that the error terms contain $1$-pair blinks in
the complement of the link $L-(C \cup_{j} C_j)$. The $4$-term 
relation will be
used, in the form of figure \ref{4Tv3}, in the proof of theorem 
\ref{thm.as2b}. 
\end{remark}

\begin{figure}[htpb]
$$\printname{4Tv3}
	\setlength{\unitlength}{0.03\standardunitlength}
	\begin{array}{c}  \hspace{-1.7mm}
        	\raisebox{-8pt}{\begingroup\makeatletter\ifx\SetFigFont\undefined
% extract first six characters in \fmtname
\def\x#1#2#3#4#5#6#7\relax{\def\x{#1#2#3#4#5#6}}%
\expandafter\x\fmtname xxxxxx\relax \def\y{splain}%
\ifx\x\y   % LaTeX or SliTeX?
\gdef\SetFigFont#1#2#3{%
  \ifnum #1<17\tiny\else \ifnum #1<20\small\else
  \ifnum #1<24\normalsize\else \ifnum #1<29\large\else
  \ifnum #1<34\Large\else \ifnum #1<41\LARGE\else
     \huge\fi\fi\fi\fi\fi\fi
  \csname #3\endcsname}%
\else
\gdef\SetFigFont#1#2#3{\begingroup
  \count@#1\relax \ifnum 25<\count@\count@25\fi
  \def\x{\endgroup\@setsize\SetFigFont{#2pt}}%
  \expandafter\x
    \csname \romannumeral\the\count@ pt\expandafter\endcsname
    \csname @\romannumeral\the\count@ pt\endcsname
  \csname #3\endcsname}%
\fi
\fi\endgroup
\begin{picture}(12846,3055)(0,-10)
\thicklines
\put(1605.500,166.500){\arc{1423.815}{2.9401}{5.0008}}
\put(4803.000,-81.000){\arc{1951.333}{3.2494}{4.4052}}
\put(6931.077,1356.692){\arc{1290.298}{1.1271}{3.4039}}
\put(9875.732,1423.709){\arc{1301.019}{1.5212}{3.2964}}
\put(1069,773){\ellipse{1530}{1530}}
\put(3758,774){\ellipse{1530}{1530}}
\put(6458,774){\ellipse{1530}{1530}}
\put(9158,774){\ellipse{1530}{1530}}
\put(458,2574){\ellipse{900}{900}}
\put(1658,2574){\ellipse{900}{900}}
\put(3158,2574){\ellipse{900}{900}}
\put(4358,2574){\ellipse{900}{900}}
\put(5858,2574){\ellipse{900}{900}}
\put(7058,2574){\ellipse{900}{900}}
\put(8558,2574){\ellipse{900}{900}}
\put(9758,2574){\ellipse{900}{900}}
\path(308,3024)(1058,24)
\path(3026,3028)(3776,28)
\path(5690,3019)(6440,19)
\path(8408,3024)(9158,24)
\path(7058,2949)(7133,3024)
\path(9758,2949)(9758,3024)
\path(1658,3024)	(1649.032,2985.005)
	(1640.764,2948.725)
	(1626.213,2883.761)
	(1614.128,2828.010)
	(1604.289,2780.374)
	(1596.476,2739.753)
	(1590.470,2705.050)
	(1583.000,2649.000)

\path(1583,2649)	(1577.776,2583.053)
	(1575.599,2542.510)
	(1574.176,2499.750)
	(1573.857,2456.896)
	(1574.989,2416.072)
	(1583.000,2349.000)

\path(1583,2349)	(1608.904,2278.063)
	(1628.515,2236.794)
	(1650.560,2194.211)
	(1673.535,2152.227)
	(1695.936,2112.754)
	(1733.000,2049.000)

\path(1733,2049)	(1763.988,1998.039)
	(1804.831,1934.501)
	(1847.259,1871.963)
	(1883.000,1824.000)

\path(1883,1824)	(1933.155,1776.064)
	(2000.622,1719.960)
	(2065.528,1659.626)
	(2108.000,1599.000)

\path(2108,1599)	(2117.796,1547.402)
	(2116.707,1484.700)
	(2111.265,1422.898)
	(2108.000,1374.000)

\path(2108,1374)	(2112.354,1308.806)
	(2116.163,1268.755)
	(2119.610,1226.404)
	(2121.605,1183.752)
	(2121.061,1142.800)
	(2108.000,1074.000)

\path(2108,1074)	(2072.803,1032.997)
	(2033.000,999.000)

\path(2033,999)	(1996.670,960.555)
	(1958.000,924.000)

\path(1958,924)	(1907.210,894.634)
	(1865.415,874.729)
	(1808.000,849.000)

\path(4358,3024)	(4349.032,2985.005)
	(4340.764,2948.725)
	(4326.213,2883.761)
	(4314.128,2828.010)
	(4304.289,2780.374)
	(4296.476,2739.753)
	(4290.470,2705.050)
	(4283.000,2649.000)

\path(4283,2649)	(4277.776,2583.053)
	(4275.599,2542.510)
	(4274.176,2499.750)
	(4273.857,2456.896)
	(4274.989,2416.072)
	(4283.000,2349.000)

\path(4283,2349)	(4308.904,2278.063)
	(4328.515,2236.794)
	(4350.560,2194.211)
	(4373.535,2152.227)
	(4395.936,2112.754)
	(4433.000,2049.000)

\path(4433,2049)	(4463.988,1998.039)
	(4504.831,1934.501)
	(4547.259,1871.963)
	(4583.000,1824.000)

\path(4583,1824)	(4633.155,1776.064)
	(4700.622,1719.960)
	(4765.528,1659.626)
	(4808.000,1599.000)

\path(4808,1599)	(4817.796,1547.402)
	(4816.708,1484.700)
	(4811.265,1422.898)
	(4808.000,1374.000)

\path(4808,1374)	(4812.354,1308.806)
	(4816.163,1268.755)
	(4819.610,1226.404)
	(4821.605,1183.752)
	(4821.061,1142.800)
	(4808.000,1074.000)

\path(4808,1074)	(4772.802,1032.997)
	(4733.000,999.000)

\path(4733,999)	(4696.670,960.555)
	(4658.000,924.000)

\path(4658,924)	(4607.210,894.634)
	(4565.415,874.729)
	(4508.000,849.000)

\path(9758,2949)	(9749.032,2910.005)
	(9740.764,2873.725)
	(9726.213,2808.761)
	(9714.128,2753.010)
	(9704.289,2705.374)
	(9696.476,2664.753)
	(9690.470,2630.050)
	(9683.000,2574.000)

\path(9683,2574)	(9677.776,2508.053)
	(9675.599,2467.510)
	(9674.176,2424.750)
	(9673.857,2381.896)
	(9674.989,2341.072)
	(9683.000,2274.000)

\path(9683,2274)	(9708.904,2203.063)
	(9728.515,2161.794)
	(9750.560,2119.211)
	(9773.535,2077.227)
	(9795.936,2037.754)
	(9833.000,1974.000)

\path(9833,1974)	(9863.988,1923.039)
	(9904.831,1859.501)
	(9947.259,1796.963)
	(9983.000,1749.000)

\path(9983,1749)	(10033.155,1701.064)
	(10100.622,1644.960)
	(10165.528,1584.626)
	(10208.000,1524.000)

\path(10208,1524)	(10217.796,1472.402)
	(10216.708,1409.700)
	(10211.265,1347.898)
	(10208.000,1299.000)

\path(10208,1299)	(10212.354,1233.806)
	(10216.163,1193.755)
	(10219.610,1151.404)
	(10221.605,1108.752)
	(10221.061,1067.800)
	(10208.000,999.000)

\path(10208,999)	(10172.802,957.997)
	(10133.000,924.000)

\path(10133,924)	(10096.670,885.555)
	(10058.000,849.000)

\path(10058,849)	(10007.210,819.634)
	(9965.415,799.729)
	(9908.000,774.000)

\path(7058,2949)	(7049.032,2910.005)
	(7040.764,2873.725)
	(7026.213,2808.761)
	(7014.128,2753.010)
	(7004.289,2705.374)
	(6996.476,2664.753)
	(6990.470,2630.050)
	(6983.000,2574.000)

\path(6983,2574)	(6977.776,2508.053)
	(6975.599,2467.510)
	(6974.176,2424.750)
	(6973.857,2381.896)
	(6974.989,2341.072)
	(6983.000,2274.000)

\path(6983,2274)	(7008.904,2203.063)
	(7028.515,2161.794)
	(7050.560,2119.211)
	(7073.535,2077.227)
	(7095.936,2037.754)
	(7133.000,1974.000)

\path(7133,1974)	(7163.988,1923.039)
	(7204.831,1859.501)
	(7247.259,1796.963)
	(7283.000,1749.000)

\path(7283,1749)	(7333.155,1701.064)
	(7400.622,1644.960)
	(7465.528,1584.626)
	(7508.000,1524.000)

\path(7508,1524)	(7517.796,1472.402)
	(7516.707,1409.700)
	(7511.265,1347.898)
	(7508.000,1299.000)

\path(7508,1299)	(7512.354,1233.806)
	(7516.163,1193.755)
	(7519.610,1151.404)
	(7521.605,1108.752)
	(7521.061,1067.800)
	(7508.000,999.000)

\path(7508,999)	(7472.802,957.997)
	(7433.000,924.000)

\path(7433,924)	(7396.670,885.555)
	(7358.000,849.000)

\path(7358,849)	(7307.210,819.634)
	(7265.415,799.729)
	(7208.000,774.000)

\put(2333,1749){\makebox(0,0)[lb]{$-$}}
\put(4958,1749){\makebox(0,0)[lb]{$=$}}
\put(7658,1749){\makebox(0,0)[lb]{$-$}}
\put(11108,1749){\makebox(0,0)[lb]{$+$\text{error terms}}}
\end{picture} }
        	\hspace{-1.9mm}
	\end{array}
 $$
\caption{The third  version of the $4$-term relation on $\cal M$, 
with
the notation of remark   \ref{rem.4T3}.}\lbl{4Tv3}
\end{figure}
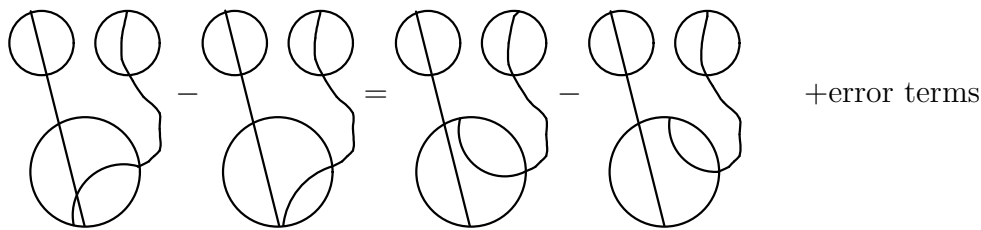

\subsection{Proof of theorem \ref{thm.as2b}}
\lbl{sub.1}

This section is devoted to the proof of theorem \ref{thm.as2b}.
The proof is long, and rather involved.
For the convenience of the reader, we divide the proof in three
propositions, each of which needs independent arguments.

All links considered in this section are
\ASA\ links in $S^3$.
Let us begin by introducing a  definition that is  useful
in this section. 

\begin{definition}
\lbl{def.near}
An \ASA\ link $(L,f)$ (in $S^3$)
 is {near} to another one $(L',f')$ (of not necessarily the same number
of components)
if there is a finite set of \ASA\ $L''$ containing $L'$ such that:
\begin{equation}
[S^3, L, f] = \sum_{L''} [S^3, L'', f''] \in \Fas {\ast}
\end{equation}
We caution the reader that 
 being ``near to'' is not a symmetric  relation. Note also 
that any link is near to an arbitrary sublink of it.
A note on transitivity of the relation `` near to'': if $(L,f)$ is near to
$(L',f')$, which itself is near to a third one $(L'',f'')$, then it is
not clear that $(L,f)$ is near to $(L'',f'')$. {\em However}, if 
$(L'',f'')$ is obtained from $(L', f')$ (and $(L',f')$ is obtained from
$(L,f)$) using the equalities of figures \ref{f.onecross}, 
\ref{f.twocross}
and \ref{3bands}, then transitivity holds. We will use transitivity 
freely
in the proof of proposition \ref{prop.1}, because of the above note.
  As a variation, we call a finite linear combination of 
\ASA\ links $(L_i,f_i)$ near to an \ASA\ link 
link $(L',f')$  if there is a  finite number of \ASA\ links $(L'',f'')$
that include $(L',f')$ such that:
\begin{equation}
\sum_i [S^3, L_i, f_i] = \sum_{L''} [S^3, L'', f''] \in \Fas {\ast}
\end{equation}
\end{definition}

With the above terminology we have the following:
\begin{proposition}
\lbl{prop.1}
For every \ASA\ $4m$-component link in $S^3$ there is a trivial 
$L_{tr}(m)$
$m$-component link such that $(L,f)$ is near to $(L_{tr}(m),f(m))$.
\end{proposition}

\begin{pf}[of proposition \ref{prop.1}]
We begin by remarking that the statement in proposition \ref{prop.1}
is a finiteness statement, and not one using
downward induction. Furthermor, using figure \ref{plusminus} it 
follows
that if the above proposition holds for one choice of unit framings, 
then
it holds for all. We will therefore omit mentioning the framings in 
the
proof given below.

We divide the proof of the proposition  in two steps. 

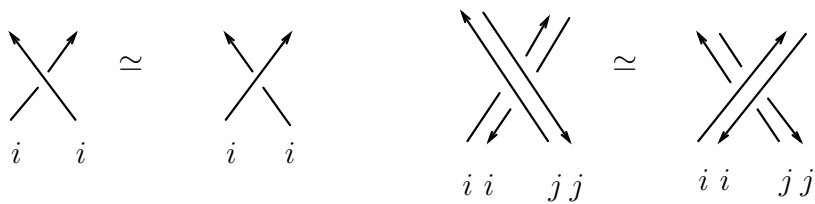
\begin{figure}[htpb]
$$\printname{sequiv}
	\setlength{\unitlength}{0.03\standardunitlength}
	\begin{array}{c}  \hspace{-1.7mm}
        	\raisebox{-8pt}{\begingroup\makeatletter\ifx\SetFigFont\undefined
% extract first six characters in \fmtname
\def\x#1#2#3#4#5#6#7\relax{\def\x{#1#2#3#4#5#6}}%
\expandafter\x\fmtname xxxxxx\relax \def\y{splain}%
\ifx\x\y   % LaTeX or SliTeX?
\gdef\SetFigFont#1#2#3{%
  \ifnum #1<17\tiny\else \ifnum #1<20\small\else
  \ifnum #1<24\normalsize\else \ifnum #1<29\large\else
  \ifnum #1<34\Large\else \ifnum #1<41\LARGE\else
     \huge\fi\fi\fi\fi\fi\fi
  \csname #3\endcsname}%
\else
\gdef\SetFigFont#1#2#3{\begingroup
  \count@#1\relax \ifnum 25<\count@\count@25\fi
  \def\x{\endgroup\@setsize\SetFigFont{#2pt}}%
  \expandafter\x
    \csname \romannumeral\the\count@ pt\expandafter\endcsname
    \csname @\romannumeral\the\count@ pt\endcsname
  \csname #3\endcsname}%
\fi
\fi\endgroup
\begin{picture}(11124,2613)(0,-10)
\thicklines
\path(3861.000,2155.500)(3912.000,2286.000)(3801.000,2200.500)
\path(3912,2286)(3012,1086)
\path(3120.982,2197.942)(3012.000,2286.000)(3059.952,2154.350)
\path(3012,2286)(3387,1761)
\path(912,1086)(12,2286)
\path(123.000,2200.500)(12.000,2286.000)(63.000,2155.500)
\path(537,1761)(912,2286)
\path(864.048,2154.350)(912.000,2286.000)(803.018,2197.942)
\path(12,1086)(387,1536)
\path(3537,1536)(3912,1011)
\path(7212,1986)(7512,2511)
\path(7477.580,2375.182)(7512.000,2511.000)(7412.462,2412.392)
\path(9612,786)(10812,2286)
\path(10756.949,2157.157)(10812.000,2286.000)(10698.384,2204.009)
\path(6418.086,2494.474)(6312.000,2586.000)(6355.683,2452.872)
\path(6312,2586)(7512,786)
\path(9718.086,2194.474)(9612.000,2286.000)(9655.683,2152.872)
\path(9612,2286)(10062,1611)
\path(6837,1461)(6387,786)
\path(9912,2286)(10212,1836)
\path(7362,1761)(7812,2511)
\path(6612,2586)(7812,786)
\path(7705.914,877.526)(7812.000,786.000)(7768.317,919.128)
\path(10587,1386)(11037,786)
\path(10926.000,871.500)(11037.000,786.000)(10986.000,916.500)
\path(9967.051,914.843)(9912.000,786.000)(10025.616,867.991)
\path(9912,786)(11112,2286)
\path(6730.683,919.128)(6687.000,786.000)(6793.086,877.526)
\path(6687,786)(6987,1236)
\path(10437,1236)(10737,786)
\put(1512,1761){\makebox(0,0)[lb]{$\simeq$}}
\put(3012,486){\makebox(0,0)[lb]{$i$}}
\put(3837,486){\makebox(0,0)[lb]{$i$}}
\put(912,486){\makebox(0,0)[lb]{$i$}}
\put(12,486){\makebox(0,0)[lb]{$i$}}
\put(8412,1761){\makebox(0,0)[lb]{$\simeq$}}
\put(9612,111){\makebox(0,0)[lb]{$i$}}
\put(9912,111){\makebox(0,0)[lb]{$i$}}
\put(10737,111){\makebox(0,0)[lb]{$j$}}
\put(11037,111){\makebox(0,0)[lb]{$j$}}
\put(7812,36){\makebox(0,0)[lb]{$j$}}
\put(7512,36){\makebox(0,0)[lb]{$j$}}
\put(6612,36){\makebox(0,0)[lb]{$i$}}
\put(6312,36){\makebox(0,0)[lb]{$i$}}
\end{picture} }
        	\hspace{-1.9mm}
	\end{array}
 $$
\caption{Two local moves that generate the relation of surgery 
equivalence.
Here arcs labeled by the same letter ($i$ or $j$)
belong to the same link  
component.}\lbl{sequiv}
\end{figure}

\begin{itemize}
\item{{\bf Step 1}} 
\hspace{0.5cm \ref{3bands}} $L$ is near to a (finite 
linear combination of) $L(\Gamma)$
where $L(\Gamma)$ are links obtained by trivalent vertex oriented 
graphs
with $4m$ edges as in \cite{Ga}, \cite{GL1}. 
\end{itemize} 
 
\begin{pf}
The proof is   similar to the proof of \cite[theorem 4]{GL2}. 
Note that theorem \cite[theorem 4]{GL2} states that
$[S^3,L,f] \in \Gas {4m}$ can be written as a linear combination
of terms of the form $[S^3,L(\Gamma),f]$ where  
$\Gamma$ are 
trivalent vertex oriented graphs of $4m$ edges. In that theorem,
we first alter $L$ by a surgery equivalent one (where surgery
equivalence is the relation generated by the local moves of figure
\ref{sequiv}). Using a relation shown in graphical form in figure
\ref{3bands} we   then replace $L$ by a linear combination of
$L(\Gamma)$ for trivalent vertex oriented graphs of $4m$ edges.

Our present claim (of step $1$) follows from the same argument
(sketched above) that shows \cite[theorem 4]{GL2}, after we use 
the equations (shown in graphical notation) in figures 
\ref{f.onecross}, \ref{f.twocross} and \ref{3bands}. Note that
in figures \ref{f.onecross}, \ref{f.twocross} and \ref{3bands}
the extra terms that are present  give links that
{\em contain} the links obtained by trivalent graphs.
\end{pf}  

\begin{itemize}
\item{{\bf Step 2}} 
\hspace{0.5cm} If $L(\Gamma)$ is a link obtained by a trivalent 
graph
of $4m$ edges, then $L(\Gamma)$ is near a trivial $m$-component 
link $L_{tr}(m)$.
\end{itemize}

\begin{pf}
This essentially follows from lemma 3.4 of \cite{Ga}. For 
completeness,
we repeat
the argument here. Take a forest $Forest$ 
of $\Gamma$  containing all of the vertices of   $\Gamma$,
and consider $L(Forest)$. This is a sublink of $L(\Gamma)$ and 
therefore $L(\Gamma)$ is near to $L(Forest)$
and an Euler characteristic argument shows 
that
$L(Forest)$ has at least $m$ components.
\end{pf}
The proof of proposition \ref{prop.1} is complete.
\end{pf}

We also have the following  proposition, that depends crucially 
on the
existence of the $4$-term relation shown in equation \eqref{2TM}.

\begin{proposition}
\lbl{prop.tubing}
There is a positive constant $c \in \BN$ such that
for every \ASA\ link $(L,f)$ in $S^3$ that contains a sublink 
$L_{tr}(c m^3) \cup K$
with the following properties:
\begin{itemize}
\item
     $L_{tr}(c m^3)$ is a trivial link of $c m^3$ components which 
bounds
a disjoint union of discs $\cup_i D_i$.
\item
     Each disc $D_i$ intersects the knot $K$ in two points, and 
intersects
no other component of $L$.
\end{itemize}
then, we have that $[S^3,L,f] \in \Fb m$, and therefore (using remark 
\ref{rem.bbl}) that $[S^3,L,f] \in \Fbl m$ as well.
\end{proposition}

Before we give the proof of it, let us introduce one more definition
that will be useful in stating the proof. Recall (see e.g. \cite{B-N1}) 
the 
combinatorial notion of a {\em chord diagram} with support on a 
circle.

\begin{definition}
\lbl{def.bCD}
A chord diagram is called $m$-{\em boundary} if it contains $m$ 
nonintersecting
chords. For an example, see figure \ref{boundaryCD}.
\end{definition}

\begin{figure}[htpb]
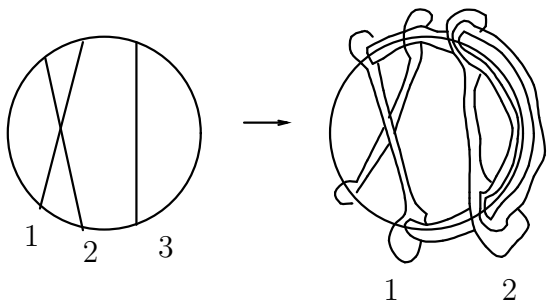

$$\printname{boundaryCD}
	\setlength{\unitlength}{0.03\standardunitlength}
	\begin{array}{c}  \hspace{-1.7mm}
        	\raisebox{-8pt}{\input draws/boundaryCD.tex }
        	\hspace{-1.9mm}
	\end{array}
 $$
\caption{Shown on the left is a chord diagram with $3$ chords. Two 
of them
($1$ and $3$) do not intersect thus the chord diagram is  $2$-
boundary. 
On the
right is shown the result of tubing it along the two nonintersecting 
chords.
The result is a boundary link of two components in the complement 
of a
two component link.}\lbl{boundaryCD}
\end{figure}

\begin{pf}
Let $L$ and $K$ be as in the statement of proposition 
\ref{prop.tubing},
and let $CD_K$ denote the associated chord diagram of $K$, .i.e., the
ordered set of points of the intersection of $K$ with the disk $D_i$.
The two points of intersection of $D_i \cap K$ will be represented 
by
a chord of $CD_K$, as in the theory of \fti s of knots.

We now give the proof in six steps:

\begin{itemize}
\item{{\bf Step 1}} 
\hspace{0.5cm} If $CD_K$ is $m$-boundary, then $[S^3,L,f] \in \Fb 
m$.
\end{itemize}
 
\begin{pf}
Indeed, tube the $m$ discs that represent $m$ nonintersecting 
chords,
using an innermost circle argument. For an example, see figure 
\ref{boundaryCD}.
The result is a boundary $m$ component sublink $L_b$ of $L$ that 
bounds 
surfaces in the {\em complement} of $L-L_b$, from which it follows 
easily that $[S^3,L,f] \in \Fb m$. Note that the result is independent
of the fact that $K$ may be knotted; it only depends on the 
associated
chord diagram of $K$.
\end{pf}

Now, if $CD_K$ is not an $m$-boundary chord diagram,
we have the following steps.

\begin{itemize}
\item{{\bf Step 2}} 
\hspace{0.5cm} We can always assume that the framing in each of 
the 
components of $L_{tr}(c m^4)$ is $+1$ or $-1$ as we please.
\end{itemize}

\begin{pf}
Use the equation in figure \ref{plusminus}.
\end{pf}

Before we state the next step, let us introduce some useful 
terminology.
We say that a chord diagram  is {\em represented by an $m$-tower} 
if it can 
be written
{\em using the $4$-term relation} (see figure \ref{4Tv1})
 as a  (finite)
 linear combination of $m$-boundary chord diagrams.
With this terminology we have the next step:

\begin{itemize}
\item{{\bf Step 3}} 
\hspace{0.5cm} If a chord diagram contains at least $c m^3$ chords 
(where
$c > 1$ fixed positive integer), and is represented by an $m-1$
tower, then it is represented by an $m$-tower.
\end{itemize}

\begin{pf}
Without loss of generality, we may assume that the chord diagram 
contains an $m-1$ tower, $T_{m-1}$. 
The end points of the $m-1$ chords of $T_{m-1}$ partition 
the
external circle of the chord diagram in $2(m-1)$ arcs, see figure 
\ref{bcircle}. Each of the rest of the chords of the chord diagram
will
begin and end in one of these arcs. If one of these chords begins and
ends in the same arc, then it, together with $T_{m-1}$, is a set of
$m$ nonintersecting chords. If not, there are $\binom{2(m-
1)}{2}$
many possibilities for the beginnings and ends of the extra arcs. 
Using
the pigeonhole principle, if 
\begin{equation}
\lbl{eq.pig}
c m^3 - (m-1) > \binom{2(m-1)}{2} (m-1) + 1
\end{equation}
we can always assume that there are $m$ many of the rest of the 
chords 
that all
begin in one arc, and end in another.
Let $\alpha, \beta$ be two arcs in which (at least) $m$ of the rest 
of the
chords 
begin and end. Let us declare one of these  chords to be 
{\em``special''}
if it has one end in $\alpha$ and the other in $\beta$;
 otherwise declare  it 
{\em ``nonspecial''}. 

Let us look at the  chords that begin in the arc $\alpha$.
According to step $4$ below, we can always move a ``special'' over a
``nonspecial'' one, and according to step $5$ below, we can always 
move
a ``special'' over another ``special''. After doing so, we may assume
that the ``special'' $m$ chords do not intersect, and therefore, form
an $m$-tower. 
\end{pf}

\begin{figure}[htpb]
$$\printname{bcircle}
	\setlength{\unitlength}{0.025\standardunitlength}
	\begin{array}{c}  \hspace{-1.7mm}
        	\raisebox{-8pt}{\begingroup\makeatletter\ifx\SetFigFont\undefined
% extract first six characters in \fmtname
\def\x#1#2#3#4#5#6#7\relax{\def\x{#1#2#3#4#5#6}}%
\expandafter\x\fmtname xxxxxx\relax \def\y{splain}%
\ifx\x\y   % LaTeX or SliTeX?
\gdef\SetFigFont#1#2#3{%
  \ifnum #1<17\tiny\else \ifnum #1<20\small\else
  \ifnum #1<24\normalsize\else \ifnum #1<29\large\else
  \ifnum #1<34\Large\else \ifnum #1<41\LARGE\else
     \huge\fi\fi\fi\fi\fi\fi
  \csname #3\endcsname}%
\else
\gdef\SetFigFont#1#2#3{\begingroup
  \count@#1\relax \ifnum 25<\count@\count@25\fi
  \def\x{\endgroup\@setsize\SetFigFont{#2pt}}%
  \expandafter\x
    \csname \romannumeral\the\count@ pt\expandafter\endcsname
    \csname @\romannumeral\the\count@ pt\endcsname
  \csname #3\endcsname}%
\fi
\fi\endgroup
\begin{picture}(6544,2367)(0,-10)
\thicklines
\put(487.500,2269.000){\arc{2493.116}{0.1206}{1.2966}}
\put(1633.594,-252.875){\arc{1954.327}{3.9728}{5.7498}}
\put(3120.259,1295.293){\arc{1653.657}{2.0016}{3.9513}}
\put(4087.500,2269.000){\arc{2493.116}{0.1206}{1.2966}}
\put(6720.259,1295.293){\arc{1653.657}{2.0016}{3.9513}}
\put(5308.594,-252.875){\arc{1954.327}{3.9728}{5.7498}}
\put(4150.338,1374.068){\arc{1664.201}{5.6083}{7.4203}}
\put(3764.062,469.000){\arc{3105.165}{5.3197}{6.5772}}
\put(4776.884,2462.664){\arc{2111.247}{0.3316}{1.8362}}
\put(3915.726,1965.371){\arc{3133.648}{6.1850}{7.4197}}
\path(1613.464,2290.524)(1725.000,2344.000)(1601.421,2349.303)
\path(653.340,1265.782)(675.000,1144.000)(713.197,1261.648)
\put(1989.474,1053.210){\arc{2635.212}{3.2106}{4.5103}}
\path(809.358,478.167)(900.000,394.000)(859.630,510.920)
\path(721.463,879.365)(675.000,994.000)(662.051,870.986)
\put(2137.500,1200.250){\arc{2953.943}{2.5642}{3.0015}}
\put(1875,1069){\ellipse{2122}{2122}}
\put(5475,1069){\ellipse{2122}{2122}}
\put(0,544){\makebox(0,0)[lb]{$\alpha$}}
\put(150,1894){\makebox(0,0)[lb]{$\beta$}}
\end{picture} }
        	\hspace{-1.9mm}
	\end{array}
 $$
\caption{Shown on the left is a $3$ tower, and two arcs $\alpha, 
\beta$ of 
the external circle of the chord diagram. Shown on the right are 
some
of the extra chords.}\lbl{bcircle}
\end{figure}
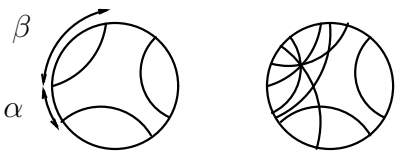

\begin{itemize}
\item{{\bf Step 4}} 
\hspace{0.5cm} With the notation of  step $3$, we can always move
a ``special'' chord over a ``nonspecial'' one. 
\end{itemize}

\begin{pf}
The proof uses the $4$-term relation. Fix a ``special'' chord and
move a ``nonspecial'' in four ways around the end of the ``special'' 
one.
Of the resulting four terms, two of them no longer have the
``nonspecial'' chord, and the two others move the ``special'' one
over the nonspecial one.  See figure \ref{4Tmove1}.  
\end{pf}

\begin{figure}[htpb]
$$\printname{4Tmove1}
	\setlength{\unitlength}{0.02\standardunitlength}
	\begin{array}{c}  \hspace{-1.7mm}
        	\raisebox{-8pt}{\begingroup\makeatletter\ifx\SetFigFont\undefined
% extract first six characters in \fmtname
\def\x#1#2#3#4#5#6#7\relax{\def\x{#1#2#3#4#5#6}}%
\expandafter\x\fmtname xxxxxx\relax \def\y{splain}%
\ifx\x\y   % LaTeX or SliTeX?
\gdef\SetFigFont#1#2#3{%
  \ifnum #1<17\tiny\else \ifnum #1<20\small\else
  \ifnum #1<24\normalsize\else \ifnum #1<29\large\else
  \ifnum #1<34\Large\else \ifnum #1<41\LARGE\else
     \huge\fi\fi\fi\fi\fi\fi
  \csname #3\endcsname}%
\else
\gdef\SetFigFont#1#2#3{\begingroup
  \count@#1\relax \ifnum 25<\count@\count@25\fi
  \def\x{\endgroup\@setsize\SetFigFont{#2pt}}%
  \expandafter\x
    \csname \romannumeral\the\count@ pt\expandafter\endcsname
    \csname @\romannumeral\the\count@ pt\endcsname
  \csname #3\endcsname}%
\fi
\fi\endgroup
\begin{picture}(13735,2620)(0,-10)
\thicklines
\put(-16.602,2804.929){\arc{3301.319}{0.2660}{1.2359}}
\put(1502.210,-470.532){\arc{2237.935}{4.0635}{5.6451}}
\put(3563.500,1452.250){\arc{2787.388}{2.3086}{3.7259}}
\put(3583.398,2804.929){\arc{3301.319}{0.2660}{1.2359}}
\put(5102.210,-470.532){\arc{2237.935}{4.0635}{5.6451}}
\put(7163.500,1452.250){\arc{2787.388}{2.3086}{3.7259}}
\put(7183.398,2804.929){\arc{3301.319}{0.2660}{1.2359}}
\put(8702.210,-470.532){\arc{2237.935}{4.0635}{5.6451}}
\put(10763.500,1452.250){\arc{2787.388}{2.3086}{3.7259}}
\put(10783.398,2804.929){\arc{3301.319}{0.2660}{1.2359}}
\put(12302.210,-470.532){\arc{2237.935}{4.0635}{5.6451}}
\put(14363.500,1452.250){\arc{2787.388}{2.3086}{3.7259}}
\put(-78.323,2004.120){\arc{2769.304}{6.0701}{7.3420}}
\put(3521.677,2004.120){\arc{2769.304}{6.0701}{7.3420}}
\put(7121.677,2004.120){\arc{2769.304}{6.0701}{7.3420}}
\put(10721.678,2004.120){\arc{2769.303}{6.0701}{7.3420}}
\put(9825.000,2972.000){\arc{3354.102}{1.3909}{2.6779}}
\put(13338.462,2366.231){\arc{2226.954}{1.2163}{3.1468}}
\path(641.937,332.187)(750.000,272.000)(682.974,375.959)
\path(401.621,1051.205)(375.000,1172.000)(341.645,1052.889)
\put(1562.500,1138.667){\arc{2375.935}{2.3239}{3.1697}}
\path(344.607,1441.901)(375.000,1322.000)(404.607,1442.098)
\path(1456.916,2560.175)(1575.000,2597.000)(1453.479,2620.077)
\put(1647.917,1326.167){\arc{2545.846}{3.1383}{4.6551}}
\put(1726,1209){\ellipse{2402}{2402}}
\put(5326,1209){\ellipse{2402}{2402}}
\put(8926,1209){\ellipse{2402}{2402}}
\put(12526,1209){\ellipse{2402}{2402}}
\path(526,946)(2926,1396)
\path(4275,722)(6525,1472)
\put(0,572){\makebox(0,0)[lb]{$\alpha$}}
\put(225,2372){\makebox(0,0)[lb]{$\beta$}}
\put(3300,1322){\makebox(0,0)[lb]{$-$}}
\put(6975,1322){\makebox(0,0)[lb]{$=$}}
\put(10575,1322){\makebox(0,0)[lb]{$-$}}
\end{picture} }
        	\hspace{-1.9mm}
	\end{array}
 $$
\caption{The $4$-term relation with a fixed ``special'' chord and a 
moving
``nonspecial'' chord. For the convenience, the arcs $\alpha, \beta$
are shown, too. The first two terms show a ``nonspecial'' chord after 
and before passing a special chord (in the $\alpha$ arc). The two last 
terms have 
no ``nonspecial'' chord. }\lbl{4Tmove1}
\end{figure}
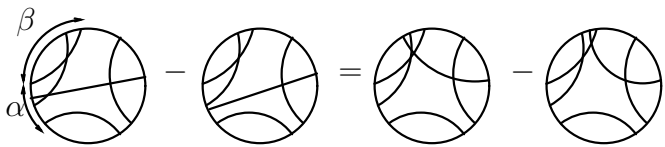

\begin{itemize}
\item{{\bf Step 5}} 
\hspace{0.5cm} With the notation of  step $3$, we can always move
a ``special'' chord over a   ``special'' one. 
\end{itemize}

\begin{pf}
The same as in step $4$, see figure \ref{4Tmove2}.
\end{pf}

\begin{itemize}
\item{{\bf Step 6}} 
\hspace{0.5cm} Induction.
\end{itemize}

We can now finish the proof of proposition \ref{prop.1} as follows:
obviously, a chord diagram contains (and therefore, is represented 
by)
a $1$-tower. Using step $3$, and induction, proposition \ref{prop.1}
follows.
\end{pf}

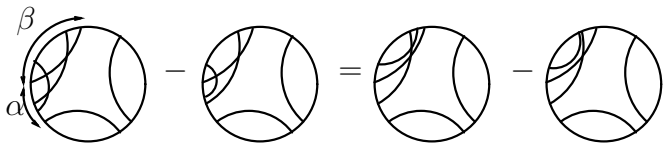
\begin{figure}[htpb]
$$\printname{4Tmove2}
	\setlength{\unitlength}{0.02\standardunitlength}
	\begin{array}{c}  \hspace{-1.7mm}
        	\raisebox{-8pt}{\begingroup\makeatletter\ifx\SetFigFont\undefined
% extract first six characters in \fmtname
\def\x#1#2#3#4#5#6#7\relax{\def\x{#1#2#3#4#5#6}}%
\expandafter\x\fmtname xxxxxx\relax \def\y{splain}%
\ifx\x\y   % LaTeX or SliTeX?
\gdef\SetFigFont#1#2#3{%
  \ifnum #1<17\tiny\else \ifnum #1<20\small\else
  \ifnum #1<24\normalsize\else \ifnum #1<29\large\else
  \ifnum #1<34\Large\else \ifnum #1<41\LARGE\else
     \huge\fi\fi\fi\fi\fi\fi
  \csname #3\endcsname}%
\else
\gdef\SetFigFont#1#2#3{\begingroup
  \count@#1\relax \ifnum 25<\count@\count@25\fi
  \def\x{\endgroup\@setsize\SetFigFont{#2pt}}%
  \expandafter\x
    \csname \romannumeral\the\count@ pt\expandafter\endcsname
    \csname @\romannumeral\the\count@ pt\endcsname
  \csname #3\endcsname}%
\fi
\fi\endgroup
\begin{picture}(13735,2620)(0,-10)
\thicklines
\put(-16.602,2804.929){\arc{3301.319}{0.2660}{1.2359}}
\put(1502.210,-470.532){\arc{2237.935}{4.0635}{5.6451}}
\put(3563.500,1452.250){\arc{2787.388}{2.3086}{3.7259}}
\put(3583.398,2804.929){\arc{3301.319}{0.2660}{1.2359}}
\put(5102.210,-470.532){\arc{2237.935}{4.0635}{5.6451}}
\put(7163.500,1452.250){\arc{2787.388}{2.3086}{3.7259}}
\put(7183.398,2804.929){\arc{3301.319}{0.2660}{1.2359}}
\put(8702.210,-470.532){\arc{2237.935}{4.0635}{5.6451}}
\put(10763.500,1452.250){\arc{2787.388}{2.3086}{3.7259}}
\put(10783.398,2804.929){\arc{3301.319}{0.2660}{1.2359}}
\put(12302.210,-470.532){\arc{2237.935}{4.0635}{5.6451}}
\put(14363.500,1452.250){\arc{2787.388}{2.3086}{3.7259}}
\put(-78.323,2004.120){\arc{2769.304}{6.0701}{7.3420}}
\put(3521.677,2004.120){\arc{2769.304}{6.0701}{7.3420}}
\put(7121.677,2004.120){\arc{2769.304}{6.0701}{7.3420}}
\put(10721.678,2004.120){\arc{2769.303}{6.0701}{7.3420}}
\put(259.500,1145.000){\arc{1297.142}{5.2651}{7.1587}}
\put(4062.500,1284.500){\arc{728.869}{5.0993}{7.4671}}
\put(7868.750,2375.125){\arc{1512.513}{0.0041}{1.6618}}
\put(11471.739,2104.609){\arc{1124.410}{5.9339}{7.9819}}
\path(335.042,1364.061)(375.000,1247.000)(394.831,1369.093)
\path(1456.074,2562.992)(1575.000,2597.000)(1454.061,2622.958)
\put(1616.803,1351.508){\arc{2492.387}{3.0576}{4.6788}}
\path(564.549,402.683)(675.000,347.000)(603.747,448.109)
\path(393.134,974.643)(375.000,1097.000)(333.421,980.505)
\put(1248.214,1011.286){\arc{1754.821}{2.2827}{3.2394}}
\put(1726,1209){\ellipse{2402}{2402}}
\put(5326,1209){\ellipse{2402}{2402}}
\put(8926,1209){\ellipse{2402}{2402}}
\put(12526,1209){\ellipse{2402}{2402}}
\put(0,572){\makebox(0,0)[lb]{$\alpha$}}
\put(225,2372){\makebox(0,0)[lb]{$\beta$}}
\put(3300,1322){\makebox(0,0)[lb]{$-$}}
\put(6975,1322){\makebox(0,0)[lb]{$=$}}
\put(10575,1322){\makebox(0,0)[lb]{$-$}}
\end{picture} }
        	\hspace{-1.9mm}
	\end{array}
 $$
\caption{The $4$-term relation with a fixed ``special'' chord and a 
moving
``special'' chord. For the convenience, the arcs $\alpha, \beta$
are shown, too. The first two terms show a ``special'' chord after and 
before passing a fixed ``special'' chord (in the $\alpha$ arc). 
The two last terms have chords that begin and end in the $\beta$ 
arc, and
therefore, by the discussion of step $3$, we can find a 
$4$-tower. }\lbl{4Tmove2}
\end{figure}

We also have the following proposition, similar to, but different
from proposition \ref{prop.tubing}:

\begin{proposition}
\lbl{prop.moretubing}
There is an increasing function $h: \BN \to \BN$ with the following
property:
for every \ASA\ link $L$  in $S^3$ (with framing $f$) which
contains a sublink $L_{tr}(h(m)) \cup L'$
with the following properties:
\begin{itemize}
\item
     $L_{tr}(h(m))$ is a trivial link of $h(m)$ components which 
bounds
a disjoint union of discs $\cup_i D_i$.
\item
     Each disc $D_i$ intersects the link $L'$ in either two or
four points. Moreover, the intersections of the 
disc $D_i$ with $L'$ come in pairs with opposite orientation for
each component of $L'$. Furthermore, in case a disc $D_i$ intersects
$L'$ in four points, we assume that these four points do not lie 
in the same component of $L'$. See figure \ref{Dintersection}.
\end{itemize}
Then  we have that $[S^3,L,f] \in \Fbl m$.
\end{proposition}

\begin{figure}[htpb]
$$\printname{Dintersection}
	\setlength{\unitlength}{0.03\standardunitlength}
	\begin{array}{c}  \hspace{-1.7mm}
        	\raisebox{-8pt}{\begingroup\makeatletter\ifx\SetFigFont\undefined
% extract first six characters in \fmtname
\def\x#1#2#3#4#5#6#7\relax{\def\x{#1#2#3#4#5#6}}%
\expandafter\x\fmtname xxxxxx\relax \def\y{splain}%
\ifx\x\y   % LaTeX or SliTeX?
\gdef\SetFigFont#1#2#3{%
  \ifnum #1<17\tiny\else \ifnum #1<20\small\else
  \ifnum #1<24\normalsize\else \ifnum #1<29\large\else
  \ifnum #1<34\Large\else \ifnum #1<41\LARGE\else
     \huge\fi\fi\fi\fi\fi\fi
  \csname #3\endcsname}%
\else
\gdef\SetFigFont#1#2#3{\begingroup
  \count@#1\relax \ifnum 25<\count@\count@25\fi
  \def\x{\endgroup\@setsize\SetFigFont{#2pt}}%
  \expandafter\x
    \csname \romannumeral\the\count@ pt\expandafter\endcsname
    \csname @\romannumeral\the\count@ pt\endcsname
  \csname #3\endcsname}%
\fi
\fi\endgroup
\begin{picture}(5116,1903)(0,-10)
\thicklines
\put(908,457){\ellipse{1800}{900}}
\put(4208,457){\ellipse{1800}{900}}
\path(308,1507)(308,307)
\path(278.000,427.000)(308.000,307.000)(338.000,427.000)
\path(1208,1507)(1208,307)
\path(1178.000,427.000)(1208.000,307.000)(1238.000,427.000)
\path(3908,1507)(3908,307)
\path(3878.000,427.000)(3908.000,307.000)(3938.000,427.000)
\path(638.000,1387.000)(608.000,1507.000)(578.000,1387.000)
\path(608,1507)(608,307)
\path(1538.000,1387.000)(1508.000,1507.000)(1478.000,1387.000)
\path(1508,1507)(1508,307)
\path(4538.000,1387.000)(4508.000,1507.000)(4478.000,1387.000)
\path(4508,1507)(4508,307)
\put(458,1732){\makebox(0,0)[lb]{$\alpha$}}
\put(3608,1732){\makebox(0,0)[lb]{$\alpha$}}
\put(4283,1732){\makebox(0,0)[lb]{$\alpha$}}
\put(983,1732){\makebox(0,0)[lb]{$\beta$}}
\put(1433,1732){\makebox(0,0)[lb]{$\beta$}}
\put(8,1732){\makebox(0,0)[lb]{$\alpha$}}
\end{picture} }
        	\hspace{-1.9mm}
	\end{array}
 $$
\caption{Shown here are the two allowable  types of intersections 
of  a disc
$D_i$ with the components of $L'$. Here $\alpha, \beta$ denote 
components of
$L'$ and we {\em assume} that $\alpha \neq 
\beta$.}\lbl{Dintersection}
\end{figure}
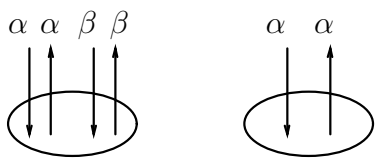

\begin{remark}
\lbl{rem.tubing}
Before we give the proof of the above proposition let us point out 
that 
the assumptions are weaker than those of proposition 
\ref{prop.tubing}.
As a result, the conclusion is weaker than that of proposition
\ref{prop.tubing}, in the sense that  $[S^3,L,f]$ lies
in $\Fbl m$ and not necessarily in $\Fb m$. Note also that the proof 
of proposition \ref{prop.moretubing} shows that the function $h$ is
constructible, e.g. we can take $h(m)=c m^{13}$ for some constant $c$.
\end{remark}

\begin{pf}[of proposition \ref{prop.moretubing}]

Let $L, L'$ be as in the statement of proposition 
\ref{prop.moretubing}. We begin
by introducing the associated {\em chord diagram} $CD_{L'}$ of $L'$ 
relative to the
union of discs $D_i$. The chord diagram $CD_{L'}$ consists of 
external circles
(one per component of $L'$) and chords (as many as the number of 
discs $D_i$).
There are two types of chords: ones that intersect the external 
circles in two points, (called
of type $I$) and
the ones (called of type $II$)
that intersect the external circles in four points (however two 
points are in one
circle and two are in another). For an example see figure 
\ref{newCD}. Note that these
chord diagrams are similar but different from  the chord diagrams 
on links.
At any rate, they include as a special case the chord diagrams 
considered in 
proposition \ref{prop.tubing}.

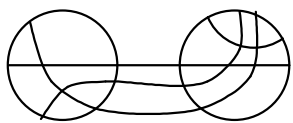
\begin{figure}[htpb]
$$\printname{newCD}
	\setlength{\unitlength}{0.03\standardunitlength}
	\begin{array}{c}  \hspace{-1.7mm}
        	\raisebox{-8pt}{\begingroup\makeatletter\ifx\SetFigFont\undefined
% extract first six characters in \fmtname
\def\x#1#2#3#4#5#6#7\relax{\def\x{#1#2#3#4#5#6}}%
\expandafter\x\fmtname xxxxxx\relax \def\y{splain}%
\ifx\x\y   % LaTeX or SliTeX?
\gdef\SetFigFont#1#2#3{%
  \ifnum #1<17\tiny\else \ifnum #1<20\small\else
  \ifnum #1<24\normalsize\else \ifnum #1<29\large\else
  \ifnum #1<34\Large\else \ifnum #1<41\LARGE\else
     \huge\fi\fi\fi\fi\fi\fi
  \csname #3\endcsname}%
\else
\gdef\SetFigFont#1#2#3{\begingroup
  \count@#1\relax \ifnum 25<\count@\count@25\fi
  \def\x{\endgroup\@setsize\SetFigFont{#2pt}}%
  \expandafter\x
    \csname \romannumeral\the\count@ pt\expandafter\endcsname
    \csname @\romannumeral\the\count@ pt\endcsname
  \csname #3\endcsname}%
\fi
\fi\endgroup
\begin{picture}(3946,1561)(0,-10)
\thicklines
\put(3443.000,1718.000){\arc{1398.463}{0.9530}{2.7452}}
\put(773,773){\ellipse{1530}{1530}}
\put(3173,773){\ellipse{1530}{1530}}
\path(23,773)(3923,773)
\path(323,1373)	(333.958,1328.175)
	(344.676,1284.957)
	(355.169,1243.307)
	(365.450,1203.188)
	(375.532,1164.562)
	(385.430,1127.391)
	(404.728,1057.264)
	(423.452,992.504)
	(441.712,932.809)
	(459.619,877.878)
	(477.283,827.409)
	(494.812,781.098)
	(512.318,738.644)
	(547.696,664.100)
	(584.297,601.357)
	(623.000,548.000)

\path(623,548)	(674.780,493.176)
	(741.337,436.430)
	(778.884,407.991)
	(818.599,379.858)
	(859.972,352.293)
	(902.495,325.558)
	(945.659,299.914)
	(988.954,275.625)
	(1031.873,252.951)
	(1073.906,232.156)
	(1114.544,213.501)
	(1153.278,197.249)
	(1223.000,173.000)

\path(1223,173)	(1295.075,154.822)
	(1334.425,146.592)
	(1375.779,138.940)
	(1418.983,131.868)
	(1463.883,125.380)
	(1510.326,119.480)
	(1558.157,114.171)
	(1607.223,109.456)
	(1657.370,105.340)
	(1708.445,101.824)
	(1760.294,98.913)
	(1812.763,96.611)
	(1865.697,94.920)
	(1918.945,93.845)
	(1972.351,93.387)
	(2025.763,93.552)
	(2079.025,94.343)
	(2131.986,95.762)
	(2184.490,97.813)
	(2236.384,100.500)
	(2287.515,103.826)
	(2337.728,107.795)
	(2386.870,112.410)
	(2434.788,117.675)
	(2481.326,123.592)
	(2526.333,130.166)
	(2569.653,137.399)
	(2611.134,145.296)
	(2650.621,153.859)
	(2723.000,173.000)

\path(2723,173)	(2761.469,185.624)
	(2803.271,201.207)
	(2847.799,219.558)
	(2894.443,240.485)
	(2942.596,263.797)
	(2991.651,289.303)
	(3040.998,316.812)
	(3090.031,346.134)
	(3138.141,377.076)
	(3184.720,409.448)
	(3229.160,443.058)
	(3270.854,477.716)
	(3343.568,549.410)
	(3398.000,623.000)

\path(3398,623)	(3430.819,688.068)
	(3456.500,761.414)
	(3466.733,802.016)
	(3475.264,845.676)
	(3482.120,892.724)
	(3487.329,943.490)
	(3490.918,998.303)
	(3492.915,1057.493)
	(3493.348,1121.389)
	(3492.243,1190.321)
	(3489.628,1264.618)
	(3487.763,1303.882)
	(3485.531,1344.611)
	(3482.935,1386.846)
	(3479.979,1430.628)
	(3476.666,1475.999)
	(3473.000,1523.000)

\path(473,23)	(505.727,73.721)
	(536.575,120.660)
	(565.680,163.982)
	(593.180,203.852)
	(643.917,273.894)
	(689.881,332.105)
	(732.173,379.803)
	(771.891,418.305)
	(810.134,448.932)
	(848.000,473.000)

\path(848,473)	(900.748,495.490)
	(964.649,512.338)
	(1035.966,524.490)
	(1110.961,532.895)
	(1185.899,538.501)
	(1257.043,542.255)
	(1320.655,545.105)
	(1373.000,548.000)

\path(1373,548)	(1438.790,550.199)
	(1513.203,548.615)
	(1553.273,546.637)
	(1595.056,543.994)
	(1638.403,540.779)
	(1683.167,537.086)
	(1729.200,533.008)
	(1776.353,528.638)
	(1824.479,524.070)
	(1873.430,519.398)
	(1923.058,514.714)
	(1973.216,510.114)
	(2023.755,505.689)
	(2074.528,501.534)
	(2125.386,497.742)
	(2176.182,494.406)
	(2226.768,491.620)
	(2276.997,489.478)
	(2326.719,488.073)
	(2375.788,487.498)
	(2424.056,487.848)
	(2471.374,489.215)
	(2517.595,491.692)
	(2562.571,495.375)
	(2606.154,500.355)
	(2648.196,506.727)
	(2688.550,514.584)
	(2727.067,524.020)
	(2798.000,548.000)

\path(2798,548)	(2860.166,579.690)
	(2925.969,623.724)
	(2992.745,677.438)
	(3057.830,738.170)
	(3118.562,803.255)
	(3172.276,870.031)
	(3216.310,935.834)
	(3248.000,998.000)

\path(3248,998)	(3261.335,1036.598)
	(3270.860,1079.841)
	(3276.575,1129.269)
	(3278.480,1186.419)
	(3276.575,1252.829)
	(3274.194,1289.987)
	(3270.860,1330.037)
	(3266.574,1373.171)
	(3261.335,1419.581)
	(3255.144,1469.460)
	(3248.000,1523.000)

\end{picture} }
        	\hspace{-1.9mm}
	\end{array}
 $$
\caption{An example of a chord diagram considered in proposition 
\ref{prop.moretubing}.}\lbl{newCD}
\end{figure}

Two chords  intersect if there is an external circle that
they both touch, such that the four intersection points of that circle 
with
the two chords is in the order $1212$. See figure \ref{CDintersect}. 
In analogy to definition \ref{def.bCD} we call a chord diagram $m$-
{\em boundary} if it
contains $m$ nonintersecting chords. 
The motivation for considering nonintersecting chords, is the fact 
that they can be tubed,
and therefore produce boundary links in the complement of $L'$. 

\begin{figure}[htpb]
$$\printname{CDintersect}
	\setlength{\unitlength}{0.03\standardunitlength}
	\begin{array}{c}  \hspace{-1.7mm}
        	\raisebox{-8pt}{\begingroup\makeatletter\ifx\SetFigFont\undefined
% extract first six characters in \fmtname
\def\x#1#2#3#4#5#6#7\relax{\def\x{#1#2#3#4#5#6}}%
\expandafter\x\fmtname xxxxxx\relax \def\y{splain}%
\ifx\x\y   % LaTeX or SliTeX?
\gdef\SetFigFont#1#2#3{%
  \ifnum #1<17\tiny\else \ifnum #1<20\small\else
  \ifnum #1<24\normalsize\else \ifnum #1<29\large\else
  \ifnum #1<34\Large\else \ifnum #1<41\LARGE\else
     \huge\fi\fi\fi\fi\fi\fi
  \csname #3\endcsname}%
\else
\gdef\SetFigFont#1#2#3{\begingroup
  \count@#1\relax \ifnum 25<\count@\count@25\fi
  \def\x{\endgroup\@setsize\SetFigFont{#2pt}}%
  \expandafter\x
    \csname \romannumeral\the\count@ pt\expandafter\endcsname
    \csname @\romannumeral\the\count@ pt\endcsname
  \csname #3\endcsname}%
\fi
\fi\endgroup
\begin{picture}(10220,2434)(0,-10)
\thicklines
\put(608,1807){\ellipse{1200}{1200}}
\put(4208,607){\ellipse{1200}{1200}}
\put(3008,1807){\ellipse{1200}{1200}}
\put(5408,1807){\ellipse{1200}{1200}}
\put(7808,1807){\ellipse{1200}{1200}}
\put(9608,1807){\ellipse{1200}{1200}}
\path(158,2182)(983,1357)
\path(1058,2182)(158,1432)
\path(2558,2182)(4658,232)
\path(5858,2182)(3758,232)
\path(7208,1807)(10208,1807)
\path(7808,2407)(7808,1207)
\end{picture} }
        	\hspace{-1.9mm}
	\end{array}
 $$
\caption{Examples of intersecting chords of various types: $(I,I)$, 
$(II,II)$
and $(I,II)$. }\lbl{CDintersect}
\end{figure}
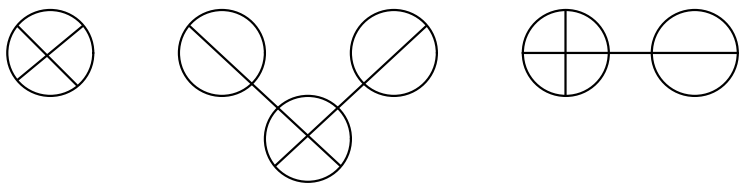

This shows the first step in the proof of proposition 
\ref{prop.moretubing}:

\begin{itemize}
\item{{\bf Step 1}} 
\hspace{0.5cm} If $CD_{L'}$ is $m$-boundary, then $[S^3,L,f] \in \Fb 
m$, and therefore
(using remark \ref{rem.bbl}) $[S^3,L,f] \in \Fbl m$
\end{itemize}

The rest of the proof will be devoted to the proof that we can 
assume the 
hypothesis as in step $1$. It uses, like proposition, \ref{prop.tubing} 
the 
$4$-term relation in a crucial way. We sketch the proof here:

\begin{itemize}
\item{{\bf Step 2}} 
\hspace{0.5cm} We can always assume that the framing in each of 
the 
components of $L_{tr}$ is $+1$ or $-1$ as we please.
\end{itemize}

Indeed, see step $2$ of proposition \ref{prop.tubing}.
We can now define (in direct analogy with proposition 
\ref{prop.tubing}) the notion of 
a chord diagram containing an $m$-{tower}. Let us concentrate on 
the
chord diagram $CD_{L'}$.  With the above terminology we have the 
following:

\begin{itemize}
\item{{\bf Step 3}} 
\hspace{0.5cm} If the chord diagram $CD_{L'}$
 has an external circle which touches
 $g_0(m) >> m$ (for some function $g_0$), then it contains an 
$m$-tower,
and therefore (by tubing) we conclude that $[S^3,L,f] \in \Fbl m$.
\end{itemize}

Before we give the proof, let us point out that the function
$g_0$ above, (and the functions $g_1, g_2, h$ to be introduced 
later), can be constructed explicitly.

\begin{pf}
The proof uses the first, second and third version of the $4$-term 
relation 
on $\M$. 
(see figures \ref{4Tv2} and \ref{4Tv3}). 
Consider such an external circle, say $C_1$. 
{\em Ignoring the subleading terms}, apply the first, second and third 
version
of the $4$-term relation (just as in step 3 of proposition 
\ref{prop.tubing})
in order to get $g_1(m)$ many chords, where $g_0 (m)>cg_1 
(m)^3$, as in Step 3 of proposition~\ref{prop.tubing},
which are nonintersecting as far as their ends in the external circle
$C_1$ are concerned. Concentrate on the $g_1(m)$ many chords that 
touch
the external circle $C_1$. Call these chords {\em preferred}.
 Consider all other external circles that these
chords touch. 
Now we consider two cases:
\begin{description}
\item[Case 1]
There are at least $m$ such other external circles.
\end{description}
 Then we can create,
using at least $m$ (of the $g_1(m)$ many preferred 
chords) an $m$-tower. See figure \ref{mtower1}.
\begin{description}
\item[Case 2]
Assume there are at most $m$ such external circles. 
\end{description}
Then, if $g_1(m) > m g_2(m)$, there is at least one circle $C_2$ 
containing
$g_2(m)$ many preferred chords (that lie on the circle $C_1$ and 
$C_2$).
{\em Ignoring} the subleading terms once again, by applying the 
$4$-term relation, we can reach a linear combination of chord 
diagrams
with a $3m$-tower provided $g_2 (m)>c(3m)^3$.

\begin{figure}[htpb]
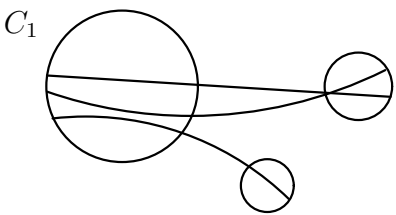

$$\printname{mtower1}
	\setlength{\unitlength}{0.03\standardunitlength}
	\begin{array}{c}  \hspace{-1.7mm}
        	\raisebox{-8pt}{\input draws/mtower1.tex }
        	\hspace{-1.9mm}
	\end{array}
 $$
\caption{An example of an external circle $C_1$ and $3$ chords that 
touch it
and form a $3$-tower as far as this circle is 
concerned.}\lbl{mtower1}
\end{figure}

This is all good, except that we used the $4$-term relation twice, 
and
we ignored the subleading terms in the $4$-term relation twice.  
We need to deal with the subleading terms, too. Let's consider the 
subleading terms in the second case. The subleading terms of the 
$4$-term relation 
contain  an extra $1$-pair  blink in the complement of the rest of 
the link.
However, the Seifert surface of this $1$-pair blink may 
intersect the discs (i.e., whose 
boundaries
are the chords of the chord diagram) of at most two chords. 
We call such chords {\em marked}.
 See figure \ref{3arcs3}. Therefore
each time we apply the $4$-term relation, the subleading terms have 
an extra $1$-pair blink 
and
we mark at most two chords. Note that without loss of generality, 
the $4$-term
relation is applied to {\em disjoint} balls (embedded in $S^3$), and 
thus the
$1$-pair blinks  bound surfaces disjoint from each other, and from 
the
rest of the components of the link. 
 Therefore, if in case 2 we apply the $4$-term
relation more than $m$ times, the resulting subleading terms lie 
already in $\Fbl m$. If on the other hand, we apply the $4$-term 
relation
at most $m$ times in order to create a $3m$ tower, 
this means that we mark at most $2m$ chords (of the
preferred ones) and therefore have  a subtower of $3m-2m=m$ 
chords.
Similarly, we can deal with the subleading terms of the $4$-term
relation in the beginning of the proof.
This concludes the proof of Step 3.
\end{pf}

The above argument shows that
ignoring  the subleading terms in the $4$-term relation
does not affect the validity of our arguments. In the rest of the 
proof of
proposition \ref{prop.moretubing} we will ignore
such subleading terms. Due to step 3, let us assume that
every external circle of the chord diagram $CD_{L'}$ touches at most
$g_0(m)$ many chords. If $CD_{L'}$ has $h(m) >> g_0(m)$ chords, 
since
every chord touches at most two circles, it implies that the number
$n$ of external circles  satisfies $n >>m$. 

\begin{itemize}
\item{{\bf Step 4}} 
\hspace{0.5cm} In this case, we have $[S^3,L,f] \in \Fbl m$.
\end{itemize}

\begin{pf}
Fix a circle, and choose and chord $c_1$ that lies on the chosen 
circle. 
The chord  touches at most two circles, and these circles
have at most $2 g_0(m)$ many 
 other chords that touch them. Mark all 
the 
circles that these chords touch and
tube the chord $c_1$, see figure \ref{mtower2}. Now consider the
rest of the circles (remembering that the number of circles is much
greater than $m$), and proceed as above.
\end{pf}

This concludes the proof of proposition \ref{prop.moretubing}.
Note that all the functions mentioned taking values in $\BN$ are
constructible. In fact, we leave it as an exercise to the
reader to show that we can take $g_2(m)= c_2 m^3, g_1(m)=c_1 m^4, g_0(m)=
c_0 m^{12}$ and  $h(m)=c m^{13}$, for some constants 
$c_0, c_1, c_2, c$.
\end{pf}

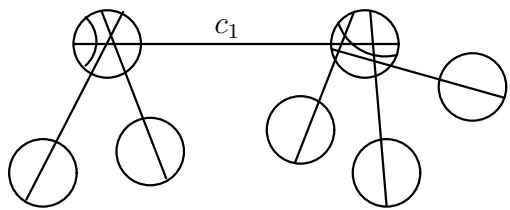
\begin{figure}[htpb]
$$\printname{mtower2}
	\setlength{\unitlength}{0.03\standardunitlength}
	\begin{array}{c}  \hspace{-1.7mm}
        	\raisebox{-8pt}{\begingroup\makeatletter\ifx\SetFigFont\undefined
% extract first six characters in \fmtname
\def\x#1#2#3#4#5#6#7\relax{\def\x{#1#2#3#4#5#6}}%
\expandafter\x\fmtname xxxxxx\relax \def\y{splain}%
\ifx\x\y   % LaTeX or SliTeX?
\gdef\SetFigFont#1#2#3{%
  \ifnum #1<17\tiny\else \ifnum #1<20\small\else
  \ifnum #1<24\normalsize\else \ifnum #1<29\large\else
  \ifnum #1<34\Large\else \ifnum #1<41\LARGE\else
     \huge\fi\fi\fi\fi\fi\fi
  \csname #3\endcsname}%
\else
\gdef\SetFigFont#1#2#3{\begingroup
  \count@#1\relax \ifnum 25<\count@\count@25\fi
  \def\x{\endgroup\@setsize\SetFigFont{#2pt}}%
  \expandafter\x
    \csname \romannumeral\the\count@ pt\expandafter\endcsname
    \csname @\romannumeral\the\count@ pt\endcsname
  \csname #3\endcsname}%
\fi
\fi\endgroup
\begin{picture}(6964,2777)(0,-10)
\thicklines
\put(782.000,2318.500){\arc{903.120}{5.4390}{7.1273}}
\put(5259.016,2795.113){\arc{1372.544}{1.3160}{2.8243}}
\put(1382,2281){\ellipse{948}{948}}
\put(4982,2281){\ellipse{948}{948}}
\put(482,481){\ellipse{948}{948}}
\put(1982,781){\ellipse{948}{948}}
\put(4082,1081){\ellipse{948}{948}}
\put(5282,481){\ellipse{948}{948}}
\put(6482,1681){\ellipse{948}{948}}
\path(932,2281)(5432,2281)
\path(1607,2731)(257,106)
\path(1307,2731)(2207,406)
\path(4832,2731)(4007,631)
\path(5057,2731)(5282,31)
\path(4532,2206)(6932,1531)
\put(2882,2431){\makebox(0,0)[lb]{$c_1$}}
\end{picture} }
        	\hspace{-1.9mm}
	\end{array}
 $$
\caption{An example of a chord $c_1$ to be tubed and of the external 
circles 
that it marks.}\lbl{mtower2}
\end{figure}

Combining propositions \ref{prop.1}, \ref{prop.tubing} and 
\ref{prop.moretubing}
enables us to give a proof of theorem \ref{thm.as2b} as follows:

\begin{pf}[of theorem \ref{thm.as2b}]

Let $h$ be the function as in proposition \ref{prop.moretubing}. We 
want to
show that for every non-negative integer $m$, we have that $\Fas {4 
h(m)}
\subseteq \Fbl m$. We know that $\Fas {4 h(m)}$ is spanned by all 
elements
of the form $[M,L,f]$ where $(L,f)$ is an \ASA\ $4h(m)$-
component 
link in an
\ihs\ $M$. 
For the convenience of the reader, we give  the proof in three 
steps:

\begin{itemize}
\item{{\bf Step 1}} 
\hspace{0.5cm} We may assume that $M=S^3$.
\end{itemize}

Indeed, it follows from the facts that:
 (i) every \ihs\ can be obtained by surgery
on an \ASA\ link $L''$ in $S^3$, (ii) the fundamental equation 
\eqref{eq.fundamental} and (iii) upward induction
on the number of components of $L''$. See also 
\cite[step 1, theorem 1]{GL1}.

>From now on, we assume that $(L,f)$ is an \ASA\ in $S^3$. 

\begin{itemize}
\item{{\bf Step 2}} 
\hspace{0.5cm} We may assume that $L$ contains a trivial sublink 
$L_{tr}(h(m))$
of $h(m)$
components.
\end{itemize}

Indeed, this is nothing but proposition \ref{prop.1}.

\begin{itemize}
\item{{\bf Step 3}} 
\hspace{0.5cm} We may assume that the trivial sublink 
$L_{tr}(h(m))$ bounds
a disjoint union of discs satisfying the properties of proposition 
\ref{prop.moretubing}.
\end{itemize}

Indeed, since $L_{tr}(h(m))$ is a trivial link it bounds a disjoint 
union of discs
$U_i D_i$. Let $L'$ be the sublink of $L$ that the union of the discs 
$D_i$
intersect. Since $L$ is an \ASA\ link, the intersections of the 
components of $L'$ with
each disc $D_i$ come in pairs with opposite orientations. 
Furthermore, using the equation
shown in graphical notation in figure \ref{3bands} (for each of the 
discs $D_i$), we 
may assume that each disc intersects $L'$ in two or four points. 
Furthermore, in case 
a disc $D_i$ intersects the same component of $L'$ in four points, 
then figure \ref{4points}
shows that we can replace such intersections by a linear 
combination of discs that
intersect that component in two points only. This finishes the proof 
of step $3$, and together
with proposition \ref{prop.moretubing} implies the proof of theorem 
\ref{thm.as2b}.
\end{pf}

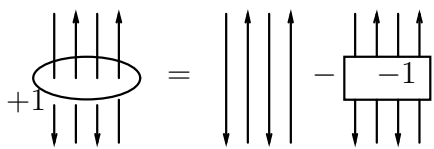
\begin{figure}[htpb]
$$\printname{4points}
	\setlength{\unitlength}{0.03\standardunitlength}
	\begin{array}{c}  \hspace{-1.7mm}
        	\raisebox{-8pt}{\begingroup\makeatletter\ifx\SetFigFont\undefined
% extract first six characters in \fmtname
\def\x#1#2#3#4#5#6#7\relax{\def\x{#1#2#3#4#5#6}}%
\expandafter\x\fmtname xxxxxx\relax \def\y{splain}%
\ifx\x\y   % LaTeX or SliTeX?
\gdef\SetFigFont#1#2#3{%
  \ifnum #1<17\tiny\else \ifnum #1<20\small\else
  \ifnum #1<24\normalsize\else \ifnum #1<29\large\else
  \ifnum #1<34\Large\else \ifnum #1<41\LARGE\else
     \huge\fi\fi\fi\fi\fi\fi
  \csname #3\endcsname}%
\else
\gdef\SetFigFont#1#2#3{\begingroup
  \count@#1\relax \ifnum 25<\count@\count@25\fi
  \def\x{\endgroup\@setsize\SetFigFont{#2pt}}%
  \expandafter\x
    \csname \romannumeral\the\count@ pt\expandafter\endcsname
    \csname @\romannumeral\the\count@ pt\endcsname
  \csname #3\endcsname}%
\fi
\fi\endgroup
\begin{picture}(5937,1839)(0,-10)
\thicklines
\put(1125,912){\ellipse{1500}{600}}
\path(675,1812)(675,912)
\path(1275,1812)(1275,912)
\path(975,537)(975,12)
\path(1575,612)(1575,12)
\path(4725,1212)(4725,612)(5925,612)
	(5925,1212)(4725,1212)
\path(4875,1812)(4875,1212)
\path(5475,1812)(5475,1212)
\path(5775,612)(5775,12)
\path(5175,612)(5175,12)
\path(675,537)(675,12)
\path(645.000,132.000)(675.000,12.000)(705.000,132.000)
\path(1275,537)(1275,12)
\path(1245.000,132.000)(1275.000,12.000)(1305.000,132.000)
\path(3075,1812)(3075,12)
\path(3045.000,132.000)(3075.000,12.000)(3105.000,132.000)
\path(3675,1812)(3675,12)
\path(3645.000,132.000)(3675.000,12.000)(3705.000,132.000)
\path(4875,612)(4875,12)
\path(4845.000,132.000)(4875.000,12.000)(4905.000,132.000)
\path(5475,612)(5475,12)
\path(5445.000,132.000)(5475.000,12.000)(5505.000,132.000)
\path(1005.000,1692.000)(975.000,1812.000)(945.000,1692.000)
\path(975,1812)(975,912)
\path(1605.000,1692.000)(1575.000,1812.000)(1545.000,1692.000)
\path(1575,1812)(1575,912)
\path(3405.000,1692.000)(3375.000,1812.000)(3345.000,1692.000)
\path(3375,1812)(3375,12)
\path(4005.000,1692.000)(3975.000,1812.000)(3945.000,1692.000)
\path(3975,1812)(3975,12)
\path(5205.000,1692.000)(5175.000,1812.000)(5145.000,1692.000)
\path(5175,1812)(5175,1212)
\path(5805.000,1692.000)(5775.000,1812.000)(5745.000,1692.000)
\path(5775,1812)(5775,1212)
\put(2250,837){\makebox(0,0)[lb]{$=$}}
\put(4275,837){\makebox(0,0)[lb]{$-$}}
\put(5175,837){\makebox(0,0)[lb]{$-1$}}
\put(0,462){\makebox(0,0)[lb]{$+1$}}
\end{picture} }
        	\hspace{-1.9mm}
	\end{array}
 $$
\caption{An equality in $\cal M$. The $-1$ in the box indicates a full 
twist.
Notice that all arcs lie in the same link component, and that the two 
links
shown on the right hand side of the figure are homotopic, and 
therefore, 
by doing a number of crossing changes we can rewrite the right hand 
side as
a finite sum of terms each of which contains at least one disc that 
intersects
the link component in two points.}\lbl{4points}
\end{figure}

\subsection{Proof of theorem \ref{thm.as2bl}}
\lbl{sub.2}

In this section we prove theorem \ref{thm.as2bl}.

\begin{pf}[of theorem \ref{thm.as2bl}]
We will first show that for every $m$ we have that:
\begin{equation}
\lbl{eq.show1}
\Fas {3m} \subseteq \Fbl {2m} + \Fas {3m+1}
\end{equation}
 This is equivalent
to showing that $\Fbl {2m}$ spans  the graded space $\Gas {3m}$.
However, we know a set of generators for this graded space, namely
$[L(\Gamma)]$, where $\Gamma$ is a trivalent vertex-oriented 
graph with
$3m$ edges, and $L(\Gamma)$ is the associated $3m$-component 
\as\ link,
see \cite{Oh}, and \cite{GL1}. We will use a slightly different set
of generators of this graded space namely $[L(\Gamma_w)]$, where
\begin{equation}
[L(\Gamma_w)]= \sum_{v} [L(\Gamma(v))]
\end{equation}
where the sum is over all subsets of vertices of $\Gamma$, and 
$\Gamma(v)$
is the result of breaking the vertices of $\Gamma$ in $v$ according
to figure \ref{whitevertex}. Note that  such vertices were called 
{\em white} in \cite{GO2}.

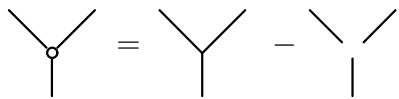
\begin{figure}[htpb]
$$\printname{whitevertex}
	\setlength{\unitlength}{0.03\standardunitlength}
	\begin{array}{c}  \hspace{-1.7mm}
        	\raisebox{-8pt}{\begingroup\makeatletter\ifx\SetFigFont\undefined
% extract first six characters in \fmtname
\def\x#1#2#3#4#5#6#7\relax{\def\x{#1#2#3#4#5#6}}%
\expandafter\x\fmtname xxxxxx\relax \def\y{splain}%
\ifx\x\y   % LaTeX or SliTeX?
\gdef\SetFigFont#1#2#3{%
  \ifnum #1<17\tiny\else \ifnum #1<20\small\else
  \ifnum #1<24\normalsize\else \ifnum #1<29\large\else
  \ifnum #1<34\Large\else \ifnum #1<41\LARGE\else
     \huge\fi\fi\fi\fi\fi\fi
  \csname #3\endcsname}%
\else
\gdef\SetFigFont#1#2#3{\begingroup
  \count@#1\relax \ifnum 25<\count@\count@25\fi
  \def\x{\endgroup\@setsize\SetFigFont{#2pt}}%
  \expandafter\x
    \csname \romannumeral\the\count@ pt\expandafter\endcsname
    \csname @\romannumeral\the\count@ pt\endcsname
  \csname #3\endcsname}%
\fi
\fi\endgroup
\begin{picture}(5424,1239)(0,-10)
\thicklines
\put(612,612){\ellipse{150}{150}}
\path(2112,1212)(2712,612)
\path(3312,1212)(2712,612)
\path(2712,612)(2712,12)
\path(4212,1212)(4662,762)
\path(5412,1212)(4962,762)
\path(4812,537)(4812,12)
\path(12,1212)(537,687)
\path(1212,1212)(687,687)
\path(612,537)(612,12)
\put(1512,612){\makebox(0,0)[lb]{$=$}}
\put(3687,612){\makebox(0,0)[lb]{$-$}}
\end{picture} }
        	\hspace{-1.9mm}
	\end{array}
 $$
\caption{The definition of a white vertex. Note that each of the 
graphs
represent unit-framed algebraically split links in 
$S^3$.}\lbl{whitevertex}
\end{figure}

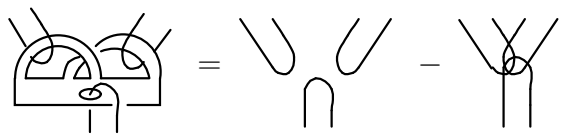
\begin{figure}[htpb]
$$\printname{blink2white}
	\setlength{\unitlength}{0.03\standardunitlength}
	\begin{array}{c}  \hspace{-1.7mm}
        	\raisebox{-8pt}{\begingroup\makeatletter\ifx\SetFigFont\undefined
% extract first six characters in \fmtname
\def\x#1#2#3#4#5#6#7\relax{\def\x{#1#2#3#4#5#6}}%
\expandafter\x\fmtname xxxxxx\relax \def\y{splain}%
\ifx\x\y   % LaTeX or SliTeX?
\gdef\SetFigFont#1#2#3{%
  \ifnum #1<17\tiny\else \ifnum #1<20\small\else
  \ifnum #1<24\normalsize\else \ifnum #1<29\large\else
  \ifnum #1<34\Large\else \ifnum #1<41\LARGE\else
     \huge\fi\fi\fi\fi\fi\fi
  \csname #3\endcsname}%
\else
\gdef\SetFigFont#1#2#3{\begingroup
  \count@#1\relax \ifnum 25<\count@\count@25\fi
  \def\x{\endgroup\@setsize\SetFigFont{#2pt}}%
  \expandafter\x
    \csname \romannumeral\the\count@ pt\expandafter\endcsname
    \csname @\romannumeral\the\count@ pt\endcsname
  \csname #3\endcsname}%
\fi
\fi\endgroup
\begin{picture}(7674,1760)(0,-10)
\thicklines
\put(688.000,761.000){\arc{1200.000}{3.1416}{6.2832}}
\put(688.000,761.000){\arc{900.000}{3.1416}{6.2832}}
\put(1475.500,678.500){\arc{988.863}{4.3234}{6.1155}}
\put(1594.818,849.636){\arc{1051.416}{3.8995}{6.4526}}
\put(1137,537){\ellipse{300}{150}}
\path(238,761)(763,761)
\path(913,761)(1138,761)
\path(1288,761)(1963,761)
\path(87,762)(87,387)(1437,387)
\path(1662,387)(2112,387)(2112,762)
\path(1137,312)(1137,12)
\path(12,1587)(237,1212)
\path(2262,1437)(2037,1137)
\path(3687,912)(3237,1587)
\path(4887,912)(5337,1587)
\path(4137,612)(4137,87)
\path(6745,900)(6295,1575)
\path(7212,912)(7662,1587)
\path(6912,912)(6912,87)
\path(7287,387)(7287,87)
\path(988,1061)	(929.600,1039.966)
	(887.477,1021.625)
	(838.000,986.000)

\path(838,986)	(816.512,954.189)
	(798.021,909.601)
	(780.770,846.963)
	(772.060,807.226)
	(763.000,761.000)

\path(1063,986)	(1012.604,938.776)
	(988.000,911.000)

\path(988,911)	(958.634,860.210)
	(938.729,818.415)
	(913.000,761.000)

\path(1137,537)	(1176.511,604.788)
	(1212.337,650.040)
	(1287.000,687.000)

\path(1287,687)	(1356.916,677.344)
	(1421.055,641.707)
	(1474.417,591.217)
	(1512.000,537.000)

\path(1512,537)	(1529.491,467.053)
	(1530.220,426.053)
	(1527.548,383.576)
	(1522.932,341.527)
	(1517.830,301.812)
	(1512.000,237.000)

\path(1512,237)	(1512.000,203.952)
	(1512.000,158.951)
	(1512.000,96.725)
	(1512.000,57.505)
	(1512.000,12.000)

\path(312,1737)	(343.100,1692.057)
	(369.792,1653.237)
	(411.713,1591.331)
	(441.277,1546.010)
	(462.000,1512.000)

\path(462,1512)	(500.791,1449.800)
	(524.703,1411.687)
	(549.045,1370.932)
	(571.881,1329.062)
	(591.277,1287.599)
	(605.295,1248.070)
	(612.000,1212.000)

\path(612,1212)	(611.166,1156.395)
	(600.127,1092.821)
	(576.276,1032.586)
	(537.000,987.000)

\path(537,987)	(501.795,972.551)
	(460.882,970.286)
	(387.000,987.000)

\path(387,987)	(355.729,1009.564)
	(312.000,1062.000)

\path(1887,1662)	(1855.900,1617.057)
	(1829.208,1578.237)
	(1787.288,1516.331)
	(1757.723,1471.010)
	(1737.000,1437.000)

\path(1737,1437)	(1698.209,1374.800)
	(1674.297,1336.687)
	(1649.955,1295.932)
	(1627.119,1254.062)
	(1607.723,1212.599)
	(1593.705,1173.070)
	(1587.000,1137.000)

\path(1587,1137)	(1587.834,1081.395)
	(1598.872,1017.821)
	(1622.724,957.586)
	(1662.000,912.000)

\path(1662,912)	(1697.205,897.551)
	(1738.118,895.286)
	(1812.000,912.000)

\path(1812,912)	(1843.271,934.564)
	(1887.000,987.000)

\path(3687,1587)	(3718.100,1542.057)
	(3744.792,1503.237)
	(3786.712,1441.331)
	(3816.277,1396.010)
	(3837.000,1362.000)

\path(3837,1362)	(3875.791,1299.800)
	(3899.703,1261.687)
	(3924.045,1220.932)
	(3946.881,1179.062)
	(3966.277,1137.599)
	(3980.295,1098.070)
	(3987.000,1062.000)

\path(3987,1062)	(3986.170,1006.395)
	(3975.131,942.821)
	(3951.277,882.586)
	(3912.000,837.000)

\path(3912,837)	(3876.795,822.551)
	(3835.882,820.286)
	(3762.000,837.000)

\path(3762,837)	(3730.729,859.564)
	(3687.000,912.000)

\path(4887,1587)	(4855.900,1542.057)
	(4829.208,1503.237)
	(4787.288,1441.331)
	(4757.723,1396.010)
	(4737.000,1362.000)

\path(4737,1362)	(4698.209,1299.800)
	(4674.297,1261.687)
	(4649.955,1220.932)
	(4627.119,1179.062)
	(4607.723,1137.599)
	(4593.705,1098.070)
	(4587.000,1062.000)

\path(4587,1062)	(4587.830,1006.395)
	(4598.869,942.821)
	(4622.723,882.586)
	(4662.000,837.000)

\path(4662,837)	(4697.205,822.551)
	(4738.118,820.286)
	(4812.000,837.000)

\path(4812,837)	(4843.271,859.564)
	(4887.000,912.000)

\path(4137,612)	(4176.511,679.788)
	(4212.337,725.040)
	(4287.000,762.000)

\path(4287,762)	(4356.916,752.344)
	(4421.055,716.707)
	(4474.417,666.217)
	(4512.000,612.000)

\path(4512,612)	(4529.491,542.053)
	(4530.220,501.053)
	(4527.547,458.576)
	(4522.932,416.527)
	(4517.830,376.812)
	(4512.000,312.000)

\path(4512,312)	(4512.000,278.952)
	(4512.000,233.951)
	(4512.000,171.725)
	(4512.000,132.505)
	(4512.000,87.000)

\path(6762,1587)	(6793.100,1542.057)
	(6819.792,1503.237)
	(6861.712,1441.331)
	(6891.277,1396.010)
	(6912.000,1362.000)

\path(6912,1362)	(6950.791,1299.800)
	(6974.703,1261.687)
	(6999.045,1220.932)
	(7021.881,1179.062)
	(7041.277,1137.599)
	(7055.295,1098.070)
	(7062.000,1062.000)

\path(7062,1062)	(7061.170,1006.395)
	(7050.131,942.821)
	(7026.277,882.586)
	(6987.000,837.000)

\path(6987,837)	(6951.795,822.551)
	(6910.882,820.286)
	(6837.000,837.000)

\path(6837,837)	(6805.729,859.564)
	(6762.000,912.000)

\path(7212,1587)	(7180.900,1542.057)
	(7154.208,1503.237)
	(7112.287,1441.331)
	(7082.723,1396.010)
	(7062.000,1362.000)

\path(7062,1362)	(7023.209,1299.800)
	(6999.297,1261.687)
	(6974.955,1220.932)
	(6952.119,1179.062)
	(6932.723,1137.599)
	(6918.705,1098.070)
	(6912.000,1062.000)

\path(6912,1062)	(6912.830,1006.395)
	(6923.869,942.821)
	(6947.723,882.586)
	(6987.000,837.000)

\path(6987,837)	(7022.205,822.551)
	(7063.118,820.286)
	(7137.000,837.000)

\path(7137,837)	(7168.271,859.564)
	(7212.000,912.000)

\path(6912,912)	(6951.511,979.788)
	(6987.337,1025.040)
	(7062.000,1062.000)

\path(7062,1062)	(7131.916,1052.344)
	(7196.055,1016.707)
	(7249.417,966.217)
	(7287.000,912.000)

\path(7287,912)	(7304.491,842.053)
	(7305.220,801.053)
	(7302.547,758.576)
	(7297.932,716.527)
	(7292.830,676.812)
	(7287.000,612.000)

\path(7287,612)	(7287.000,578.952)
	(7287.000,533.951)
	(7287.000,471.725)
	(7287.000,432.505)
	(7287.000,387.000)

\put(5712,837){\makebox(0,0)[lb]{$-$}}
\put(2637,837){\makebox(0,0)[lb]{$=$}}
\end{picture} }
        	\hspace{-1.9mm}
	\end{array}
 $$
\caption{This figure represents a special case of equation 
\eqref{eq.fundamental} in graphical notation. The present identity 
holds
in $\cal G^{as}_{\ast } \cal M$
On the left shown is a $1$-pair blink, which (after surgery)
 corresponds in $\cal M$
to the difference of two terms. The first term is shown on the first 
part on
 right, and the second term is (surgically equivalent to) the result of
blowing down the $1$-pair blink. Notice that each of the two 
components of
the blink are an unknot and can be blown down in any order.  
}\lbl{blink2white}
\end{figure}

Using the identity in figure \ref{blink2white}
we see that summing over each  white vertex (in the sum of 
$[L(\Gamma_w)]$)
is equivalent to summing over a $1$-pair blink. Since the graphs 
$\Gamma$
are trivalent with $3m$ edges, (and therefore of $2m$ vertices), we 
have
that $[L(\Gamma_w)] \in \Fbl {2m}$.
This finishes the proof of equation \ref{eq.show1}.
Now interpolating equation \eqref{eq.show1} implies that:
$\Fas {3m} \subseteq \widehat{\Fbl {2m}}$, which finishes the first
part of the theorem.

The second part follows immediately from the first, using theorems 
\ref{thm.bl2as} and \ref{thm.as2b}.
\end{pf}

\subsection{Proof of theorem \ref{thm.bl2as}}
\lbl{sub.3}

This section is devoted to the proof of theorem \ref{thm.bl2as}.

 The proof presented here is similar to the proof of theorem $1$ of
\cite{GL2}. It uses primary and secondary induction as well as 
the identities (for not-necessarily \as\ links) of section 
\ref{sub.review}. For the convenience of the reader, we separate
 the proof into $5$ steps.
We begin with some definitions that will be useful. A triple of links 
$T=(L,L_b ,L_{bl})$ in an \ihs\ $M$ consists of an algebraically split 
link $L$, a boundary link $L_b$ and a blink $L_{bl}$ such that each 
component of $L_b$ and pair of $L_{bl}$ bounds a connected 
oriented Seifert surface in $M$ and these surfaces are disjoint from 
each other and from $L$. Such a (disconnected) surface is called
an  {\em admissible} Seifert 
surface for $T$. If $k=|L|, n=|L_b |, m=|L_{bl}|$, then we call $T$ a 
$(k,n,m)-clink$. The {\em genus} $g(T)$ is the minimal total genus of 
an admissible Seifert surface of $T$. An admissible framing for $T$ 
is 
one which is unit on $L\cup L_b$ and unit Seifert-framing
 on $L_{bl}$. We can 
then define $[M,T,f] \in \M$ in the usual way. We will prove that 
\begin{equation}
\lbl{eq.nnn}
[M,T,f]\in\Fas{3n+3m/2}
\end{equation}
 Note that in \cite[Theorem 1]{GL2} we proved this fact in 
the special case  $m=0$. Note also that the present theorem 
\ref{thm.bl2as}
is the case $n=0$ of equation \eqref{eq.nnn}.
The argument for equation \eqref{eq.nnn} is a generalization of 
that in \cite{GL2}, proceeding by primary downward induction on 
$k(T)=k$ and secondary upward induction on $g(T)$.

\begin{itemize}
\item{{\bf Step 1}} 
\hspace{0.5cm}We may assume that $M=S^3$.
\end{itemize}

The proof follows from the following $3$ facts:
\begin{itemize}
\item
     Every \ihs\ $M$ can be  converted, by surgery on an \ASA\ link 
$L'$,
into $S^3$.
\item
We may assume that $L'$
above can be chosen so that the Seifert 
surfaces bounded by $L_b\cup L_{bl}$ are
disjoint from $L'$.  This follows by general position, since the 
Seifert surfaces are contained in a regular neighborhood of some 
embedded graph in $M$ and we can perturb $L'$ away from this graph. 
\item
     Equation \eqref{eq.fundamental}
\end{itemize}
and upward induction on the number of components of $L'$, see also
step $1$ in theorem $1$ of \cite{GL2}.

Suppose now that $(L,L_b ,L_{bl})$ is a $(k,n,m)$-clink in $S^3$.
If $k\geq 3n+3m/2$ we are done by definition.
If $g(T)=0$ (i.e., $L_{bl}$ is an unblink, see figure \ref{f.unblink})
we are also done, since in this case we have that 
$[S^3,T,f]=0$. Indeed, if $L_{bl}$ is a unit-Seifert framed
 unblink in an \ihs\ $M$, then
$[M,L_{bl},f]=M - M_{L_{bl},f}$, and by applying Kirby moves to a band
of a genus $0$ surface that $L_{bl}$ bounds, we deduce that 
$M_{L_{bl},f}$
is diffeomorphic to $M$, and thus $[M,L_{bl},f]=0$. 
This begins the induction.

\begin{itemize}
\item{{\bf Step 2}}
\hspace{0.5cm}We may assume that every component of $L$ 
is
unknotted. 
\end{itemize}

Equation \eqref{eq.fundamental} 
implies that the change of $[S^3, T,f]$ before
and after a crossing change in the same component of $L$
can be written as $[S^3, T', f \cup \pm 1]$ where $T'=(L\cup C,L_b 
,L_{bl})$ and $C$
is a circle that encloses the crossing to be changed.
Since $k(T')>k(T)$, by using the primary
inductive hypothesis we can change crossings of components of $L 
$, and thus assume that each component of $L $
is unknotted.

\begin{itemize}
\item{{\bf Step 3}}
\hspace{0.5cm}Suppose that $L_b\cup L_{bl} =\partial\Sigma$, 
where 
$\Sigma$ is an admissible
 Seifert surface for $T$. 
We may assume that $\Sigma$ is embedded in a standard, 
almost planar (except for the necessary band crossings) way.  
See Figure \ref{genus1}.
\end{itemize}

 This follows using the same argument as in 
\cite{GL2} by introducing extra components into $L$ in order 
to change band crossings.

Let $ \{K_i \} $ (for $1 \leq i \leq k(T) $)
 denote the components of $L $. Since
by step $2$  
they are unknotted, we may choose embedded disks $D_i$ so 
that $K_i =\partial D_i$. 
Furthermore, since $\Sigma$ is just a thickening of a wedge of 
circles,
 we may choose  the 
$D_i$ so that their intersections with $\Sigma$ consist 
of a number of transverse penetrations of the interiors of 
the $D_i$ by the bands of $\Sigma$. See Figure \ref{f4}. We 
will be interested in counting the number of ''band 
penetrations''.

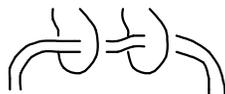
\begin{figure}[htpb]
$$\printname{f4}
	\setlength{\unitlength}{0.03\standardunitlength}
	\begin{array}{c}  \hspace{-1.7mm}
        	\raisebox{-8pt}{\begingroup\makeatletter\ifx\SetFigFont\undefined
% extract first six characters in \fmtname
\def\x#1#2#3#4#5#6#7\relax{\def\x{#1#2#3#4#5#6}}%
\expandafter\x\fmtname xxxxxx\relax \def\y{splain}%
\ifx\x\y   % LaTeX or SliTeX?
\gdef\SetFigFont#1#2#3{%
  \ifnum #1<17\tiny\else \ifnum #1<20\small\else
  \ifnum #1<24\normalsize\else \ifnum #1<29\large\else
  \ifnum #1<34\Large\else \ifnum #1<41\LARGE\else
     \huge\fi\fi\fi\fi\fi\fi
  \csname #3\endcsname}%
\else
\gdef\SetFigFont#1#2#3{\begingroup
  \count@#1\relax \ifnum 25<\count@\count@25\fi
  \def\x{\endgroup\@setsize\SetFigFont{#2pt}}%
  \expandafter\x
    \csname \romannumeral\the\count@ pt\expandafter\endcsname
    \csname @\romannumeral\the\count@ pt\endcsname
  \csname #3\endcsname}%
\fi
\fi\endgroup
\begin{picture}(3021,1231)(0,-10)
\thicklines
\path(987,1208)	(1037.646,1152.924)
	(1080.540,1104.628)
	(1116.339,1062.232)
	(1145.704,1024.858)
	(1187.766,961.659)
	(1212.000,908.000)

\path(1212,908)	(1222.304,869.566)
	(1230.250,824.116)
	(1235.546,773.953)
	(1237.898,721.377)
	(1237.012,668.693)
	(1232.596,618.200)
	(1224.356,572.202)
	(1212.000,533.000)

\path(1212,533)	(1176.488,467.584)
	(1123.181,399.267)
	(1058.034,341.567)
	(987.000,308.000)

\path(987,308)	(928.400,305.852)
	(865.646,321.189)
	(807.319,348.681)
	(762.000,383.000)

\path(762,383)	(737.964,414.653)
	(718.624,459.189)
	(702.222,521.880)
	(694.573,561.682)
	(687.000,608.000)

\path(87,533)	(139.437,592.485)
	(179.573,634.970)
	(237.000,683.000)

\path(237,683)	(285.066,706.222)
	(347.228,728.495)
	(410.526,746.771)
	(462.000,758.000)

\path(462,758)	(511.957,761.573)
	(574.841,761.176)
	(637.555,759.191)
	(687.000,758.000)

\path(687,758)	(736.573,758.000)
	(799.500,758.000)
	(862.427,758.000)
	(912.000,758.000)

\path(912,758)	(938.014,758.000)
	(987.000,758.000)

\path(612,1208)	(608.473,1152.036)
	(607.298,1110.725)
	(612.000,1058.000)

\path(612,1058)	(626.344,1023.431)
	(649.496,983.000)
	(672.650,942.569)
	(687.000,908.000)

\path(687,908)	(689.348,881.634)
	(687.000,833.000)

\path(1962,1208)	(2012.646,1152.924)
	(2055.540,1104.628)
	(2091.339,1062.232)
	(2120.704,1024.858)
	(2162.766,961.659)
	(2187.000,908.000)

\path(2187,908)	(2197.304,869.566)
	(2205.250,824.116)
	(2210.546,773.953)
	(2212.898,721.377)
	(2212.012,668.693)
	(2207.596,618.200)
	(2199.356,572.202)
	(2187.000,533.000)

\path(2187,533)	(2151.488,467.584)
	(2098.181,399.267)
	(2033.034,341.567)
	(1962.000,308.000)

\path(1962,308)	(1903.400,305.852)
	(1840.646,321.189)
	(1782.319,348.681)
	(1737.000,383.000)

\path(1737,383)	(1712.964,414.653)
	(1693.624,459.189)
	(1677.222,521.880)
	(1669.573,561.682)
	(1662.000,608.000)

\path(1362,608)	(1418.223,605.621)
	(1459.620,604.827)
	(1512.000,608.000)

\path(1512,608)	(1562.622,622.227)
	(1624.500,645.500)
	(1686.378,668.773)
	(1737.000,683.000)

\path(1737,683)	(1789.380,686.173)
	(1830.777,685.379)
	(1887.000,683.000)

\path(1362,758)	(1418.223,755.621)
	(1459.620,754.827)
	(1512.000,758.000)

\path(1512,758)	(1562.622,772.227)
	(1624.500,795.500)
	(1686.378,818.773)
	(1737.000,833.000)

\path(1737,833)	(1789.380,836.173)
	(1830.777,835.379)
	(1887.000,833.000)

\path(1587,1208)	(1583.473,1152.036)
	(1582.298,1110.725)
	(1587.000,1058.000)

\path(1587,1058)	(1601.344,1023.431)
	(1624.496,983.000)
	(1647.650,942.569)
	(1662.000,908.000)

\path(1662,908)	(1664.348,881.634)
	(1662.000,833.000)

\path(2337,833)	(2382.918,819.593)
	(2422.434,807.777)
	(2484.896,788.041)
	(2529.660,772.035)
	(2562.000,758.000)

\path(2562,758)	(2614.588,729.314)
	(2678.595,689.735)
	(2740.554,646.788)
	(2787.000,608.000)

\path(2787,608)	(2826.045,561.809)
	(2869.421,500.274)
	(2909.086,436.352)
	(2937.000,383.000)

\path(2937,383)	(2959.936,317.991)
	(2971.389,277.530)
	(2982.281,234.676)
	(2992.201,191.613)
	(3000.735,150.523)
	(3012.000,83.000)

\path(3012,83)	(3013.193,56.881)
	(3012.000,8.000)

\path(2337,608)	(2382.757,593.746)
	(2422.158,581.326)
	(2484.529,561.106)
	(2529.385,545.584)
	(2562.000,533.000)

\path(2562,533)	(2596.252,519.881)
	(2639.205,502.798)
	(2712.000,458.000)

\path(2712,458)	(2740.144,410.192)
	(2761.867,348.001)
	(2777.407,284.560)
	(2787.000,233.000)

\path(2787,233)	(2789.384,199.695)
	(2789.119,157.771)
	(2787.000,83.000)

\path(2787,83)	(2787.000,56.986)
	(2787.000,8.000)

\path(87,533)	(55.494,478.802)
	(33.664,438.080)
	(12.000,383.000)

\path(12,383)	(8.473,349.429)
	(8.865,307.535)
	(12.000,233.000)

\path(12,233)	(12.000,180.965)
	(12.000,139.482)
	(12.000,83.000)

\path(162,83)	(157.967,127.793)
	(155.087,166.505)
	(152.783,228.324)
	(155.087,273.731)
	(162.000,308.000)

\path(162,308)	(186.728,364.410)
	(225.667,428.375)
	(270.273,488.402)
	(312.000,533.000)

\path(312,533)	(383.768,577.246)
	(426.218,596.184)
	(462.000,608.000)

\path(462,608)	(512.354,613.290)
	(575.194,612.702)
	(637.687,609.763)
	(687.000,608.000)

\path(687,608)	(731.068,608.000)
	(791.070,608.000)
	(829.243,608.000)
	(874.037,608.000)
	(926.329,608.000)
	(987.000,608.000)

\end{picture} }
        	\hspace{-1.9mm}
	\end{array}
 $$
\caption{A band of a surface penetrating two pieces of 
discs.}\lbl{f4}
\end{figure}

\begin{itemize}
\item{{\bf Step 4}}
\hspace{0.5cm}We may assume that every band of $\Sigma$ 
penetrates 
at least one $D_i$.
\end{itemize}

\begin{pf}
This follows by an argument similar to that in \cite{GL2}. If 
some band penetrates no $D_i$, then we may arrange, as in 
\cite{GL2}, that the circle in $\Sigma$ going through that band 
bounds a disk in $S^3$ disjoint from $L\cup\Sigma$. Then, 
depending on the band, we can either reduce the genus, as in 
\cite{GL2}, or remove one component of that pair. thus we have 
turned one of the blink components into a bounding component and 
eliminated the band with no penetrations. If $T'$ is the new clink, 
then it is easy to see that $[S^3 ,T',f|T']=[S^3 ,T,f]$.
\end{pf}

\begin{itemize}
\item{{\bf Step 5}} 
\hspace{0.5cm}We may assume that each disc  $D_i$ has at most 
two band penetrations.
\end{itemize}

This follows precisely as in \cite{GL2}.

\begin{itemize}
\item{{\bf Step 6}} If any component $\Sigma_j$ of the Seifert 
surface of $L_b$ has genus one and a band of $\Sigma_j$ penetrates 
only one disk $D_i$, then we may assume that $D_i$ is penetrated by 
no other bands of $\Sigma$.
\end{itemize}

Again this will follow by the same argument as in \cite{GL2}.

We can now complete the proof of theorem \ref{thm.bl2as}  by 
counting the band penetrations. Suppose that $T$ is a $(k,n,m)$-
clink satisfying all the assertions of the previous steps. Since $[S^3 
,T,f]\in\Fas{k+n}$, it suffices to show that $k\geq 2n+3m/2$. Let 
$b=$number of penetrations of $\cup D_i$ by bands of $\Sigma$. Set 
$n=n_0 +n_1$, where $n_0$ is the number of Seifert surface 
components for $L_b$ of genus one. Set $2n_0 =n_0^{\p}+n_0^{\p\p}$ 
where $n_0^{\p}$ is the number of bands of these genus one Seifert 
surfaces with only one disk penetration.

Now it follows from Steps 4 and 6 that
$$b\geq 3m+4n_1 +n_0^{\p}+2n_0^{\p\p} $$
If we write $k=n_0^{\p}+k'$ then it follows from Step 5 that $b\leq 
n_0^{\p}+2k'$. Combining these two inequalities gives
$$n_0^{\p}+2k'\geq 3m+4n_1 +n_0^{\p}+2n_0^{\p\p} $$
But then we have 
$$2k=2n_0^{\p}+2k'\geq 3m+4n_1 +2n_0^{\p}+2n_0^{\p\p}=3m+4n$$ 
which was to be proved.
The proof of theorem \ref{thm.bl2as} is complete.
\qed

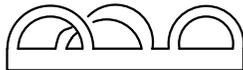
\begin{figure}[htpb]
$$\printname{genus1}
	\setlength{\unitlength}{0.03\standardunitlength}
	\begin{array}{c}  \hspace{-1.7mm}
        	\raisebox{-8pt}{\begingroup\makeatletter\ifx\SetFigFont\undefined
% extract first six characters in \fmtname
\def\x#1#2#3#4#5#6#7\relax{\def\x{#1#2#3#4#5#6}}%
\expandafter\x\fmtname xxxxxx\relax \def\y{splain}%
\ifx\x\y   % LaTeX or SliTeX?
\gdef\SetFigFont#1#2#3{%
  \ifnum #1<17\tiny\else \ifnum #1<20\small\else
  \ifnum #1<24\normalsize\else \ifnum #1<29\large\else
  \ifnum #1<34\Large\else \ifnum #1<41\LARGE\else
     \huge\fi\fi\fi\fi\fi\fi
  \csname #3\endcsname}%
\else
\gdef\SetFigFont#1#2#3{\begingroup
  \count@#1\relax \ifnum 25<\count@\count@25\fi
  \def\x{\endgroup\@setsize\SetFigFont{#2pt}}%
  \expandafter\x
    \csname \romannumeral\the\count@ pt\expandafter\endcsname
    \csname @\romannumeral\the\count@ pt\endcsname
  \csname #3\endcsname}%
\fi
\fi\endgroup
\begin{picture}(3325,951)(0,-10)
\thicklines
\put(612.000,312.000){\arc{1200.000}{3.1416}{6.2832}}
\put(612.000,312.000){\arc{900.000}{3.1416}{6.2832}}
\put(1399.500,229.500){\arc{988.863}{4.3234}{6.1155}}
\put(1518.818,400.636){\arc{1051.416}{3.8995}{6.4526}}
\put(2787.000,382.312){\arc{1059.375}{3.0085}{6.4163}}
\put(2787.000,380.750){\arc{762.500}{2.9603}{6.4645}}
\path(12,312)(12,12)(3312,12)(3312,312)
\path(162,312)(687,312)
\path(837,312)(1062,312)
\path(1212,312)(1887,312)
\path(2037,312)(2262,312)
\path(2412,312)(3162,312)
\path(912,612)	(853.600,590.966)
	(811.477,572.625)
	(762.000,537.000)

\path(762,537)	(740.512,505.189)
	(722.021,460.601)
	(704.770,397.963)
	(696.060,358.226)
	(687.000,312.000)

\path(987,537)	(936.604,489.776)
	(912.000,462.000)

\path(912,462)	(882.634,411.210)
	(862.729,369.415)
	(837.000,312.000)

\end{picture} }
        	\hspace{-1.9mm}
	\end{array}
 $$
\caption{An example of a $1$-pair blink that bounds a genus $1$
surface. Note that the surface has $3$ bands.}\lbl{genus1}
\end{figure}

\begin{remark}
The bound obtained in theorem~\ref{thm.bl2as} is sharp.
Indeed, if $T$ is as in figure \ref{2genus1}
then $[S^3, T, f] \in \Fas 3$, by theorem~\ref{thm.bl2as}, but we 
claim that $[S^3, T, f] \not\in \Fas 4$. Indeed, using figures
\ref{whitevertex} and \ref{blink2white} we see that 
figure \ref{2genus1} represents the element $[\Theta] \in \Fas 3$,
where $[\Theta] \in \M$ is the element represented by the 
 trivalent graph $\Theta$, with white vertices. But this element of 
$\Fas 3$
is nontrivial in $\Gas 3$, see \cite[proposition 2.13]{GO2}, or 
proposition
\cite[theorem 6]{GL1}.
This implies in particular that the analogue of the
key lemma $2.1$ of \cite{GL2} for blinks is false, and that
the last step $6$ of theorem $1$ of \cite{GL2} would be false
for blinks.
\end{remark}

\begin{figure}[htpb]
 $$\printname{2genus1}
	\setlength{\unitlength}{0.03\standardunitlength}
	\begin{array}{c}  \hspace{-1.7mm}
        	\raisebox{-8pt}{\begingroup\makeatletter\ifx\SetFigFont\undefined
% extract first six characters in \fmtname
\def\x#1#2#3#4#5#6#7\relax{\def\x{#1#2#3#4#5#6}}%
\expandafter\x\fmtname xxxxxx\relax \def\y{splain}%
\ifx\x\y   % LaTeX or SliTeX?
\gdef\SetFigFont#1#2#3{%
  \ifnum #1<17\tiny\else \ifnum #1<20\small\else
  \ifnum #1<24\normalsize\else \ifnum #1<29\large\else
  \ifnum #1<34\Large\else \ifnum #1<41\LARGE\else
     \huge\fi\fi\fi\fi\fi\fi
  \csname #3\endcsname}%
\else
\gdef\SetFigFont#1#2#3{\begingroup
  \count@#1\relax \ifnum 25<\count@\count@25\fi
  \def\x{\endgroup\@setsize\SetFigFont{#2pt}}%
  \expandafter\x
    \csname \romannumeral\the\count@ pt\expandafter\endcsname
    \csname @\romannumeral\the\count@ pt\endcsname
  \csname #3\endcsname}%
\fi
\fi\endgroup
\begin{picture}(3325,2441)(0,-10)
\thicklines
\put(612.000,2114.000){\arc{1200.000}{6.2832}{9.4248}}
\put(612.000,2114.000){\arc{900.000}{6.2832}{9.4248}}
\put(1399.500,2196.500){\arc{988.863}{0.1676}{1.9598}}
\put(1518.818,2025.363){\arc{1051.415}{6.1138}{8.6669}}
\put(2787.000,2043.688){\arc{1059.375}{6.1500}{9.5579}}
\put(2787.000,2045.250){\arc{762.500}{6.1019}{9.6061}}
\put(612.000,312.000){\arc{1200.000}{3.1416}{6.2832}}
\put(612.000,312.000){\arc{900.000}{3.1416}{6.2832}}
\put(1399.500,229.500){\arc{988.863}{4.3234}{6.1155}}
\put(1518.818,400.636){\arc{1051.416}{3.8995}{6.4526}}
\put(2787.000,382.312){\arc{1059.375}{3.0085}{6.4163}}
\put(2787.000,380.750){\arc{762.500}{2.9603}{6.4645}}
\path(12,2114)(12,2414)(3312,2414)(3312,2114)
\path(162,2114)(687,2114)
\path(837,2114)(1062,2114)
\path(1212,2114)(1887,2114)
\path(2037,2114)(2262,2114)
\path(2412,2114)(3162,2114)
\path(12,312)(12,12)(3312,12)(3312,312)
\path(162,312)(687,312)
\path(837,312)(1062,312)
\path(1212,312)(1887,312)
\path(2037,312)(2262,312)
\path(2412,312)(3162,312)
\path(912,1814)	(853.600,1835.034)
	(811.477,1853.375)
	(762.000,1889.000)

\path(762,1889)	(740.512,1920.811)
	(722.021,1965.399)
	(704.770,2028.037)
	(696.060,2067.774)
	(687.000,2114.000)

\path(987,1889)	(936.604,1936.224)
	(912.000,1964.000)

\path(912,1964)	(882.634,2014.790)
	(862.729,2056.585)
	(837.000,2114.000)

\path(912,612)	(853.600,590.966)
	(811.477,572.625)
	(762.000,537.000)

\path(762,537)	(740.512,505.189)
	(722.021,460.601)
	(704.770,397.963)
	(696.060,358.226)
	(687.000,312.000)

\path(987,537)	(936.604,489.776)
	(912.000,462.000)

\path(912,462)	(882.634,411.210)
	(862.729,369.415)
	(837.000,312.000)

\path(536,1738)	(576.901,1800.224)
	(611.000,1813.000)

\path(611,1813)	(682.420,1757.290)
	(706.315,1708.400)
	(724.111,1651.634)
	(737.144,1591.600)
	(746.747,1532.910)
	(754.254,1480.173)
	(761.000,1438.000)

\path(761,1438)	(766.419,1400.319)
	(770.097,1354.703)
	(772.130,1303.804)
	(772.614,1250.279)
	(771.646,1196.781)
	(769.324,1145.965)
	(765.743,1100.487)
	(761.000,1063.000)

\path(761,1063)	(754.905,1013.377)
	(748.804,950.717)
	(740.892,880.590)
	(729.365,808.570)
	(712.417,740.229)
	(688.243,681.138)
	(655.039,636.871)
	(611.000,613.000)

\path(611,613)	(577.801,625.446)
	(536.000,688.000)

\path(461,1438)	(442.185,1394.448)
	(426.507,1356.622)
	(403.681,1295.511)
	(390.765,1249.395)
	(386.000,1213.000)

\path(386,1213)	(388.560,1165.208)
	(400.741,1103.965)
	(410.989,1066.102)
	(424.302,1022.239)
	(440.899,971.498)
	(461.000,913.000)

\path(1436,1738)	(1476.901,1800.224)
	(1511.000,1813.000)

\path(1511,1813)	(1582.420,1757.290)
	(1606.315,1708.400)
	(1624.111,1651.634)
	(1637.144,1591.600)
	(1646.747,1532.910)
	(1654.254,1480.173)
	(1661.000,1438.000)

\path(1661,1438)	(1666.419,1400.319)
	(1670.097,1354.703)
	(1672.130,1303.804)
	(1672.614,1250.279)
	(1671.646,1196.781)
	(1669.324,1145.965)
	(1665.743,1100.487)
	(1661.000,1063.000)

\path(1661,1063)	(1654.905,1013.377)
	(1648.804,950.717)
	(1640.892,880.590)
	(1629.365,808.570)
	(1612.417,740.229)
	(1588.243,681.138)
	(1555.039,636.871)
	(1511.000,613.000)

\path(1511,613)	(1477.801,625.446)
	(1436.000,688.000)

\path(2808,1731)	(2848.901,1793.224)
	(2883.000,1806.000)

\path(2883,1806)	(2954.420,1750.290)
	(2978.315,1701.400)
	(2996.111,1644.634)
	(3009.144,1584.600)
	(3018.747,1525.910)
	(3026.254,1473.173)
	(3033.000,1431.000)

\path(3033,1431)	(3038.419,1393.319)
	(3042.097,1347.703)
	(3044.130,1296.804)
	(3044.614,1243.279)
	(3043.646,1189.781)
	(3041.324,1138.965)
	(3037.743,1093.487)
	(3033.000,1056.000)

\path(3033,1056)	(3026.905,1006.377)
	(3020.804,943.717)
	(3012.892,873.590)
	(3001.365,801.570)
	(2984.417,733.229)
	(2960.243,674.138)
	(2927.039,629.871)
	(2883.000,606.000)

\path(2883,606)	(2849.801,618.446)
	(2808.000,681.000)

\path(1361,1513)	(1342.185,1469.448)
	(1326.507,1431.622)
	(1303.681,1370.511)
	(1290.765,1324.395)
	(1286.000,1288.000)

\path(1286,1288)	(1288.560,1240.208)
	(1300.741,1178.965)
	(1310.989,1141.102)
	(1324.302,1097.239)
	(1340.899,1046.498)
	(1361.000,988.000)

\path(2711,1513)	(2692.185,1469.448)
	(2676.507,1431.622)
	(2653.681,1370.511)
	(2640.765,1324.395)
	(2636.000,1288.000)

\path(2636,1288)	(2638.560,1240.208)
	(2650.741,1178.965)
	(2660.989,1141.102)
	(2674.302,1097.239)
	(2690.899,1046.498)
	(2711.000,988.000)

\end{picture} }
        	\hspace{-1.9mm}
	\end{array}
 $$
 \caption{A special case of a $2$-pair blink $L_{bl}$
union a $3$-component \as\ link $L$. The result $[S^3, T, 
f]$
lies in $\Fas 3$ and is non-trivial in $\Gas 3$.}
 \lbl{2genus1}
\end{figure}
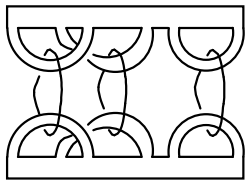

\section{Surgery equivalence and the Seifert matrix}
\lbl{sec.smatrix}

In this section  we prove theorem \ref{thm.seifert}.
This theorem suggests that finite type invariants in the sense of 
blinks, i.e., corresponding to the I-adic filtration of the Torelli group, 
according to Theorem 2, should be expressible in terms of Seifert 
matrix invariants of the associated blinks. Thus Alexander 
polynomial type invariants rather than Jones polynomial type 
invariants should suffice. This is an intriguing consequence which 
we hope to exploit in future work.

It is clear that surgery equivalent $n$-pair blinks (with 
framings which correspond) represent the same element of $\Gbl n$.

Recall the notion of Seifert matrix of an oriented boundary link.  
If $L=(K_1 
,\cdots ,K_n )$ and $K_i =\partial V_i$ where $\{ V_i \}$ are 
disjoint oriented surfaces in a \ihs\ $M$, then the Seifert pairing is 
the collection of bilinear pairings $\sigma_{ij}:H_1 (V_i )\times 
H_1 (V_j )\to\Z$ defined by $\sigma_{ij}(\alpha ,\beta 
)=\lk(\alpha_+ ,\beta )$, the linking number, where $\alpha_+ \in 
H_1(M-\cup_i V_i )$ represents the translate of $\alpha$ off $V_i$ 
in the positive normal direction. We can represent the Seifert 
pairing by a square integral matrix $A$ (the Seifert matrix) divided 
into blocks, each of which represents one of the $\sigma_{ij}$. 
There is an explicit algebraic description of the relation between 
two Seifert matrices of the same link corresponding to different 
choices of $\{ V_i \}$ and different bases of the homology. This 
description uses the notion of S-equivalence and the action of a 
certain group of automorphisms of the free group (see 
\cite{Kstudent} for the 
details). 

The definition of a Seifert pairing of a blink is exactly the same 
using Seifert surfaces of the blink as defined in 
Definition~\ref{def.blink}. The Seifert pairing is again represented 
by a square matrix $A$ separated into blocks representing the 
$\sigma_{ij}$. The relation between two Seifert matrices of the 
same blink will be generally similar to that for boundary links, but 
more complicated, since it is permissible to change the orientation 
any of the Seifert surfaces. We do not want to explore this question 
now and our formulation of Theorem~\ref{thm.seifert} allows us to 
avoid it.

\subsection{Proof of theorem \ref{thm.seifert}}
\lbl{sub.sei1}

We first show that surgery equivalent blinks admit equal Seifert 
matrices. Suppose $(M,L)$ and $(M',L')$ are surgery equivalent. We 
may assume that they are related by a single blink surgery, i.e., if 
$\Sigma$ is a Seifert surface for $L$ in $M$, then there is a $1$-
pair blink $(l,l')$ in $M$ so that $l-l'=\partial\sigma$ where 
$\sigma\subseteq M-\Sigma$ is a Seifert surface for $(l,l')$, and 
that $(M',L')=(M_{(l,l')},L)$ using some unit Seifert-framing of 
$(l,l')$. 
Now we may also regard $\Sigma$ as a Seifert surface for $L'$ in 
$M'$ and so we need to show that if $\alpha ,\beta\in H_1(\Sigma )$, 
then $\lk_M (\alpha_+ ,\beta )=\lk_{M'} (\alpha_+ ,\beta )$ where 
$\lk_M$ denotes the linking number in $M$. In general, given two 
disjoint simple closed curves $\xi ,\eta$ in $M-\Sigma$, we show 
that their linking number in $M$ is the same as in $M'$. Suppose 
$\eta$ bounds an orientable surface $\Lambda$ in $M$. By general 
position $\sigma\cap\Lambda$ is a collection of proper curves in 
$\sigma$ and so, homologically, $\eta =\lk (\xi ,\eta )m+r(m_l 
+m_{l'})$ in $M-\xi -l-l'$, where $m$ is a meridian of $\xi$, $m_l$ 
and  
$m_{l'}$ are meridians of $l$ and $l'$, and $r$ is some integer. Thus 
it 
suffices to observe that $m_l+m_{l'}$ is homologically trivial in 
$M'-\xi$. But, by the definition of a unit Seifert framing 
$m_l^{\p}=m_l+\epsilon l$ and $m_{l'}^{\p}=m_{l'}-\epsilon l'$, for 
some 
$\epsilon =\pm 1$, and so 
$$m_l +m_{l'}=m_l^{\p}+m_{l'}^{\p}+\epsilon (l'-l)
=m_l^{\p}+m_{l'}^{\p}+\epsilon\partial\Sigma =0 \in H_1 (M'-\xi )$$

\subsection{Conclusion of proof}
\lbl{sub.sei2}
Now suppose that two blinks $(M,L)$ and $(M',L')$ admit the same 
Seifert matrix. We may first of all assume that $M=M'=S^3$ since, by 
Proposition~\ref{thm.boundary}, we can convert any $3$-manifold 
into $S^3$ by surgery on a boundary link, which we can assume is far 
away from any other given link. By the observation in 
Remark~\ref{rem.bbl} this is the same as surgery on some blink. Let 
$\S ,\S'$ be Seifert surfaces for $L,L'$ which give identical Seifert 
matrix $A$. Since $A-A^T$ is the intersection matrix of $\S$ and 
$\S'$, we conclude that $\S$ and $\S'$ are diffeomorphic. 

Now it is an easy consequence of Smale theory that the regular 
homotopy type of an embedding of a bounded surface in $S^3$ is 
determined by the twisting numbers of the bands mod 2. Since these 
twisting numbers are determined by the Seifert matrix, it follows 
that $\S$ and $\S'$ are regularly homotopic. A {\em regular 
homotopy} of 
$\S$ consists of a sequence of isotopies and crossings of bands (see 
figure \ref{fA}), so we only have to show that the band crossings can 
be 
achieved by blink surgeries. 

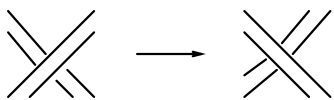
\begin{figure}[htpb]
$$\printname{fA}
	\setlength{\unitlength}{0.03\standardunitlength}
	\begin{array}{c}  \hspace{-1.7mm}
        	\raisebox{-8pt}{\begingroup\makeatletter\ifx\SetFigFont\undefined
% extract first six characters in \fmtname
\def\x#1#2#3#4#5#6#7\relax{\def\x{#1#2#3#4#5#6}}%
\expandafter\x\fmtname xxxxxx\relax \def\y{splain}%
\ifx\x\y   % LaTeX or SliTeX?
\gdef\SetFigFont#1#2#3{%
  \ifnum #1<17\tiny\else \ifnum #1<20\small\else
  \ifnum #1<24\normalsize\else \ifnum #1<29\large\else
  \ifnum #1<34\Large\else \ifnum #1<41\LARGE\else
     \huge\fi\fi\fi\fi\fi\fi
  \csname #3\endcsname}%
\else
\gdef\SetFigFont#1#2#3{\begingroup
  \count@#1\relax \ifnum 25<\count@\count@25\fi
  \def\x{\endgroup\@setsize\SetFigFont{#2pt}}%
  \expandafter\x
    \csname \romannumeral\the\count@ pt\expandafter\endcsname
    \csname @\romannumeral\the\count@ pt\endcsname
  \csname #3\endcsname}%
\fi
\fi\endgroup
\begin{picture}(4524,1239)(0,-10)
\thicklines
\path(12,12)(1212,1212)
\path(312,12)(1212,912)
\path(1212,12)(837,387)
\path(537,612)(12,1212)
\path(12,912)(387,462)
\path(687,237)(912,12)
\path(3312,1212)(4512,12)
\path(3312,912)(4212,12)
\path(4512,1212)(3987,612)
\path(3762,387)(3312,12)
\path(4212,1212)(3837,762)
\path(3612,537)(3312,312)
\path(1812,612)(2712,612)
\path(2592.000,582.000)(2712.000,612.000)(2592.000,642.000)
\end{picture} }
        	\hspace{-1.9mm}
	\end{array}
 $$
\caption{An illustration of  a crossing of bands.}\lbl{fA}
\end{figure}

Let $b_1 ,b_2$ be any two bands of 
$\S$, possibly the same band. Each time we encounter a crossing of 
$b_1$ with $b_2$ there is a corresponding change in the Seifert 
matrix. Since the net change in the Seifert matrix must be zero, we 
conclude that there are an equal number of crossings of $b_1$ and 
$b_2$ in each of the two directions. Now suppose that our regular 
homotopy actually breaks up into a sequence of isotopies and {\em 
double band crossings}, where we define a double band crossing as 
two simultaneous crossings of a single pair of bands in opposite 
directions (see figure \ref{fB}). 
We illustrate in figure \ref{fC} that a double band 
crossing can be achieved by a blink surgery. 

\begin{figure}[htpb]
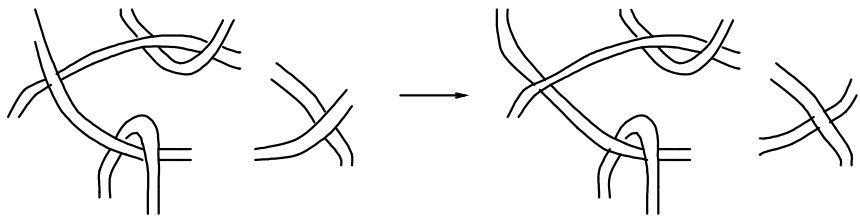

$$\printname{fB}
	\setlength{\unitlength}{0.03\standardunitlength}
	\begin{array}{c}  \hspace{-1.7mm}
        	\raisebox{-8pt}{\input draws/fB.tex }
        	\hspace{-1.9mm}
	\end{array}
 $$
\caption{A double crossing change of bands in opposite 
directions.}\lbl{fB}
\end{figure}

\begin{figure}[htpb]
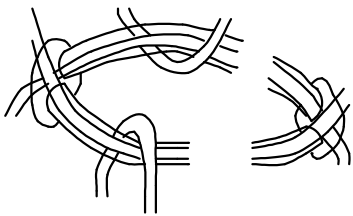

$$\printname{fC}
	\setlength{\unitlength}{0.03\standardunitlength}
	\begin{array}{c}  \hspace{-1.7mm}
        	\raisebox{-8pt}{\input draws/fC.tex }
        	\hspace{-1.9mm}
	\end{array}
 $$
\caption{A $1$-pair blink that achieves the double crossing change
of the bands of figure \ref{fB}.}\lbl{fC}
\end{figure}

We need to do two 
oppositely framed surgeries on circles around the two bands at the 
crossing points. The surface bounded by these two circles is 
obtained by taking the small twice punctured disks bounded by these 
circles and connecting the punctures of one with the punctures of 
the other by tubes along the two band segments connecting the 
crossing points.

So it suffices to show that we can find a regular homotopy from 
$\S$ to $\S'$ consisting of isotopies and double band crossings. Let 
us consider $\S$ and $\S'$ as disks $D$ and $D'$ with bands 
attached. By a preliminary isotopy we can assume that $D=D'$ and 
even a bit more, that $\S$ and $\S'$ coincide in a neighborhood of 
$D=D'$. Then we can choose a regular neighborhood $N$ of $D$ and 
assume that we have a regular homotopy which moves the bands of 
$\S$ onto those of $\S'$ in the complement of $N$ and is stationary 
inside $N$. Let us modify this regular homotopy of $\S$ by 
performing some additional band moves inside $N$. Every time a band 
crossing occurs (outside $N$) let us introduce a crossing of the same 
bands, but in the opposite direction, inside $N$ (see figures \ref{fD}
 and \ref{fE}). 
Thus every band crossing in the original homotopy is replaced by a 
double band crossing and so our new regular homotopy will be a 
sequence of double band crossings. This new homotopy now consists 
of two independent parts: the {\em original homotopy} outside $N$ 
and the 
new part inside $N$. We have complete freedom in how we perform 
each of the band crossings inside $N$ and so, since the number of 
crossings of each pair of bands in the two directions is equal, we 
can choose the corresponding crossings introduced inside $N$ to 
cancel each other. Thus the net effect will be to leave $\S\cap N$ 
unchanged. In other words our {\em modified regular homotopy} will 
have 
the same result as the original regular homotopy, i.e., to move $\S$ 
onto $\S'$. Since the modified regular homotopy is a sequence of 
double band crossings, this completes the proof.
\qed

\begin{figure}[htpb]
$$\printname{fD}
	\setlength{\unitlength}{0.03\standardunitlength}
	\begin{array}{c}  \hspace{-1.7mm}
        	\raisebox{-8pt}{\begingroup\makeatletter\ifx\SetFigFont\undefined
% extract first six characters in \fmtname
\def\x#1#2#3#4#5#6#7\relax{\def\x{#1#2#3#4#5#6}}%
\expandafter\x\fmtname xxxxxx\relax \def\y{splain}%
\ifx\x\y   % LaTeX or SliTeX?
\gdef\SetFigFont#1#2#3{%
  \ifnum #1<17\tiny\else \ifnum #1<20\small\else
  \ifnum #1<24\normalsize\else \ifnum #1<29\large\else
  \ifnum #1<34\Large\else \ifnum #1<41\LARGE\else
     \huge\fi\fi\fi\fi\fi\fi
  \csname #3\endcsname}%
\else
\gdef\SetFigFont#1#2#3{\begingroup
  \count@#1\relax \ifnum 25<\count@\count@25\fi
  \def\x{\endgroup\@setsize\SetFigFont{#2pt}}%
  \expandafter\x
    \csname \romannumeral\the\count@ pt\expandafter\endcsname
    \csname @\romannumeral\the\count@ pt\endcsname
  \csname #3\endcsname}%
\fi
\fi\endgroup
\begin{picture}(9013,1838)(0,-10)
\thicklines
\path(7726,1811)(7876,1361)
\path(7576,1736)(7726,1361)
\path(5776,161)(6076,161)
\path(6226,161)(7426,161)
\path(7576,161)(7951,161)
\path(8101,161)(9001,161)
\path(7875,1362)(7950,1212)
\path(7725,1362)(7800,1212)
\path(8250,1662)(8400,1662)
\path(8250,1512)(8400,1512)
\path(4125,762)(5025,762)
\path(4905.000,732.000)(5025.000,762.000)(4905.000,792.000)
\path(2175,1662)(2325,1212)
\path(2025,1587)(2175,1212)
\path(225,12)(525,12)
\path(675,12)(1875,12)
\path(2025,12)(2400,12)
\path(2550,12)(3450,12)
\path(7951,1211)	(7968.870,1164.828)
	(7985.352,1121.832)
	(8000.499,1081.846)
	(8014.367,1044.705)
	(8038.485,978.303)
	(8058.145,921.305)
	(8073.787,872.394)
	(8085.851,830.252)
	(8101.000,761.000)

\path(8101,761)	(8105.051,723.427)
	(8106.953,677.738)
	(8107.201,626.663)
	(8106.291,572.930)
	(8104.720,519.269)
	(8102.984,468.407)
	(8101.579,423.075)
	(8101.000,386.000)

\path(8101,386)	(8101.000,352.952)
	(8101.000,307.951)
	(8101.000,245.725)
	(8101.000,206.505)
	(8101.000,161.000)

\path(7801,1211)	(7818.897,1156.834)
	(7835.401,1106.413)
	(7850.567,1059.544)
	(7864.451,1016.034)
	(7877.107,975.692)
	(7888.590,938.325)
	(7908.257,871.749)
	(7923.893,814.766)
	(7935.935,765.839)
	(7944.825,723.430)
	(7951.000,686.000)

\path(7951,686)	(7954.071,636.136)
	(7953.730,573.241)
	(7952.024,510.476)
	(7951.000,461.000)

\path(7951,461)	(7951.000,416.932)
	(7951.000,356.930)
	(7951.000,318.757)
	(7951.000,273.963)
	(7951.000,221.671)
	(7951.000,161.000)

\path(7501,761)	(7518.217,701.028)
	(7532.754,649.235)
	(7544.830,604.741)
	(7554.666,566.668)
	(7568.496,506.266)
	(7576.000,461.000)

\path(7576,461)	(7579.589,416.634)
	(7580.785,356.533)
	(7580.486,318.384)
	(7579.589,273.665)
	(7578.093,221.497)
	(7576.000,161.000)

\path(7351,761)	(7359.968,722.005)
	(7368.236,685.725)
	(7382.787,620.761)
	(7394.872,565.010)
	(7404.711,517.374)
	(7412.524,476.753)
	(7418.530,442.050)
	(7426.000,386.000)

\path(7426,386)	(7428.154,352.800)
	(7428.872,307.749)
	(7428.154,245.573)
	(7427.257,206.416)
	(7426.000,161.000)

\path(6075,162)	(6079.928,201.437)
	(6084.696,238.099)
	(6093.861,303.651)
	(6102.715,359.755)
	(6111.476,407.509)
	(6129.603,482.362)
	(6150.000,537.000)

\path(6150,537)	(6173.482,580.505)
	(6205.417,629.910)
	(6243.581,682.748)
	(6285.750,736.549)
	(6329.701,788.843)
	(6373.208,837.163)
	(6414.050,879.038)
	(6450.000,912.000)

\path(6450,912)	(6512.839,959.006)
	(6550.882,984.026)
	(6592.455,1009.742)
	(6636.894,1035.921)
	(6683.535,1062.325)
	(6731.715,1088.721)
	(6780.769,1114.871)
	(6830.033,1140.541)
	(6878.844,1165.496)
	(6926.537,1189.498)
	(6972.449,1212.315)
	(7015.915,1233.708)
	(7056.271,1253.444)
	(7125.000,1287.000)

\path(7125,1287)	(7160.251,1303.892)
	(7201.104,1322.675)
	(7249.097,1344.008)
	(7305.769,1368.551)
	(7372.657,1396.963)
	(7410.413,1412.825)
	(7451.299,1429.902)
	(7495.509,1448.275)
	(7543.234,1468.028)
	(7594.667,1489.242)
	(7650.000,1512.000)

\path(6225,162)	(6236.566,223.087)
	(6247.067,275.678)
	(6256.722,320.650)
	(6265.751,358.883)
	(6282.810,418.646)
	(6300.000,462.000)

\path(6300,462)	(6342.356,534.429)
	(6371.167,576.423)
	(6402.847,619.530)
	(6435.735,661.696)
	(6468.165,700.866)
	(6525.000,762.000)

\path(6525,762)	(6566.395,797.658)
	(6618.849,837.661)
	(6678.955,880.099)
	(6743.306,923.062)
	(6808.494,964.644)
	(6871.110,1002.933)
	(6927.748,1036.021)
	(6975.000,1062.000)

\path(6975,1062)	(7011.293,1079.501)
	(7056.381,1098.938)
	(7107.377,1119.519)
	(7161.398,1140.457)
	(7215.555,1160.964)
	(7266.966,1180.249)
	(7312.742,1197.524)
	(7350.000,1212.000)

\path(7350,1212)	(7405.082,1234.034)
	(7439.544,1247.819)
	(7480.084,1264.035)
	(7527.801,1283.122)
	(7583.793,1305.518)
	(7649.160,1331.665)
	(7685.702,1346.281)
	(7725.000,1362.000)

\path(7800,1587)	(7844.583,1604.177)
	(7883.145,1618.686)
	(7944.844,1640.576)
	(7990.371,1654.428)
	(8025.000,1662.000)

\path(8025,1662)	(8058.436,1665.572)
	(8103.566,1666.762)
	(8165.663,1665.572)
	(8204.722,1664.084)
	(8250.000,1662.000)

\path(7875,1437)	(7929.198,1468.506)
	(7969.920,1490.336)
	(8025.000,1512.000)

\path(8025,1512)	(8058.830,1517.288)
	(8104.091,1519.050)
	(8166.057,1517.288)
	(8204.952,1515.084)
	(8250.000,1512.000)

\path(525,12)	(526.762,58.936)
	(528.787,102.596)
	(531.100,143.142)
	(533.730,180.742)
	(540.051,247.756)
	(547.969,304.957)
	(557.703,353.664)
	(569.473,395.195)
	(600.000,462.000)

\path(600,462)	(649.274,527.158)
	(715.730,593.072)
	(753.864,625.797)
	(794.451,658.092)
	(836.875,689.750)
	(880.523,720.566)
	(924.778,750.333)
	(969.029,778.845)
	(1012.659,805.895)
	(1055.054,831.277)
	(1095.600,854.785)
	(1133.683,876.212)
	(1200.000,912.000)

\path(1200,912)	(1265.055,942.225)
	(1303.782,957.987)
	(1345.794,974.003)
	(1390.458,990.138)
	(1437.136,1006.258)
	(1485.194,1022.228)
	(1533.997,1037.914)
	(1582.910,1053.180)
	(1631.297,1067.892)
	(1678.523,1081.917)
	(1723.952,1095.118)
	(1766.950,1107.361)
	(1806.880,1118.512)
	(1875.000,1137.000)

\path(1875,1137)	(1948.607,1155.633)
	(1992.007,1166.038)
	(2038.895,1176.950)
	(2088.586,1188.201)
	(2140.400,1199.626)
	(2193.652,1211.059)
	(2247.660,1222.335)
	(2301.741,1233.288)
	(2355.213,1243.753)
	(2407.392,1253.563)
	(2457.595,1262.553)
	(2505.141,1270.557)
	(2549.345,1277.410)
	(2589.526,1282.947)
	(2625.000,1287.000)

\path(2625,1287)	(2680.330,1290.594)
	(2714.853,1291.493)
	(2755.414,1291.793)
	(2803.110,1291.493)
	(2859.041,1290.594)
	(2924.305,1289.097)
	(2960.780,1288.123)
	(3000.000,1287.000)

\path(675,12)	(679.130,51.384)
	(683.206,88.001)
	(691.307,153.482)
	(699.522,209.544)
	(708.071,257.284)
	(727.050,332.193)
	(750.000,387.000)

\path(750,387)	(792.152,446.023)
	(852.844,510.248)
	(918.363,569.098)
	(975.000,612.000)

\path(975,612)	(1010.707,632.406)
	(1055.434,653.405)
	(1106.267,674.465)
	(1160.291,695.051)
	(1214.591,714.631)
	(1266.253,732.672)
	(1312.360,748.639)
	(1350.000,762.000)

\path(1350,762)	(1394.293,777.727)
	(1448.747,796.330)
	(1509.989,816.803)
	(1574.648,838.140)
	(1639.350,859.334)
	(1700.724,879.380)
	(1755.398,897.270)
	(1800.000,912.000)

\path(1800,912)	(1844.464,927.113)
	(1898.803,945.935)
	(1959.764,967.131)
	(2024.093,989.366)
	(2088.534,1011.306)
	(2149.836,1031.614)
	(2204.742,1048.957)
	(2250.000,1062.000)

\path(2250,1062)	(2286.802,1071.318)
	(2331.945,1081.976)
	(2382.696,1093.314)
	(2436.326,1104.668)
	(2490.103,1115.375)
	(2541.295,1124.775)
	(2587.171,1132.204)
	(2625.000,1137.000)

\path(2625,1137)	(2680.330,1140.594)
	(2714.853,1141.493)
	(2755.414,1141.793)
	(2803.110,1141.493)
	(2859.041,1140.594)
	(2924.305,1139.097)
	(2960.780,1138.123)
	(3000.000,1137.000)

\path(2400,1062)	(2417.870,1015.828)
	(2434.352,972.832)
	(2449.499,932.846)
	(2463.367,895.705)
	(2487.485,829.303)
	(2507.145,772.305)
	(2522.787,723.394)
	(2534.851,681.252)
	(2550.000,612.000)

\path(2550,612)	(2554.051,574.427)
	(2555.953,528.738)
	(2556.201,477.663)
	(2555.291,423.930)
	(2553.720,370.269)
	(2551.984,319.407)
	(2550.579,274.075)
	(2550.000,237.000)

\path(2550,237)	(2550.000,203.952)
	(2550.000,158.951)
	(2550.000,96.725)
	(2550.000,57.505)
	(2550.000,12.000)

\path(2250,1062)	(2267.897,1007.834)
	(2284.401,957.413)
	(2299.567,910.544)
	(2313.451,867.034)
	(2326.107,826.692)
	(2337.590,789.325)
	(2357.258,722.749)
	(2372.893,665.766)
	(2384.935,616.839)
	(2393.825,574.430)
	(2400.000,537.000)

\path(2400,537)	(2403.071,487.136)
	(2402.730,424.241)
	(2401.024,361.476)
	(2400.000,312.000)

\path(2400,312)	(2400.000,267.932)
	(2400.000,207.930)
	(2400.000,169.757)
	(2400.000,124.963)
	(2400.000,72.671)
	(2400.000,12.000)

\path(1950,612)	(1967.217,552.028)
	(1981.754,500.235)
	(1993.830,455.741)
	(2003.666,417.668)
	(2017.496,357.266)
	(2025.000,312.000)

\path(2025,312)	(2028.589,267.634)
	(2029.785,207.533)
	(2029.486,169.384)
	(2028.589,124.665)
	(2027.093,72.497)
	(2025.000,12.000)

\path(1800,612)	(1808.968,573.005)
	(1817.236,536.725)
	(1831.787,471.761)
	(1843.872,416.010)
	(1853.711,368.374)
	(1861.524,327.753)
	(1867.530,293.050)
	(1875.000,237.000)

\path(1875,237)	(1877.154,203.800)
	(1877.872,158.749)
	(1877.154,96.573)
	(1876.257,57.416)
	(1875.000,12.000)

\put(5551,536){\makebox(0,0)[lb]{$N$}}
\put(0,387){\makebox(0,0)[lb]{$N$}}
\end{picture} }
        	\hspace{-1.9mm}
	\end{array}
 $$
\caption{An illustration of an original homotopy.}\lbl{fD}
\end{figure}
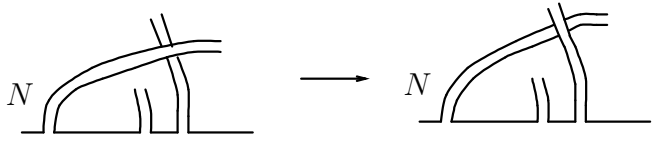

\begin{figure}[htpb]
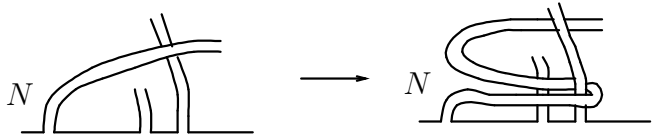

$$\printname{fE}
	\setlength{\unitlength}{0.03\standardunitlength}
	\begin{array}{c}  \hspace{-1.7mm}
        	\raisebox{-8pt}{\input draws/fE.tex }
        	\hspace{-1.9mm}
	\end{array}
 $$
\caption{An illustration of a modified  homotopy.}\lbl{fE}
\end{figure}

\section{Appendix}

\subsection{Remarks on the group $\Lg$}

If $L\sub H=H_1 (\Sigma_g )$ is the Lagrangian used to define $\Lg$, 
then we can also  define a larger group $\blg=\{ h: 
h_{\ast}|L=\text{identity}\}$, where $h_{\ast}$ is the automorphism 
of $H$ induced by $h$. Clearly $\Lg\sub\blg$ and 
$\Tg\sub\blg$. So we have a lattice of subgroups of the 
mapping class group:
$$\begin{array}{ccccc}
 & & \blg & & \\
 &\nearrow & &\nwarrow & \\
\Lg & & & &\Tg \\
 &\nwarrow & &\nearrow & \\
 & & \Lg\cap\Tg & & \\
 & & \uparrow & & \\
 & & \Kg & & 
\end{array}$$

Note that the above diagram defines a map  
$\Lg /\Lg\cap\Tg\ \to\ \blg/\Tg$. We can now show the 
following:

\begin{proposition}
\lbl{prop.lg}
The above defined map 
$\Lg /\Lg\cap\Tg\ \to\ \blg/\Tg$ is an isomorphism.
\end{proposition}
\begin{pf} 
By its definition  it follows that it is one-to-one. In order to show
that it is onto, recall first \cite{Ha}, \cite{Jo}, \cite{Mo}
the following classical short exact sequence:
\begin{equation}
1 \to \Tg \to \Gamma_{g,1} \to Sp(2g,\BZ) \to 1
\end{equation}
where the map $\Gamma_{g,1} \to Sp(2g,\BZ)$ is the map $h \to 
h_\ast$
that sends a surface diffeomorphism to its linear action on 
$H_1 (\Sigma_g, \BZ )$. We therefore have an isomorphism
$ \Gamma_{g,1}/ \Tg \simeq Sp(2g,\BZ)$. We can therefore 
identify  $\blg/\Tg$ with its image in $Sp(2g,\BZ)$, and as 
such,
$\blg/\Tg$ consists of all isometries  $\phi$ of $H_1 (\S_g )$ 
which are  the identity on $L$. 
With these preliminaries in mind, in order to show that the
map of the proposition is onto,  
it suffices to show that every isometry $\phi$ of $H_1 (\S_g )$ 
which is the identity on $L$ is induced by some product 
$\rho_1\cdots\rho_k$ of $L$-twists. If we write $H_1 (\S_g 
)=L\oplus L'$ where 
$L'$ is a Lagrangian dual to $L$ and we choose a basis $\{ e_i\}$ of 
$L$ and dual basis $\{ e'_i\}$for $L'$, then $\phi$ has a matrix 
representative:
$\pmatrix I&C\\0&I\endpmatrix$, where $C$ is a symmetric matrix. 
If $\rho$ is a Dehn twist along a simple closed curve representing 
$\sum_i\lambda_i e_i$, then it has such a matrix representative, 
where the entries of $C$ are given by 
$c_{ij}=\pm\lambda_i\lambda_j$. We can certainly realize the 
elements $e_i$ and $e_i\pm e_j$ by simple closed curves and it is 
then an easy exercise to see that any $C$ can be realized by a 
composition of Dehn twists along such curves, using the fact that
$$\pmatrix I&C\\0&I\endpmatrix\cdot\pmatrix 
I&C'\\0&I\endpmatrix =\pmatrix I&C+C'\\0&I\endpmatrix$$
\end{pf}

Notice also that $[\blg,\blg]\sub\Tg$ because $\blg
/\Tg$ is abelian.

 The next natural problem to consider is the determination of 
$\Lg\cap\Tg /\Kg$. In order to do so, we will need 
an important  homomorphism $\tau 
:\Tg\to\Lambda^3 H$ defined by D. Johnson (see \cite{Jo2}). 
We review its definition here:
If $h\in\Tg$ then,
by definition, $h_{\ast}$ is the identity on $H$. Thus, for any 
$\beta\in\pi =\pi_1 (\Sigma_g )$, we can write: 
\begin{equation}
h_{\ast}(\beta )\beta^{-1}
\equiv t(h )\cdot\beta\in \pi_2 /\pi_3 \simeq\Lambda^2 H
\end{equation}
This defines a homomorphism $t:\Tg\to\hom (H,\Lambda^2 H)$. We 
now
have the identifications:
\begin{equation}
\hom (H,\Lambda^2 H)\simeq H^{\ast}\otimes\Lambda^2 H\simeq 
H\otimes
\Lambda^2 H
\end{equation}
where the latter isomorphism uses the symplectic structure on $H$.
 Thus we obtain a homomorphism $t':\Tg\to H\otimes
\Lambda^2 H$.
Johnson showed that $\text{im}( t')\sub\Lambda^3 H$, where 
the embedding
$\Lambda^3 H\sub H\otimes
\Lambda^2 H$ is defined by:
\begin{equation}
x\wedge y\wedge z\mapsto x\otimes (y\wedge z)+ y\otimes 
(z\wedge x)+
 z\otimes (x\wedge y)
\end{equation}
This defines a homomorphism $\tau: \Tg \to \Lambda^3 H$.
Johnson showed that $\tau$ is onto and its kernel is exactly $\Kg$.

With these preliminaries in mind, 
an important step in understanding $\Lg$ is the calculation of 
$\tau (\Lg\cap\Tg )$. For example, by Proposition~\ref{prop.lg},  
$\tau (\Lg\cap\Tg )=\Lambda^3 H$ if and only if $\Lg =\blg$.
 But we now show that this is false.

\begin{proposition}
\lbl{prop.lg1}
 Suppose $h\in\Tg\cap\Lg$. Then $\tau (h)
\in\ker\{ \Lambda^3 H\to\Lambda^3 (H/L)\}$.
\end{proposition}

\begin{pf} 
Let $h\in\Tg\cap\Lg$. With the discussion of  Johnson's 
homomorphism
above, and using the the following commutative diagram:
\begin{displaymath}
\begin{CD}
\Lambda^3 H @>>> \Lambda^3 (H/L) \\
@VVV    @VVV  \\
H\otimes\Lambda^2 H  @>>> H/L\otimes\Lambda^2 (H/L)
\end{CD}
\end{displaymath}
(where the vertical arrows are both injections, and $L$ is the
Lagrangian used to define $\Lg$) it follows that
the proposition is equivalent to showing 
that $t'(h)\in\ker\{ H\otimes\Lambda^2 H\to H/L\otimes\Lambda^2 
(H/L)\}$
or, equivalently, $t(h)\in\ker\{ \hom (H,\Lambda^2 H)\to\hom 
(L,\Lambda^2  (H/L)\}$. 
In order to show that we need the following lemma:

\begin{lemma}
\lbl{lem.john}
Suppose that $\Sigma$ is a compact orientable surface with one 
boundary
component and $C$ a simple closed curve in the interior of $\Sigma$
representing, up to conjugacy, an element $\alpha\in\pi =\pi_1 
(\Sigma )$.
Let $h$ denote the homeomorphism of $\Sigma$ defined by a Dehn 
twist
along $C$. Then, for any $\beta\in\pi$, we can write
$$h_{\ast}(\beta )\beta^{-
1}=\gamma_1^{\epsilon_1}\cdots\gamma_k^{\epsilon_k}
$$ 
where each $\gamma_i$ is a conjugate of $\alpha$ and, if $[\xi ]$ 
denotes the homology class of $\xi$ for any $\xi\in\pi$, then
$\sum_i\epsilon_i =\pm [\alpha ]\cdot [\beta ]$, the intersection 
number,
 and the sign depends on the direction of the Dehn twist.
\end{lemma} 

\begin{pf} If $\lambda$ is any path in $\Sigma$ which
intersects $C$ transversely, then we can write $h\circ\lambda$ as a
 product $\lambda_1\cdot C^{\epsilon_1}\cdots\lambda_k\cdot 
C^{\epsilon_k}
\cdot\lambda_{k+1}$, for some factorization $\lambda 
=\lambda_1\cdots
\lambda_{k+1}$ as a product of paths. The $C^{\epsilon_i}$ insert 
themselves whenever $\lambda$
 crosses $C$ and $\epsilon_i$ is the sign of the intersection.
\end{pf}

We can rewrite $h_{\ast}(\beta )\beta^{-1}$ as given in 
lemma~\ref{lem.john}
in the form 
\begin{equation}
\lbl{eq.comm}
[\xi_1,\alpha^{\epsilon_1}]\alpha^{\epsilon_1}\cdots 
[\xi_k,\alpha^{\epsilon_k}]\alpha^{\epsilon_k}\equiv
(\prod_i[\xi_i,\alpha^{\epsilon_i}])\alpha^e \mod\pi_3
\end{equation}
where $e=\pm [\alpha ]\cdot [\beta ]$. In particular hold in mind the
case of $e=0$. In this case if we apply another homeomorphism 
defined by
 a Dehn twist along a curve representing $\gamma\in\pi$ so that 
$[\alpha ]\cdot [\gamma ]=0$, then $h_{\ast}(\beta )\beta^{-1}$ is
mapped to a new element which is still in the form of 
equation \eqref{eq.comm}
 with $e=0$.
Continuing in this way we obtain the following conclusion: 
suppose that $h$ is a product of Dehn twists on curves representing
elements in some Lagrangian $L$. If $\beta\in\pi$ also represents
an element of $L$, then
\begin{equation}
\lbl{eq.lag}
h_{\ast}(\beta )\beta^{-1}\equiv\prod_i [\xi_i,\alpha_i ] \mod\pi_3
\end{equation}
where $[\alpha_i ]\in L$. But this is shows that:
$$  t (h)\cdot\beta =\sum_i [\xi_i ]\wedge [\alpha_i ] $$
Since $[\alpha_i ]\in L$, the right side clearly maps to $0$ in
$\Lambda^2 (H/L)$.
\end{pf} 

Proposition \ref{prop.lg} immediately implies the following 
corollary:

\begin{corollary}
\lbl{cor.lg}
As a subgroup of the mapping class group, $\blg$ is generated 
by
$\Tg$ and $\Lg$.
\end{corollary} 

\begin{remark} 
\lbl{rem.KL}
$ t| \Tg \cap \Lg $ actually extends to a homomorphism
 $\tilde t:\blg\to\hom (L,\Lambda^2 H)$ by the same definition
as for $t$, using the defining property that, if $h\in\blg$ then
$h_{\ast}|L=$identity.  The above proof actually shows that
$\tilde t (\Lg )
\sub \hom ( L,K )$, where $K=\ker\{\Lambda^2 H\to\Lambda^2 
(H/L)\}$.
\end{remark}

\begin{question}
Is $\Lg={\tilde t}^{-1}\hom ( L,K )?$ Is $\tau (\Lg\cap\Tg )=
\ker\{\Lambda^3 H\to\Lambda^3 (H/L)\}$?
\end{question}

\subsection{The lower central series of $\blg $}

In this section we study the  
image of the lower central series of $\blg $ under Johnson's map 
into 
$\Lambda^3 H$. In particular we prove:

\begin{proposition}
\lbl{prop.lcblg}
For all $ g \geq 1$ we have:
 $(\blg )_5\sub\Kg$. In addition, for
$g \geq 1$ we have:  $(\blg )_4\not\sub\Kg$ if $g\ge 3$.
\end{proposition}

The proof of proposition~\ref{prop.lcblg} will be based on the 
following:

\begin{lemma}
\lbl{lem.lcbg}
With the notation of remark \ref{rem.KL}, we have the following:
$t'([\blg ,\blg ])\sub (L\otimes\Lambda^2 H) +  (H\otimes K)$.
\end{lemma}

%\begin{pf}[Proof of lemma~\ref{lem.lcbg}]
Using corollary \ref{cor.lg}, and the fact that
 $[\Tg ,\Tg ]\sub\Kg$, the lemma  will follow from the 
following two assertions:
\begin{enumerate}
\item $t' ([\Lg ,\Lg ])\sub H\otimes K$.
\item $t' ([\Lg ,\Tg ])\sub (H\otimes K) +  (L\otimes\Lambda^2 
H)$.
\end{enumerate}

 \begin{pf}{Proof of (1)} $[\Lg ,\Lg ]$ is generated by elements of 
the 
form
 $[h,D_C ]$, where $h\in\Lg$ and $C$ is a simple closed curve in $M$ 
which
 represents an element of $L$ and $D_C$ denotes a Dehn twist along 
$C$. 
 Now $[h, D_C ]=D_{h(C)}\circ (D_C )^{-1}$. We can apply 
lemma~\ref{lem.john} to obtain;
\begin{equation}
(D_C )_{\ast}(\beta )\equiv\beta\prod_i [\xi_i 
,\alpha^{\epsilon_i}]\alpha^d\ \mod\pi_3
\end{equation}
where $C$ represents $\alpha\in\pi$, up to conjugation, and $d=\pm 
[\alpha ]\cdot [\beta ]$,
 for any $\beta\in\pi$. Similarly we have:
\begin{equation}
(D_{h(C)} )_{\ast}(\beta )\equiv\beta\prod_i [\theta_i 
,h_{\ast}(\alpha 
)^{\epsilon'_i}]h_{\ast }
(\alpha )^e\ \mod\pi_3
\end{equation}
where $e=\pm [h_{\ast}\alpha ]\cdot [\beta ]$, for any $\beta\in\pi$.
Putting these together we get:
\begin{equation}
\lbl{eq.dc} 
(D_{h(C)} )_{\ast}\circ (D_C )_{\ast}^{-1}(\beta )\equiv\beta\prod_i 
[\theta_i ,
h_{\ast}(\alpha )^{\epsilon'_i}]\prod_i 
[\xi_i ,\alpha^{\epsilon_i}]\alpha^d h_{\ast }(\alpha )^e\
 \mod\pi_3
\end{equation}
where $d=\pm [\alpha ]\cdot [\beta ]$ again, but now $e=\pm 
[h_{\ast}\alpha ]\cdot
 (D_C )_{\ast}^{-1}[\beta ]$.

Since $h\in\Lg$ and $[\alpha ]\in L$, we have $h_{\ast}(\alpha 
)\alpha^{-1}\in\pi_2$ and
 so $\alpha$ and $h_{\ast}(\alpha )$ commute mod $\pi_3$. Thus  
equation \eqref{eq.dc} can
 be rewritten:
\begin{equation}
\lbl{eq.dc1} 
(D_{h(C)} )_{\ast}\circ (D_C )_{\ast}^{-1}(\beta )\equiv\beta\prod_i 
[\theta_i 
,
h_{\ast}(\alpha )^{\epsilon'_i}]\prod_i [\xi_i ,\alpha^{\epsilon_i}]
(h_{\ast }(\alpha )\alpha^{-1})^e\alpha^{d+e}\ \mod\pi_3
\end{equation}
This simplifies considerably to:
$$
(D_{h(C)} )_{\ast}\circ (D_C )_{\ast}^{-1}(\beta 
)\equiv\beta\alpha^{d+e}\ 
\mod\pi_2
$$
But $D_{h(C)}\circ (D_C )^{-1}\in\Tg$, since $h(C)$ is 
homologous to $C$, and so $d+e=0$. (We can assume that $[\alpha 
]\not= 0$ otherwise we already have that $D_C\in\Kg$.) 

We can apply equation \eqref{eq.lag} to write:
$$h_{\ast}(\alpha )\alpha^{-1}\equiv\prod_i [\eta_i ,\alpha_i 
]\mod\pi_3
$$
where $[\alpha_i ]\in L$. Putting this all together into 
equation \eqref{eq.dc1} we get:
\begin{equation}
\lbl{eq.dc2}
(D_{h(C)} )_{\ast}\circ (D_C )_{\ast}^{-1}(\beta )\equiv\beta\prod_i 
[\zeta_i 
,\alpha'_i ]
\end{equation}
where $[\alpha'_i ]\in L$. This translates into:
$$\tilde t ((D_{h(C)} )\circ (D_C )^{-1})\cdot [\beta ]=\sum_i 
[\zeta_i ]\wedge [\alpha'_i ]
$$
and this element lies in $K$.  
\end{pf}
\begin{pf}{Proof of (2)}
Suppose $\l\in\blg$ and $h\in\Tg$. Then $[\l ,h]=(\l h\l^{-1})h^{-1}$ 
which can be written, as an element of the abelianization $\Tg 
/(\Tg )_2$, in the additive form $\l\cdot h-h$, where we use the 
canonical 
action of the mapping class group $\Gamma_{g,1}$ on $\Tg /(\Tg 
)_2$ 
by conjugation. It is pointed out by Johnson (see e.g. \cite{Jo2}) that 
$\tau$, or $t':\Tg\to H\otimes\L^2 H$, is equivariant with respect 
to the action of $\Gamma_{g,1}$ (acting on the right side by the 
canonical action on $H$). Thus 
$$ t'([\l ,h])=t'(\l\cdot h-h)=(\l -1)\cdot t'(h) $$
If $a\otimes (a_1\wedge a_2 )\in H\otimes\L^2 H$, then 
\begin{eqnarray}
\lbl{eq.lmo}
(\l -1)\cdot (a\otimes (a_1\wedge a_2 ))&=&(\l -1)a\otimes\l  
(a_1\wedge a_2 )+\ a\otimes ((\l -1)a_1\wedge\l a_2 )\\
 & & +\ a\otimes 
(a_1\wedge (\l -1)a_2 )\nonumber
\end{eqnarray}
Now recall that, for any $\l\in\blg$, the action on $H$ satisfies:
\begin{itemize}
\item $\l |L=\text{identity}$
\item $(\l -1)(H)\sub L$
\end{itemize}
and so $(\l -1)^2 =0$.
Thus, in equation~\eqref{eq.lmo}, the terms on the right side are in 
either $L\otimes\L^2 H$ or $H\otimes K$.

This completes the proof of (2) and of lemma~\ref{lem.lcbg}.
\end{pf}

\begin{pf}{Proof of Proposition~\ref{prop.lcblg}}

We will use an argument similar to that in the proof of (2) above to 
prove the following assertions in order.
\begin{enumerate}
\item $t'((\blg )_3 )\sub (L\otimes K) +  (H\otimes\L^2 L)$
\item  $t'((\blg )_4 )\sub (L\otimes\L^2 L)$
\item $t'((\blg )_5 )=0$
\end{enumerate}
Recalling corollary \ref{cor.lg}, the above two assertions prove
proposition \ref{prop.lcblg}.

To prove (1) we apply equation~\eqref{eq.lmo}, where we can 
assume, 
by 
lemma~\ref{lem.lcbg}, either $a_1\wedge a_2\in K$ or $a\in L$ and 
we see that the terms on the right side of equation~\eqref{eq.lmo} 
lie 
in 
$(L\otimes K) + (H\otimes\L^2 L)$. Note that $(\l -1)K\sub\L^2 
L$.

To prove (2) (or (3)) we use equation~\eqref{eq.lmo} in the same 
way, 
taking into account (1) (or (2)) to tell where $a\otimes 
(a_1\wedge a_2 )$ must lie.

This completes the proof that $(\tilde{\Lg} )_5\sub\Kg$.
It remains to show that, if $g\ge 3$ then $t'((\blg )_4 )\not= 0$.

Suppose $h\in\Tg$ so that $\tau (h)=a_1\wedge a_2\wedge 
a_3\in\L^3 H$. If $\l\in\Lg$ then we have $\tau ([\l ,h])=(\l -1)\cdot 
(a_1\wedge a_2 \wedge a_3)$. But we can use the following 
analogue of equation~\eqref{eq.lmo}:
\begin{eqnarray}
\lbl{eq.lmo1}
(\l -1)\cdot (a_1\wedge a_2 \wedge a_3)&=&(\l -1)a_1\wedge\l  
a_2\wedge\l a_3+\ a_1\wedge (\l -1)a_2\wedge\l a_3 \\
 & & +\ a_1\wedge 
a_2\wedge (\l -1)a_3 \nonumber
\end{eqnarray}
Noting that $(\l -1)a_i\in L$ and $(\l -1)|L=0$, we can use 
equation~\eqref{eq.lmo1} in this way repeatedly to compute:
$$\tau ([\l ,[\l ,[\l ,h]]])= 6(\l -1)a_1\wedge (\l -1)a_2\wedge (\l 
-1)a_3
$$
Now suppose, following the conventions in the proof of 
proposition~\ref{prop.lg}, that $\{ e_i\}$ is a basis of $L$ and $\{ 
e'_i\}$ is a dual basis of $L'$. For any symmetric matrix $C$ there is 
some $\l\in\Lg$ so that:
$$\l (e_i )=e_i ,\quad \l (e'_i )=e'_i +\sum_j c_{ij}e_j $$
Let us choose $\l$ so that $C$ is the identity matrix and $a_i =e'_i$. 
Then we have 
$$\tau ([\l ,[\l ,[\l ,h]]])=6e'_1\wedge e'_2\wedge e'_3 \not= 0 $$.

This completes the proof of proposition~\ref{prop.lcblg}. 

\end{pf}

\begin{remark} 
It is not clear whether $(\Lg )_4\sub\Kg$.
\end{remark}

%%%%%%%%%%%%%%%%%%%%\inplude{refs}

\ifx\undefined\bysame
	\newcommand{\bysame}{\leavevmode\hbox 
to3em{\hrulefill}\,}
\fi


\begin{thebibliography}{[EMSS]}

\bibitem[B-N]{B-N1} D. Bar-Natan, 
        {\em On the Vassiliev knot invariants}, Topology {\bf 34} 
        (1995)        423-472.
   
\bibitem[BGRT]{BGRT} D. Bar-Natan, S. Garoufalidis, L. Rozansky, D.
        Thurston, to appear.

\bibitem[Ga1]{Ga} S. Garoufalidis,
        {\em On finite type 3-manifold invariants I}, 
        J. Knot Theory and its Ramifications {\bf 5}, 
        no. 4 (1996), p.441-462.

\bibitem[Ga2]{Ga2} \bysame,
        {\em Comparing finite type invariants of knots and integral
        homology 3-spheres}, preprint February 1996.

\bibitem[GL1]{GL1} S. Garoufalidis, J. Levine,
       {\em On finite type 3-manifold invariants II},  
       Math. Annalen, {\bf 306} (1996) p.691-718.

\bibitem[GL2]{GL2} \bysame,
       {\em On finite type 3-manifold invariants IV: comparison of 
       definitions}, 
       Math. Proc. Cambridge Phil. Society, in print. 

\bibitem[GL3]{GL3} \bysame,
       {\em Finite type 3-manifold invariants and the Torelli group I},
       Brandeis Univ. and Harvard Univ. preprint August 1996.

\bibitem[GO1]{GO1} S. Garoufalidis, T. Ohtsuki,
       {\em On finite type 3-manifold invariants III: manifold weight 
       systems}, 
       Tokyo Institute of Technology  and
       M.I.T. preprint August 1995.

\bibitem[GO2]{GO2} \bysame,
       {\em On finite type 3-manifold invariants V: rational homology 
       3-spheres }, 
       Proceedings of the Aarhus Conference, Geometry and Physics, 
         Marcel Dekker (1996) 445-457.

\bibitem[Ha]{Ha} R. Hain,
       {\em Completions of the mapping class groups and the cycle
       $C - C^{-}$},
       Contemporary Math. {\bf 150} (1993) 75-105.

\bibitem[Ha]{Ha2} R. Hain,
       {\em Infinitesimal presentations of the Torelli groups},
       preprint December 1995.

\bibitem[Jo1]{Jo} D. Johnson, 
       {\em Homeomorphisms of a surface which act trivially on 
homology},
       Proc. Amer. Math. Soc., {\bf 75} (1979) 118-125.

\bibitem[Jo2]{Jo2}\bysame ,
       {\em A survey of the Torelli group},
       Contemporary Math. {\bf 20} (1983) 163-179.

%\bibitem[Ha]{Ha} N. Habegger,
%        {\em Finite type 3-manifold invariants: a proof of a conjecture 
%        of Garoufalidis}, preprint, July 1995.

\bibitem[Ho]{Ho} J. Hoste, {\em A formula for Casson's invariant},
        Trans. A.M.S. {\bf 297} (1986) 547-562.

\bibitem[Kh]{Ko} T. Kohno,
       {\em Holonomy Lie algebras, logarithmic connections and the 
lower central series of fundamental groups}, Contemporary Math. 
{\bf 
       90}       (1990) 171-181.
 
\bibitem[Ko]{Kstudent}K. H. Ko, {\em  Seifert matrices and boundary 
link
        cobordisms}, Transactions A.M.S. {\bf 299} (1987) 657-681.

\bibitem[LMO]{LMO} T.T.Q.  Le, J. Murakami, T. Ohtsuki,
   {\em A universal quantum invariant of 3-manifolds},
   preprint, November 1995.

\bibitem[L]{L} T.T.Q.  Le, 
        {\em An invariant of \ihs s which is universal for all \fti s},
        preprint January 1996.

\bibitem[Le]{Le} J. Levine,
  {\em Surgery equivalence of links}, Topology {\bf 26} (1987) 45-
61.

\bibitem[Mo1]{Mo} S. Morita, {\em Casson's invariant for homology 3-
spheres
         and characteristic classes of vector bundles I},
         Topology, {\bf 28} (1989) 305-323.

\bibitem[Mo2]{Mo1} \bysame , {\em On the structure of the Torelli
         group and the Casson invariant},
         Topology, {\bf 30} (1991) 603-621.

%\bibitem[Oh1]{Oh1} T. Ohtsuki,
%        {\em A polynomial invariant of integral homology spheres}, 
%Math. Proc. Camb. Phil. Soc. {\bf 117} (1995) 83-112.

\bibitem[Oh]{Oh} T. Ohtsuki,
        {\em Finite type invariants of integral homology 3-spheres}, 
        J. Knot Theory and its Rami. {\bf 5} (1996) 101-115. 

\bibitem[Qu]{Qu} D. Quillen,
        {\em On the associated graded ring of a group ring},
        Journal of Algebra {\bf 10} (1968) 411-418.
        
\bibitem[Ro]{Ro1} D. Rolfsen,
        {Knots and links}, Publish or Perish, 1976.

%\bibitem[Rz1]{Rz1} L. Rozansky,
%   {\em The trivial connection contribution to Witten's invariant 
%and
%   finite type invariants of rational homology spheres}, preprint 
%   {\tt q-alg/9503011}.

%\bibitem[Rz2]{Rz2} \bysame,
%       {\em Witten's invariants of rational homology spheres at 
%prime values
%       of $K$ and the trivial connection contribution},
%       preprint {\tt q-alg/9504015}.

\bibitem[Wi]{Wi} E. Witten, 
        {\em Quantum field theory and the Jones polynomial},
        Commun.  Math. Phys.  {\bf 121} (1989) 360-376. 

\bibitem[WP]{W} E. Witten, J. Polchinski,
        {\em Evidence for heterotic-type I string duality},
        {\tt hepth 9510169}. 

\end{thebibliography}
\end{document}